\documentclass{amsart}
\usepackage{bm}
\usepackage{amssymb}
\usepackage{amsthm}
\usepackage{float}
\usepackage{tikz}
\usepackage{tikz-cd}
\usepackage{hyperref}
\usetikzlibrary{backgrounds}
\usetikzlibrary{arrows}
\usetikzlibrary{arrows.meta}
\usetikzlibrary{shapes,shapes.geometric,shapes.misc}
\usetikzlibrary{fit}
\usetikzlibrary{positioning}
\usetikzlibrary{calc}
\tikzstyle{tikzfig}=[baseline=-0.25em,scale=0.5]
\pgfdeclarelayer{edgelayer}
\pgfdeclarelayer{nodelayer}
\pgfsetlayers{background,nodelayer,edgelayer,main}
\tikzstyle{none}=[inner sep=0mm]
\tikzstyle{every loop}=[]
\tikzset{>={latex[width=1mm,length=1mm]}}
\newcommand{\drawThm}[3]{ \draw #1 %
	\ifx&#2&%
	\else
	node[style=thmref, pos=0.01] {\resizebox{5mm}{3mm}{#2}} %
	\fi
	\ifx&#3&%
	\else
	node[style=thmref, pos=0.99] {\resizebox{5mm}{3mm}{#3}}
	\fi ;
}

\tikzstyle{box}=[fill={rgb,255: red,228; green,228; blue,228}, draw=black, shape=rectangle]

\tikzstyle{reducible}=[->, dashed]
\tikzstyle{strictreducible}=[->]
\tikzstyle{nonreducible}=[->, draw=red]
\tikzstyle{uncomp}=[draw=red, <->]
\usepackage{hyperref}
\hypersetup{
	colorlinks=true,
	linkcolor=red,
	filecolor=magenta,      
	urlcolor=cyan,
	pdftitle={Overleaf Example},
	pdfpagemode=FullScreen,
}

\makeatletter
\newtheorem*{rep@theorem}{\rep@title}
\newcommand{\newreptheorem}[2]{%
\newenvironment{rep#1}[1]{%
 \def\rep@title{#2 \ref{##1}}%
 \begin{rep@theorem}}%
 {\end{rep@theorem}}}
\makeatother

\newtheorem{theorem}{Theorem}[section]
\newreptheorem{theorem}{Theorem}
\newtheorem{corollary}[theorem]{Corollary}
\newreptheorem{corollary}{Corollary}
\newtheorem{lemma}[theorem]{Lemma}

\newtheorem{proposition}[theorem]{Proposition}
\newreptheorem{proposition}{Proposition}

\newtheorem*{claimn}{Claim}

\theoremstyle{definition}
\newtheorem{definition}[theorem]{Definition} 
\newtheorem{remark}[theorem]{Remark}

\newtheorem{openquestion}[theorem]{Open Question}

\newcommand{\baire}{\omega^{\omega}}
\newcommand{\ramsey}{[\omega]^{\omega}}
\newcommand{\upto}{{\upharpoonright}}
\newcommand{\bSigma}{\bm{\Sigma}}
\newcommand{\bPi}{\bm{\Pi}}
\newcommand{\bDelta}{\bm{\Delta}}
\newcommand{\multif}{\rightrightarrows}
\newcommand{\subc}{:\subseteq}

\newcommand{\bmd}{\beta\text{-}\mathsf{mod}}
\newcommand{\homo}{\mathrm{hom}}
\newcommand{\hj}{\mathsf{HJ}}
\newcommand{\real}[1]{#1^{\mathrm{r}}}
\newcommand{\wfind}{\mathsf{wFindHS}}
\newcommand{\find}{\mathsf{FindHS}}
\newcommand{\rt}{\mathsf{RT}}
\newcommand{\wf}{\mathsf{WF}}
\newcommand{\CBaire}{\mathsf{C}_{\omega^{\omega}}}
\newcommand{\UCBaire}{\mathsf{UC}_{\omega^{\omega}}}
\newcommand{\SRT}[1]{\bSigma^0_{#1}\text{-}\mathsf{RT}}
\newcommand{\DRT}[1]{\bDelta^0_{#1}\text{-}\mathsf{RT}}
\newcommand{\dRT}[1]{\Delta^0_{#1}\text{-}\mathsf{RT}}
\newcommand{\sRT}[1]{\Sigma^0_{#1}\text{-}\mathsf{RT}}

\newcommand{\rca}{\mathsf{RCA}_0}
\newcommand{\aca}{\mathsf{ACA}_0}
\newcommand{\atr}{\mathsf{ATR}_0}
\newcommand{\pioo}{\bPi^1_1{-}\mathsf{CA}_0}

\newcommand{\achoice}{\bSigma^1_1\mathrm{Choice}}
\newcommand{\cachoice}{\bPi^1_1\mathrm{Choice}}
\newcommand{\tr}{\boldsymbol{\mathrm{Tr}}}
\newcommand{\T}{\mathrm{T}}
\newcommand{\ck}{\mathrm{CK}}
\newcommand{\id}[1]{\mathsf{id}_#1}
\newcommand{\W}{\mathrm{W}}
\newcommand{\sW}{\mathrm{sW}}
\newcommand{\ari}{\mathrm{a}}
\newcommand{\J}{\mathsf{J}}
\newcommand{\ra}[1]{#1{-}\mathsf{RA}}
\newcommand{\wstar}{\operatorname{\tilde{\star}}}
\DeclareMathOperator{\dom}{dom}
\DeclareMathOperator{\ran}{ran}

\title[The Galvin-Prikry Theorem in the Weihrauch lattice]{The Galvin-Prikry Theorem\\ in the Weihrauch lattice}
\author{Alberto Marcone}
\address
  {Dipartimento di Scienze Matematiche, Informatiche e Fisiche\\
  Universit\`a di Udine\\
  33100 Udine\\
  Italy}
\email{\href{mailto:alberto.marcone@uniud.it}{alberto.marcone@uniud.it}}
\author{Gian Marco Osso}
\address
  {Dipartimento di Scienze Matematiche, Informatiche e Fisiche\\
  Universit\`a di Udine\\
  33100 Udine\\
  Italy}
\email{\href{mailto:osso.gianmarco@spes.uniud.it}{osso.gianmarco@spes.uniud.it}}

\thanks{The authors thank Yudai Suzuki and Keita Yokoyama for useful conversations on the topics of the paper. 
The authors were partially supported by the Italian PRIN 2022 ``Models, sets and classifications'', prot.\ 2022TECZJA, funded by the European Union - Next Generation EU}

\date{\today}

\begin{document}
	
	\begin{abstract}
		This paper classifies different fragments of the Galvin-Prikry theorem, an infinite dimensional generalization of Ramsey's theorem, in terms of their uniform computational content (Weihrauch degree).
		It can be seen as a continuation of \cite{marconevalenti}, which focused on the Weihrauch classification of functions related to the open (and clopen) Ramsey theorem.
		We show that functions related to the Galvin-Prikry theorem for Borel sets of rank $n$ are strictly between the $(n+1)$-th and $n$-th iterate of the hyperjump operator $\hj$, which is in turn equivalent to the better known $\widehat{\wf}$, which corresponds to $\pioo$ in the Weihrauch lattice.
		To establish this classification we obtain the following computability theoretic result: a Turing jump ideal containing homogeneous sets for all $\Delta^0_{n+1}(X)$ sets must also contain $\hj^n(X)$. 
  We also extend our analysis to the transfinite levels of the Borel hierarchy.
  We further obtain some results about the reverse mathematics of the lightface fragments of the Galvin-Prikry theorem.
	\end{abstract}
	
	\maketitle
	
	\section{Introduction}
	
	The study of Weihrauch degrees is a framework which allows us to classify relations on so-called \emph{represented spaces} in terms of their computational content.
	Starting with \cite{GherardiMarcone}, it has become clear that this investigation can be turned into a study of the computational content of mathematical theorems as follows.
	Any theorem of the form $\forall x \in A (\varphi(x) \rightarrow \exists y \in B \,  \psi(x,y))$ can be identified, in a way reminiscent of Skolem functions, with the \emph{partial multi-valued function} $f \subc A \multif B$ (in informal discussion we abuse language and refer to these simply as functions) defined as $f(x)=\{y \in B : \psi(x,y)\}$ where $\dom(f)=\{x \in A : \varphi(x)\}$.
	If $f_1$ and $f_2$ are two functions associated with theorems $T_1$ and $T_2$ as above, we regard a Weihrauch reduction (Definition \ref{defweihrauch}) from $f_1$ to $f_2$, denoted $f_1 \leq_{\W} f_2$, as information about the computational relationship between theorems $T_1$ and $T_2$.
	Weihrauch classification has been successfully applied to a wide range of mathematical results, yielding insight on their computational properties.
	
	 This study of the computational content of theorems can be fruitfully connected with reverse mathematics, the study of mathematical results in terms of proof theoretic strength. 
	This is the case, roughly speaking, thanks to the deep ties between reverse mathematics and computability theory. 
	Since Weihrauch classification allows us to detect finer differences than those visible through the lens of reverse mathematics, it has historically been useful to zoom on particular levels of the reverse mathematical hierarchy and look at the Weihrauch degrees of their inhabitants.
	In the beginning of this line of research, most of the work was done on theorems which, in terms of reverse mathematics, require fairly weak axioms (see, e.g., \cite{GherardiMarcone}, \cite{brattkagherardi}, \cite{connectedchoice}, \cite{additiveramsey}). 
	In the last few years researchers initiated the study of Weihrauch degrees of theorems requiring strong axioms, say at the levels of $\atr$ and $\pioo$, and there have been several developments (see \cite{leafmanagement}, \cite{kmp}, \cite{gohembeddings}, \cite{marconevalenti}, \cite{dauriac}, \cite{CBWeihrauch}, \cite{li2024}).
	The present paper belongs to this line of research. In this context, it is often natural to substitute computable functions with arithmetical ones in the definition of Weihrauch reducibility, obtaining arithmetical Weihrauch degrees.
	
	We classify fragments of the Galvin-Prikry theorem, which states that all Borel subsets of $\ramsey$ are Ramsey (Definition \ref{def:ramsey}).
	For each pointclass $\Gamma$ in $\{\bSigma^0_{\alpha}, \bPi^0_{\alpha}, \bDelta^0_{\alpha} : 0 < \alpha< \omega_1\}$ we consider the restriction of the Galvin-Prikry theorem to sets in $\Gamma$.
 When $\Gamma$ is either $\bDelta^0_1$ or $\bSigma^0_1$ this restriction is equivalent to $\atr$, while for $\bDelta^0_n$ or $\bSigma^0_n$ with $1<n< \omega$ we obtain statements equivalent to $\pioo$ (see Section \ref{reversemath} for details).
 In the Weihrauch setting, we associate three different functions to the statement ``sets in $\Gamma$ are Ramsey'' (this follows the setup of \cite{marconevalenti}, which focuses on the case $n=1$) and classify their Weihrauch degrees.
	We show that two out of three of the functions related to the Ramseyness of sets in one of the classes $\{\bSigma^0_n, \bPi^0_n, \bDelta^0_n : n \in \omega\}$ lie below the $n$-th iterate of the hyperjump operator and above the ($n-1$)-th, thus showing that the strictly increasing hierarchy given by the iterates of the hyperjump corresponds to more and more sets being Ramsey.
 For the remaining function, we obtain only the lower bound. We obtain similar results also for the transfinite levels of the Borel hierarchy.
	Our classification rests on new results of computability theoretic nature in the style of the fundamental paper on Ramsey's theorem with finite exponent \cite{jockushramsey}, and draws from known results in the reverse mathematics of the Galvin-Prikry theorem.\smallskip
    
    W now describe the structure of the paper and mention the main results. In Section \ref{basicnotions} we fix some notation and introduce background material in reverse mathematics, computable analysis, and on the Galvin-Prikry theorem.
    Further, we also introduce the new operator $\wstar$, which, in some situations, represents the compositional product in the strong Weihrauch degrees.
    
    In Section \ref{bmd} we introduce the multi-valued function $\bmd$ which maps a set $X$ to the set of pairs $(W,h)$ where $W$ is a code for a countable $\beta$-model containing $X$ and $h$ is a $\Sigma^1_1$ truth predicate for $W$. We rephrase a theorem of Simpson to obtain: 
    \begin{repproposition}{degreeofbmd}
        $\bmd \equiv_{\sW} \widehat{\wf}$.
    \end{repproposition}
    
    $\beta$-models and $\Sigma^1_1$ truth predicates allow us to perform countably many sequential choices from $\bSigma^1_1$ (or $\bPi^1_1$) sets which have codes in the model we are considering. This is key in several of our constructions, notably in the proof of Proposition \ref{Prop:transfiniteub} (the technical result behind Theorems \ref{thm:transub1} and \ref{thm:transub2}).
    
    In Section \ref{secseparation} we exploit ideas implicit in \cite{everyhigherdegree} and \cite{jockushramsey} to characterize in terms of hyperjumps the Turing jump ideals (resp.\ $\mathsf{HYP}$-ideals) containing homogeneous subsets for all sets of Borel rank $n$ (resp.\ $\alpha \geq \omega$). We obtain:
    \begin{reptheorem}{reversal}
        Let $n \geq 1$, let $X \subseteq \omega$ and let $\mathcal{I}$ be a Turing jump ideal such that $X \in \mathcal{I}$ and, for any $\Delta^0_{n+1}(X)$ definable set $A$, there exists $h \in [\omega]^{\omega} \cap \mathcal{I}$ such that $[h]^{\omega} \subseteq A$ or $[h]^{\omega} \cap A = \emptyset$. Then $\hj^n(X) \in \mathcal{I}$. 
    \end{reptheorem}
    \begin{reptheorem}{thm:transreversal}
        Let $X \subseteq \omega$, let $\omega \leq \alpha \leq \omega^X_1$, and let $\mathcal{I}$ be a $\mathsf{HYP}$-ideal such that $X \in \mathcal{I}$ and, for any $\Delta^0_{1+\alpha}(X)$ definable set $A$, there exists $h \in [\omega]^{\omega} \cap \mathcal{I}$ such that $[h]^{\omega} \subseteq A$ or $[h]^{\omega} \cap A = \emptyset$. Then $\hj^{\alpha}(X) \in \mathcal{I}$. 
    \end{reptheorem}

    The last two theorems are a higher level version of Jockusch's results about the Ramsey theorem for finite exponents: as the complexity of the Borel set increases the number of hyperjumps needed to compute solutions increases, mirroring the increase of the number of Turing jumps needed to compute solutions to instances of the Ramsey theorem for larger and larger finite exponents.
    
    These results can be seen as lower bounds to the complexity of sets in ideals containing homogeneous sets for partitions of a given Borel complexity. Corresponding upper bounds for sets of finite Borel rank are known from the literature (see, e.g. Corollary \ref{Cor:degupperbounds} and the results mentioned before it). We obtain new level by level upper bounds for partitions of transfinite Borel rank:

    \begin{reptheorem}{thm:transub1}
        Let $X \subseteq \omega$, $\alpha < \omega^{X}_1$ and let $\mathcal{M}$ be an $\omega$-model of $\atr$ containing $\hj^{\omega^{1+\alpha}}(X)$. Then $\mathcal{M} \vDash \Sigma^0_{\omega+\alpha}(X)$-$\mathsf{RT}$ (all sets which are $\Sigma^0_{\omega+\alpha}$ definable from $X$ are Ramsey from the point of view of $\mathcal{M}$).	
    \end{reptheorem}
	
    \begin{reptheorem}{thm:transub2}
        Let $X \subseteq \omega$, $\alpha < \omega^{X}_1$ and let $A \subseteq \ramsey$ be a $\Sigma^0_{\omega+\alpha}(X)$ set. There exists $g \leq_{\T} \hj^{\omega^{1+\alpha}+1}(X)$ such that $[g]^{\omega} \subseteq A$ or $[g]^{\omega} \cap A = \emptyset$.
    \end{reptheorem}

	In Section \ref{upperbounds} we draw from the reverse mathematics literature to exploit chains of $\beta$-models in the computation of solutions to instances of the Galvin-Prikry theorem of finite Borel rank, thus obtaining upper bounds for the Weihrauch degrees of the functions we consider. 
    We combine these with some of the results of Section \ref{secseparation} and with results of \cite{marconevalenti} to obtain a fairly complete picture (see Figures \ref{figure1} and \ref{figure2}) of the (arithmetical) Weihrauch degree of functions related to fragments of the Galvin-Prikry theorem. The great number of incomparabilities appearing in Figure \ref{figure2} is to a large extent due to the following:
    \begin{repcorollary}{nohyperjumps}
        Let $\alpha < \omega_1$,  $\hj \nleq_{\W} \SRT {\alpha}$.
    \end{repcorollary}
    This result shows a stark contrast between the picture of the Weihrauch degrees and the arithmetical Weihrauch degrees of these functions: indeed, we have that for all $n \in \omega$, $\hj^n \leq^{\ari}_{\W} \SRT {n}$.
    The separation of Corollary \ref{nohyperjumps} is obtained (using Proposition \ref{prop:nononhyp}) exploiting an old characterization of the hyperarithmetical sets as the \emph{recursively encodable} ones (i.e.\ $A$ is hyperarithmetical if and only if every infinite set $B$ has an infinite subset $C$ with $A \leq_{\T} C$), together with a property of the relations $\SRT {\alpha}$ typical of Ramsey-like statements which we called \emph{meeting every infinite set} (Definition \ref{def:mes}).
    
 	\begin{figure}[h]
		\begin{tikzcd}
			\find_{\bSigma^0_{k+1}} &                                                                                                                   & {\widehat{\wf}^{[k+1]}}      \\
			& {\SRT 1 \star \widehat{\wf}^{[k]}} \arrow[lu, "\text{\ref{Cor:findaboveall}, \ref{Cor:findstrictabove}}"'] \arrow[ru, "\text{\ref{opencase}, \ref{Prop:hjstrictlyabove}}"'] \arrow[red, dd, bend right=60, "\text{\ref{Lem:MV415}}"']                                &                              \\
			\find_{\bPi^0_{k+1}}    & \SRT {k+1} \arrow[u, dashed, "\text{\ref{slightlylower}}"']                                                                                      & \DRT {k+1} \arrow[l, dashed] \\
			& {\CBaire \star \widehat{\wf}^{[k]} \equiv^{\ari}_{\W} \wfind_{\bPi^0_{k+1}}} \arrow[u, dashed] \arrow[lu, dashed] &                              \\
			& {\UCBaire \star \widehat{\wf}^{[k]} \equiv^{\ari}_{\W} \wfind_{\bSigma^0_{k+1}}} \arrow[u, "\text{\ref{arisepwfind}}"']                        &                              \\
			& \wfind_{\bDelta^0_{k+1}} \arrow[u, dashed] \arrow[ruuu, dashed, bend right=40]                                       &                              \\
			& {\widehat{\wf}^{[k]}} \arrow[u, "\text{\ref{computingmorehyperjumps}}"']                                                                                   &                             
		\end{tikzcd}
		\caption{\label{figure1}This table is a picture of the arithmetical Weihrauch degrees of functions related to the Galvin-Prikry theorem for Borel sets of rank $k+1$. A solid black arrow from $f$ to $g$ means $f <^{\ari}_{\W} g$, while a dashed one denotes $f \leq^{\ari}_{\W} g$ (and $g\leq^{\ari}_{\W} f$ is open) and a red one indicates $f \nleq^{\ari}_{\W} g$.  Labeled arrows reference the proofs of the reductions and/or separations. Unlabeled arrows correspond to trivial reductions and/or separations.}
	\end{figure}

	\begin{figure}[h]
		\begin{tikzcd}
			& \find_{\bSigma^0_{k+1}}                                            &                                                                          & {\widehat{\wf}^{[k+1]}}      \\
			\find_{\bPi^0_{k+1}} &                                                                    & {\SRT 1 \star \widehat{\wf}^{[k]}} \arrow[lu, "\text{\ref{Cor:upleftarrow}}"'] \arrow[ru, "\text{\ref{opencase}, \ref{Prop:hjstrictlyabove}}"']                 &                              \\
			& {\widehat{\wf}^{[k]}} \arrow[red, r, "\text{\ref{nohyperjumps}}"'] \arrow[lu, "\text{\ref{Cor:findWLB}}"'] \arrow[ru] & \SRT {k+1} \arrow[u, "\text{\ref{Cor:strict}}"']                     &                              \\
			& \wfind_{\bPi^0_{k+1}} \arrow[luu, "\text{\ref{Cor:strongsep}}"] \arrow[ru, dashed] \arrow[red, r, "\text{\ref{arisepwfind}}"']              & \wfind_{\bSigma^0_{k+1}} \arrow[u, dashed]                               & \DRT{k+1} \arrow[lu, dashed] \\
			&                                                                    & \wfind_{\bDelta^0_{k+1}} \arrow[lu] \arrow[u, dashed] \arrow[ru, dashed] \arrow [red, luu, bend left=40, dash pattern=on 59pt off 15pt, "\text{\ref{computingmorehyperjumps}}"] &                             
		\end{tikzcd}
			\caption{\label{figure2}This table is a picture of the Weihrauch degrees of functions related to the Galvin-Prikry theorem for Borel sets of rank $k$. A solid black arrow from $f$ to $g$ indicates $f <_{\W} g$, while a dashed one indicates $f \leq_{\W} g$ (and $g\leq_{\W} f$ is open) and a red one indicates $f \nleq_{\W} g$. Labeled arrows reference the proofs of the reductions and/or separations. Unlabeled arrows correspond to trivial reductions and/or separations.}
	\end{figure}
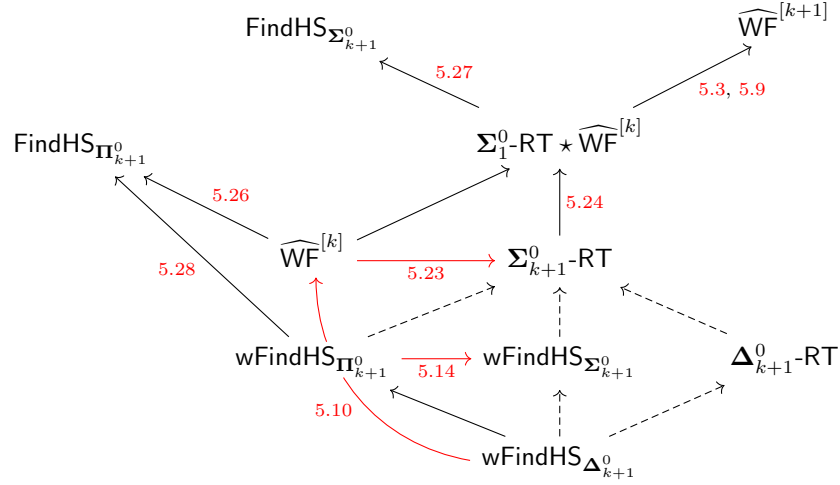
 
	The characterization obtained in Theorem \ref{reversal} is fundamental to many of our Weihrauch separations, but it is also used in Section \ref{revmath} to study the reverse mathematics of the Galvin-Prikry theorem restricted to lightface definable sets. We get:

    \begin{reptheorem}{equivalences}
		For any $k \in \omega$, the following are equivalent for $\beta$-models:
		\begin{enumerate}
			\item $\hj^k(\emptyset)$ exists,
			\item there exists a sequence of countable coded $\beta$-models $M^1 \in M^2 \in \dots \in M^k$,
			\item the lightface $\Sigma^0_{k+1}$ Galvin-Prikry theorem,
			\item the lightface $\Delta^0_{k+1}$ Galvin-Prikry theorem.
		\end{enumerate}
    \end{reptheorem}

    In fact, we can show that items 1 to 4 are actually equivalent over a theory $T \supseteq \atr$ which is strictly weaker than $\pioo$ (Corollary \ref{Cor:T}).
    This provides an example of how Weihrauch reducibility can lead to a refinement of reverse mathematical results.
	Incidentally, we note that in \cite{suzukiyokoyama}, which appeared on arXiv while this paper was being finalized, the authors also study the reverse mathematics of fragments of the lightface Galvin-Prikry theorem.
	The differences between our results, which are a consequence of our different approaches, are discussed in Section \ref{revmath}.

    \section{Basic notions}\label{basicnotions}
 
	We adopt the notational convention, common in reverse mathematics, to reserve the symbol $\omega$ for the actual set of natural numbers, while the symbol $\mathbb{N}$ refers to the universe of numbers of a given structure in the language of second order arithmetic.
	Accordingly we only use $\mathbb{N}$ in formulas in the language of second order arithmetic, and we use $\omega$ everywhere else.
	We assume familiarity with basic concepts of descriptive set theory (e.g.\ Borel and arithmetically definable subsets of Baire space, see \cite{moschovakis}) and computability theory (Turing reducibility $\leq_{\T}$, arithmetical reducibility $\leq_{\ari}$, the jump operator, Turing jump ideals, type 2 computability for functions on the Baire space, see \cite{rogers} and \cite{Weihrauch}).
 
	In particular we (sometimes silently) make ample use of pairing and projection functions to identify finite sequences of natural numbers with single natural numbers, or countably many elements of Baire space with a single one.
    We denote pairing functions with the symbol $\langle \cdot , \cdot \rangle$; for projections, if $p \in \baire$ and $W \subseteq \omega$, we write $p_i$ for the function $n \mapsto p(\langle i,n \rangle)$ and $W_i$ for the set $\{n \in \omega : \langle i,n \rangle \in W\}$.
    Given two strings $x \in \omega^{<\omega}$ and $y \in \baire \cup \omega^{<\omega}$, we use $x^{\smallfrown}y$ to denote their concatenation, and we write $x \sqsubseteq y$ to mean that $x$ is an initial segment of $y$.
    These operations are well-known and, each in the appropriate sense, effective.
    
    We use $\mathcal{O}$ to denote Kleene's set of ordinal notations. 
    We denote by $\omega^{\ck}_1$ the least ordinal which does not have an computable copy. 
    For $a \in \mathcal{O}$, we say $|a|=\alpha$ if $a$ is a notation for the ordinal $\alpha < \omega^{\ck}_1$ (recall that $|1|=0$, $|2^b|=|b|+1$ and $|3 \cdot 5^e| = \sup \{|\{e\}(n)| : n \in \omega\}$).
    Kleene's $\mathcal{O}$ can be relativized to an oracle $X$, obtaining $\mathcal{O}_X$ (the function $X \mapsto \mathcal{O}_X$ is the \emph{hyperjump}), $\omega^X_1$ and $|\cdot|_X$.
    
    For any $Y$ and $a \in \mathcal{O}_X$ we use the notation $Y^{(a)}$ and $\hj^a(Y)$ to denote the iteration of, respectively, the jump and hyperjump operator along $|a|_X$. 
    Sometimes we abuse the notation and write $Y^{(\alpha)}$ in place of $Y^{(a)}$ for some $a$ such that $|a| = \alpha$, trusting that the context disambiguates the different uses. We denote hyperarithmetical reducibility as $\leq_h$. Recall that $X \leq_h Y$ if and only if there is some $a \in \mathcal{O}_Y$ such that $X \leq_T Y^{(a)}$. We say that a subset of $\mathcal{P}(\omega)$ is a $\mathsf{HYP}$-ideal if it is closed under $\oplus$ and $\leq_h$. For more information and notational conventions about ordinal notations, hyperarithmetical sets, and the hyperjump, see \cite{sacks2017}.
	
	\subsection{The Galvin-Prikry theorem and the Ramsey space}
	
	Ramsey's theorem \cite{ramsey} states that, given any natural number $k$ and any partition of the set $[\omega]^k$ into disjoint subsets $A_0$ and $A_1$, there exists an infinite set $H \subseteq \omega$ which is homogeneous for the partition, meaning $[H]^k \subseteq A_i$ for some $i<2$.
	It is natural to ask whether a similar result holds for partitions of the set of infinite subsets of $\omega$, i.e.\ whether it is the case that, whenever $\ramsey$ is partitioned into disjoint subsets $A_0$ and $A_1$, there is an infinite set $H \subseteq \omega$ and an index $i<2$ such that $[H]^{\omega} \subseteq A_i$.
	Using the axiom of choice it is easy to see that this is not the case.
	On the other hand, Galvin and Prikry proved \cite{galvinprikry}:
	
	\begin{theorem}
		Let $\ramsey =A_0 \cup A_1$ be a partition of $\ramsey$ into Borel sets.
		There exists some $H \subseteq \omega$ which is homogeneous for the partition, i.e.\ such that $[H]^{\omega} \subseteq A_i$ for some $i < 2$.
	\end{theorem}
	
	Shortly after the publication of this proof, Silver showed that the result can be extended to partitions consisting of a $\bSigma^1_1$ and a $\bPi^1_1$ set \cite[Theorem 1]{silveranalyticramsey}.
	On the other hand, it is consistent with $\mathsf{ZFC}$ that there is a partition of $\ramsey$ into two $\Delta^1_2$ sets which admits no infinite homogeneous set (see e.g.\ page 1 of \cite{silveranalyticramsey}).
	
	Both Ramsey's theorem and Galvin and Prikry's theorem hold for partitions comprising an arbitrary finite number of pieces. Here and in the rest of the paper we focus on partitions consisting of two sets for simplicity. In this context, a partition is uniquely determined by any of the two sets.
	
	From now on we refer to the set $\ramsey$ of infinite subsets of $\omega$ as the \emph{Ramsey space}, and we frequently identify $X \in \ramsey$ with its principal function $p_X \colon \omega \rightarrow \omega$, i.e.\ with the unique function enumerating $X$ in increasing order\footnote{From the perspective of computability, we note that $X \equiv_{\T} p_X$.}.
	This identification allows us to see the Ramsey space $\ramsey$ as a closed subset of Baire space $\baire$. The topology of the Ramsey space is then generated by the subbasis $\{[s] \cap \ramsey : s \in [\omega]^{<\omega}\}$ where $[\omega]^{<\omega}$ denotes the set of finite increasing sequences of natural numbers, and, for any such sequence $s$, $[s]=\{f \in \baire : s \sqsubset f\}$.
	This topology is Polish as it is the subspace topology of a closed subset of a Polish space.
	For any $f \in \ramsey$ we denote as $[f]^\omega$ the set of infinite subsets of $\ran(f)$, where again we think of sets via their principal functions.
	An advantage of this point of view is that $[f]^\omega=\{f \circ g : g \in \ramsey\}$.
	
	Given a set $A \subseteq \ramsey$, we say that $A$ \emph{is Ramsey} if there is $f \in \ramsey$ with $[f]^{\omega} \subseteq A$ or $[f]^{\omega} \subseteq (\ramsey \setminus A)$. In the first case we say that $f$ \emph{lands in} $A$ and in the second case we say that $f$ \emph{avoids} $A$, and in both cases we say that $f$ is an \emph{infinite homogeneous set} for $A$. Given $A \subseteq \ramsey$ we denote the set of homogeneous sets for $A$ by $\mathrm{HS}(A)$.
	
	\begin{definition}\label{def:ramsey}
		Let $\Gamma \subseteq \mathcal{P}(\ramsey)$ be a class of subsets of the Ramsey space.
		We say that $\Gamma$ is \emph{Ramsey} if every $A \in \Gamma$ is Ramsey and abbreviate such statement with $\Gamma$-$\rt$.
	\end{definition}
	
	In these terms, the Galvin-Prikry theorem says that $\bDelta^1_1$-$\rt$ holds, while Silver's theorem says that $\bSigma^1_1$-$\rt$ holds.
	
	In the rest of the paper it will be convenient to work in spaces of the form $[f]^{\omega}$ for some $f \in \ramsey$ rather than in the whole space $\ramsey$.
	We show that this operation is not problematic and prove three simple lemmas to be used in the following.
	\begin{lemma}\label{behaviourofintersections}
		Let $A \subseteq \ramsey$ be any set and $f \in \ramsey$.
		If no set $g \in \ramsey$ avoids (resp.\ lands in) $A$, then no set $g \in [f]^{\omega}$ avoids (resp.\ lands in) $A \cap [f]^{\omega}$. 
	\end{lemma}
	\begin{proof}
		Assume no $g \in \ramsey$ avoids $A$, then in particular for every $g \in [f]^{\omega}$ there exists $h \in \ramsey$ such that $g \circ h \in A$.
		Clearly $g \circ h \in [f]^{\omega}$, hence no set in $[f]^{\omega}$ avoids $A \cap [f]^{\omega}$.
		The proof for sets landing in $A$ is essentially the same.
	\end{proof}
	
	Other than intersection with $[f]^{\omega}$, another way to move between subsets of $\ramsey$ and $[f]^{\omega}$ is exploiting the natural homeomorphism given by $g \mapsto f \circ g$.
	
	\begin{lemma}\label{behaviourofcompositions}
		Let $A \subseteq \ramsey$ be any set, $f \in \ramsey$ and $F \colon \ramsey \rightarrow [f]^{\omega}$ given by $h \mapsto f \circ h$.
		A set $g \in \ramsey$ lands in (resp.\ avoids) $A$ if and only if $f \circ g$ lands in (resp.\ avoids) $F[A]$. 
	\end{lemma}
	\begin{proof}
		We have that $g$ lands in $A$ if and only if $\forall h \in \ramsey \, g \circ h \in A$.
		This is the case if and only if $ \forall h \in \ramsey \, f \circ g \circ h \in \{f \circ k : k \in A\}= F[A]$ (the backward direction holds because, since $f$ is injective, if $f \circ g \circ h = f \circ k$ then $g \circ h =k$). The proof for sets avoiding $A$ is the same.
	\end{proof}
	
	\begin{corollary}\label{relativization}
		Let $B$ be a Borel subset of $\ramsey$ and let $f \in \ramsey$.
		Then there exists some $g \in [f]^{\omega}$ homogeneous for $B$.
	\end{corollary}
	\begin{proof}
		Consider the Borel set $B \cap [f]^{\omega}$ and look at $F^{-1}[(B \cap [f]^\omega)]$ (where again $F(g)=f \circ g$): this is also a Borel set in $\ramsey$, so by the Galvin-Prikry theorem there is some $h \in \ramsey$ which is homogeneous for it.
		By Lemma \ref{behaviourofcompositions} we have that $f \circ h$ is homogeneous for $B \cap [f]^{\omega}$. Now if $f \circ h$ lands in $B \cap [f]^{\omega}$, then a fortiori it lands in $B$. If $f \circ h$ avoids $B \cap [f]^{\omega}$, then $[f \circ h]^{\omega} \cap (B \cap [f]^{\omega}) = \emptyset$, which implies that $[f \circ h]^{\omega} \cap B= \emptyset$ as $[f \circ h]^{\omega} \subseteq [f]^{\omega}$.
	\end{proof}
	
	\subsection{Reverse mathematics}\label{reversemath}
	
	We introduce the notions of reverse mathematics which play a role in the rest of the paper.
	Our source is Simpson's book \cite{simpson}, to which we refer the reader for more details.
	
	Reverse mathematics is concerned with studying implications between mathematical principles over weak theories in the language of second order arithmetic.
	A model $\mathcal{M}$ for this language is a tuple $\langle M, \mathcal{S}, +_M, \cdot_M, 0_M, 1_M \rangle$.
	The theory most commonly used as base is $\rca$, consisting of the axioms for commutative, discretely ordered semirings, the axiom schema of $\bDelta^0_1$-comprehension and the axiom schema of $\bSigma^0_1$-induction.
    Other theories which feature in the paper are, in increasing order of strength, $\aca$, which is $\rca$ + the axiom schema of comprehension restricted to arithmetical formulas, $\atr$, which is $\rca$ + the axiom schema of arithmetical transfinite recursion (equivalently, the axiom asserting the existence of a Turing jump hierarchy along any well-ordered set) and lastly $\pioo$, which is $\rca$ + the axiom schema of comprehension restricted to $\bPi^1_1$ formulas.
	
	There are well known computability theoretic interpretations of the theories above, for example, over $\rca$, the theory $\aca$ is equivalent to the axiom stating that every set has a Turing jump \cite[Exercise VIII.1.12]{simpson}, while $\pioo$ is equivalent to the axiom stating that every set has a hyperjump \cite[Exercise VII.1.16]{simpson}.
	
	A structure $\langle M, \mathcal{S}, +_M, \cdot_M, 0_M, 1_M \rangle$ in the language of second order arithmetic is called an \emph{$\omega$-model} if $M=\omega$ and the ordering, algebraic operations and constants coincide with their natural interpretations in $\omega$.
	It is customary to identify an $\omega$-model with just the set $\mathcal{S} \subseteq \mathcal{P}(\omega)$.
	
	Given $\mathcal{M}=\langle M, \mathcal{S}, +_M, \cdot_M, 0_M, 1_M \rangle$, we say that $\mathcal{N}=\langle N, \mathcal{Q}, +_N, \cdot_N, 0_N, 1_N \rangle$ is a \emph{$\beta$-submodel} of $\mathcal{M}$ if $M=N$, $\mathcal{Q} \subseteq \mathcal{S}$, the interpretations of the operations, the order and the constants in the two models coincide, and for every $\Sigma^1_1$ sentence $\varphi$ with parameters from $\mathcal{Q}$, $\mathcal{M} \vDash \varphi$ if and only if $\mathcal{N} \vDash \varphi$.
    Notice that a $\beta$-submodel of an $\omega$-model is itself an $\omega$-model.
	We call \emph{$\beta$-models} the $\beta$-submodels of the $\omega$-model $\mathcal{P}(\omega)$, so by definition if $\mathcal{M}$ is a $\beta$-model and $\varphi$ is a $\Sigma^1_1$ sentence with parameters from $\mathcal{S}$, then $\mathcal{M} \vDash \varphi$ if and only if $\mathcal{P}(\omega) \vDash \varphi$.
	We refer to this property as $\Sigma^1_1$ (or $\Pi^1_1$) \emph{correctness}.
	
	Given a model $\mathcal{M}=\langle M, \mathcal{S}, +_M, \cdot_M, 0_M, 1_M \rangle$, we say that $W \in \mathcal{S}$ \emph{codes a countable $\beta$-model} if the set $\mathcal{S}'=\{W_i: i \in M\}$, where $W_i=\{n \in M : \langle i,n \rangle \in W\}$, is a $\beta$-submodel of $\mathcal{M}$.
    In this case say that a set $X \subseteq M$ belongs to $W$, and write $X \in W$, if $X=W_i$ for some $i \in M$.
	Over $\aca$, the theory $\pioo$ is equivalent to the assertion that any set $X$ belongs to some countable coded $\beta$-model $W$ \cite[Theorem VII.2.10]{simpson}. Moreover, any $\beta$-model is a model of $\atr$ \cite[Theorem VII.2.7]{simpson}.
	
	We now recall the known results in the reverse mathematics of (fragments of) the Galvin-Prikry theorem.
	First we note that, if $\Gamma$ is a class of formulas, then the schema \[\exists f \in [\mathbb{N}]^{\mathbb{N}} \left[\forall g\in [\mathbb{N}]^{\mathbb{N}} \, \varphi(f \circ g) \lor \forall g\in [\mathbb{N}]^{\mathbb{N}} \,  \neg\varphi(f \circ g)\right], \,\,\, \varphi \in \Gamma\]
	expresses the fact that every $\Gamma$-definable subset of $[\mathbb{N}]^{\mathbb{N}}$ is Ramsey.
	Mimicking what we did in the $\mathsf{ZFC}$ context in the previous section, we denote such schema as $\Gamma$-$\rt$.\footnote{Using the same notation for a $\mathsf{ZFC}$ principle and a principle in second order arithmetic does not lead to confusion as the intended meaning is always clear from the context.}
	In these terms, the Galvin-Prikry theorem for open sets is expressed in second order arithmetic by the schema $\SRT1$ where a formula is $\bSigma^0_1$ if it is $\Sigma^0_1$ with possible set parameters.
	We will use analogous notation for other classes of definable sets, e.g., $\sRT 1$ denotes the Galvin-Prikry theorem for lightface open sets (those open sets which are definable without set parameters).
	\begin{theorem}[{\cite[Theorem V.9.7]{simpson}}]\label{openramsey}
		Over $\rca$, $\SRT1$ is equivalent to $\atr$.
	\end{theorem}
	
	\begin{lemma}[{\cite[Lemma VI.6.2]{simpson}}]\label{manybetamodels}
		For all $k \in \omega$, the following is provable in $\aca$.
		Let $M^1 \in M^2 \in \dots \in M^k$ be a sequence of nested countable coded $\beta$-models.
		Then for every $j \leq i \leq k$, $M^i$ satisfies the Galvin-Prikry theorem for sets which are $\Sigma^0_{i+1-j}$ in some parameter in $M^j$.
	\end{lemma}
	
	Exploiting Simpson's proof we can obtain the following corollary about the Ramsey property for lightface classes:
	\begin{corollary}\label{simpsonreversal}
		For all $k \in \omega$, the following is provable in $\atr$: if there exists a sequence of nested countable coded $\beta$-models $M^1 \in M^2 \in \dots \in M^k$, then $\sRT {k+1}$ holds.
	\end{corollary}
	
	\begin{proof}
		The proof of \cite[Lemma VI.6.2]{simpson} actually establishes that for any set $X$ and any $A \subseteq [\mathbb{N}]^{\mathbb{N}}$ with $A$ a $\Sigma^0_{k+1}(X)$ set, if $W$ codes a countable $\beta$-model and $X \in W$, then there are a set $Y \leq_{\T} W$ and a set $f \in [\mathbb{N}]^\mathbb{N}$ such that $f \leq_{\T} Y$, $Y$ codes a $\bSigma^0_k$ set $B$ and, for all $g \in \mathrm{HS}(B)$, $f \circ g \in \mathrm{HS}(A)$.
		If we start with a computable $X$ coding a $\Sigma^0_{k+1}$ set $A \subseteq [\mathbb{N}]^{\mathbb{N}}$, we can iterate this procedure $k$ many times as long as there is a sequence of $k$ nested countable coded $\beta$-models given by $M^1, M^2, \dots  M^k$ (notice that $X \in M^1$ because $X$ is computable).
		The end result is a set $Y_k \leq_{\T} M^k$ coding a $\bSigma^0_1$ set $C$ and a set $f_k \in [\mathbb{N}]^{\mathbb{N}}$ with $f_k \leq_{\T} Y_k$ such that, for all $g \in \mathrm{HS}(C)$, $f_k \circ g \in \mathrm{HS}(A)$.
		Since by Theorem \ref{openramsey} $\atr$ proves $\SRT1$, it follows that $\mathrm{HS}(C) \neq \emptyset$, hence $\mathrm{HS}(A) \neq \emptyset$.
		This shows that the Galvin-Prikry theorem for lightface $\Sigma^0_{k+1}$ sets holds. 
	\end{proof}
	
	We relate the existence of chains of nested $\beta$ models with the existence of iterated hyperjumps.
		
	\begin{lemma}\label{smallbetamodels}
		For any $k \in \omega$, $\aca$ proves that, if $\hj^{k}(\emptyset)$ exists, then there is a sequence $M^1 \in M^2 \in \dots \in M^k$ of countable coded $\beta$-models such that for all $i \leq k$, $M^i \leq_{\T} \hj^i(\emptyset)$ (so in particular $\hj^i(\emptyset) \notin M^i$). Conversely given any sequence of nested countable coded $\beta$-models $M^1 \in M^2 \in \dots \in M^k$, we have $\hj^i(\emptyset) \leq_{\ari} M^i$ for every $i \leq k$.
	\end{lemma}
	\begin{proof}
		We prove both statements simultaneously by (external) induction on $k$. By \cite[Corollary VII.2.12]{simpson} (where the set $X$ is instantiated with $\emptyset$), $\aca$ proves that if $\hj(\emptyset)$ exists, then there is a countable coded $\beta$-model $M^1$, and by construction $M^1 \leq_{\T} \hj(\emptyset)$.
		Conversely, if $M^1$ is any countable coded $\beta$-model, then is $n \in \hj(\emptyset) \leftrightarrow \neg \exists i [(M^1)_i \text{ codes a path in }T^{\emptyset}_n]$, so $\hj(\emptyset)$ is arithmetical in $M^1$.
		
		Now we work in $\aca$, we assume that $\hj^{k+1}(\emptyset)$ exists, and we assume inductively that $M^1 \in M^2 \in \dots \in M^k$ is a sequence of countable coded $\beta$-models such that, for every $i \leq k$, $M^i \leq_{\T} \hj^i(\emptyset)$. In particular $\hj(M^k) \leq \hj(\hj^k(\emptyset))=\hj^{k+1} (\emptyset)$, so $\hj(M^k)$ exists.
		Then, again by \cite[Corollary VII.2.12]{simpson}, there is a countable coded $\beta$-model $M^{k+1}$ containing $M^k$ such that $M^{k+1} \leq_{\T} \hj(M^k)$.
		Conversely, assume there is a sequence of nested countable coded $\beta$-models $M^1 \in M^2 \in \dots \in M^k \in M^{k+1}$, and inductively assume that $\hj^i(\emptyset) \leq_{\ari} M^i$ for all $i \leq k$.
		Since $\hj^k(\emptyset) \leq_{\ari} M^k \in M^{k+1}$ and $M^{k+1}$ is closed under arithmetical equivalence, it follows that $\hj^k(\emptyset) \in M^{k+1}$.
		Since $M^{k+1}$ is a $\beta$-model, it follows that again $n \in \hj^{k+1}(\emptyset) \leftrightarrow \neg \exists i [(M^{k+1})_i \text{ codes a path in }T^{\hj^k(\emptyset)}_n]$ so in particular $\hj^{k+1}(\emptyset)$ is arithmetical in $M^{k+1}$.
	\end{proof}
	
	We can put Corollary \ref{simpsonreversal} together with Lemma \ref{smallbetamodels} to obtain a degree theoretic upper bound for homogeneous sets of $\Sigma^0_k$ sets.
	
	\begin{corollary}\label{Cor:degupperbounds}
		Let $k \in \omega$, let $A \subseteq [\omega]^{\omega}$ be a $\Sigma^0_k$ set. There exists $g \in \mathrm{HS}(A)$ such that $g \leq_{\T} \hj^k(\emptyset)$.
	\end{corollary}
	\begin{proof}
		Consider a chain of nested countable coded $\beta$-models $M^1 \in M^2 \in \dots \in M^k$ as in the statement of Lemma \ref{smallbetamodels}. By Corollary \ref{simpsonreversal}, $\exists g \in M^k$ such that $g \in \mathrm{HS}(A)$. So $g \leq_{\T}M^k \leq_{\T} \hj^k(\emptyset)$.
	\end{proof}
	
	The following Corollary is obtained by combining Corollary \ref{simpsonreversal} with (the first claim of) Lemma \ref{smallbetamodels}.
	
	\begin{corollary}[{\cite[Theorem VI.6.4]{simpson}}]
		For all $k \in \omega$, $\pioo \vdash \bSigma^0_{k}$-$\rt$.
	\end{corollary}
	
	For the converse, we have that $ \rca + \DRT 2 \vdash \pioo$ (see \cite[Lemma VI.6.1]{simpson}). Hence, for any $n$ with $2 \leq n < \omega$, the Galvin-Prikry theorem for $\bDelta^0_n$ sets is exactly as strong as the Galvin-Prikry theorem for $\bDelta^0_2$ sets, and this strength is captured by $\pioo$ (again see \cite[Theorem VI.6.4]{simpson}).
	
	\subsection{Computable analysis}
	
	We now introduce tools and notions from computable analysis which we use in the rest of the paper.
	Our sources are \cite{Weihrauch} and \cite{BrattkaMeas}.
	
	Given a set $X$, a \emph{representation} on $X$ is a partial surjective function $\delta \subc \baire \rightarrow X$, and the pair $(X,\delta)$ is called a \emph{represented space}.
	For any $x \in X$, elements of $\delta^{-1}({x})$ are called ($\delta$-) \emph{names} or \emph{codes} for $X$.
	If $\delta$ and $\delta'$ are two representations of the same set $X$, we say that $\delta$ and $\delta'$ are \emph{computably equivalent} if there are computable functionals $h,h'$ such that $\delta(h(x))=\delta'(x)$ for every $x \in \dom(\delta')$ and $\delta'(h'(x))=\delta(x)$ for every $x \in \dom(\delta)$. We think of computably equivalent representation as if they were essentially the same for the purposes of Weihrauch classification.
	If $f \subc (X,\delta_X) \multif (Y,\delta_Y)$ is a partial multi-valued function between represented spaces, we say that $F \subc \baire \rightarrow \baire$ is a \emph{realizer} of $f$ (denoted as $F \vdash f$) if, for every $x \in \delta_X^{-1}(\dom(f))$, $x \in \dom(F)$, $F(x) \in \dom(\delta_Y)$ and $\delta_Y(F(x)) \in f(\delta_X(x))$. 
	\begin{definition}\label{defweihrauch}
		Let $f \subc (X,\delta_X) \multif (Y,\delta_Y)$ and $g \subc (U,\delta_U) \multif (V,\delta_V)$ be partial multi-valued functions between represented spaces and let $\Phi, \Psi \subc \baire \rightarrow \baire$ be functions such that, for every $F \vdash f$, the function $x \mapsto \Psi(\langle x, F(\Phi(x)) \rangle)$ is a realizer of $g$.
		We say that:
		\begin{itemize}
			\item $g \leq_{\W} f$ ($g$ is \emph{Weihrauch reducible to} $f$) if $\Phi$ and $\Psi$ are computable,
			\item $g \leq^{\ari}_{\W} f$ ($g$ is \emph{arithmetically Weihrauch reducible to} $f$) if $\Phi$ and $\Psi$ are arithmetically definable.
		\end{itemize}
		If $\Phi$ and $\Psi$ are computable (arithmetical) and satisfy the stronger requirement that $x \mapsto \Psi(F(\Phi(x)))$ is a realizer of $g$ for every $F \vdash f$, then we say that $g$ is \emph{strongly} (arithmetically) Weihrauch reducible to $f$, denoted $g \leq_{\mathrm{sW}} f$ ($g \leq^{\ari}_{\mathrm{sW}} f$).
	\end{definition}
	All notions of Weihrauch reducibility give rise to preorders, and quotienting under the relation of bi-reducibility gives rise to equivalence classes called (\emph{strong, arithmetic}) \emph{Weihrauch degrees}.
	The Weihrauch degree of a function is a measure of its uniform computational content.
	
	We now recall some operations on Weihrauch degrees for future use. First, we introduce a unary operator on partial multivalued functions which is necessary to give a definition of the compositional product.
	
	\begin{definition}
		Let $f \subc (X,\delta_X) \multif (Y,\delta_Y)$ be a partial multi-valued function between represented spaces. Define $\real f \subc \baire \multif \baire$ as $\real{f}(p)= \delta_Y^{-1} \circ f \circ \delta_X(p)$.
	\end{definition}
	
	The partial multi-valued function $\real{f}$ is the (set theoretic) union of the realizers of $f$, and it is strongly Weihrauch equivalent to $f$.
	
	\begin{definition}\label{Def:operations}
		Let $f \subc (X, \delta_X) \multif (Y, \delta_Y)$ and $g \subc (U, \delta_U) \multif (V, \delta_V)$ be partial multi-valued functions between represented spaces. Define
		\begin{itemize}
			\item $f \times g$, the \emph{parallel product} of $f$ with $g$ as the function $(x,u) \mapsto f(x) \times g(u)$,
			\item $\hat{f} \subc (X^{\omega}, \delta_{X^{\omega}}) \multif (Y^{\omega}, \delta_{Y^{\omega}})$, the \emph{parallelization} of $f$, as $(x_n)_{n \in \omega} \mapsto \{(y_n)_{n \in \omega} : \forall n \, y_n \in f(x_n)\}$ with $\dom(\hat{f})=(\dom(f))^\omega$,
			\item $f \star g$, the \emph{compositional product} of $f$ and $g$, is given by the function $(p,q) \mapsto (\id {{\baire}} \times \real{f}) \circ \Phi_p \circ \real{g}(q)$ with domain $\{(p,q) \in (\baire)^2 : q \in \dom(\real{g}) \land \forall y \in \real{g}(q) \, (y \in \dom(\Phi_p) \land \Phi_p(y) \in \dom(\id {{\baire}} \times \real{f}))\}$ where $\Phi_p$ is the functional computed by a fixed universal Turing machine with oracle $p$.
		\end{itemize}
	\end{definition}
	
	These operators are all degree theoretic. We say that $f$ is \emph{parallelizable} if $f \equiv_{\W} \hat{f}$. It is immediate to see that parallelization is idempotent, so that for every $f$, $\hat{f}$ is parallelizable.
	The compositional product of $f$ and $g$ corresponds, intuitively, to applying $g$, following it with a computable operation, and then applying $f$.
	Notice that $f \star g$ is defined even when $f$ and $g$ are not composable.
    By \cite[Theorem 11.5.2]{Brattka2021} $f \star g \equiv_{\W} \max_{\leq_{\W}}\{h \circ k : h \leq_{\W} f \land k \leq_{\W} g\}$.
    Following conventional use, we denote by $f^{[n]}$ the $n$-fold compositional product of $f$ with itself. We state a simple lemma on the interaction of the compositional product with parallelization.
    
    \begin{lemma}[{\cite[Proposition 4.11 (9)]{algebraicweihrauchdegrees}}]\label{Lem:distributivity}
    	For any $f,g$, we have $\widehat{f \star g} \leq_{\W} \hat{f} \star \hat{g}$.
    \end{lemma}
    
	We refer the reader to \cite{algebraicweihrauchdegrees} and \cite{Brattka2021} for more information on these operations.

    We now introduce a new operator which in some cases corresponds to composition on the strong Weihrauch degrees. We use this operator in Section \ref{upperbounds}.
    
    \begin{definition}
        Let $f$ and $g$ be multi-valued functions between represented spaces, define $f \wstar g$ as $(p,q) \mapsto \real{f} \circ \Phi_p \circ \real{g}(q)$ with domain $\{(p,q) \in (\baire)^2 : \forall y \in \real{g}(q) \, (y \in \dom(\Phi_p) \land \Phi_p(y) \in \dom(\real{f})\}$ (recall that $\Phi_p$ is the functional computed by a fixed universal Turing machine with oracle $p$).
    \end{definition}
    
    The intuition behind $f \wstar g$ is the following: we use $g$, then perform some computation and then necessarily use $f$ on the output of this computation, \emph{forgetting its input}. The strong similarities between $\wstar$ and $\star$ are made precise by the following Lemma, where we use $\pi_2 \colon (\baire)^2 \rightarrow \baire$ to denote the projection on the second component.
    
    \begin{lemma}\label{lem:relwstarstar} Let $H=H' \times \id{{\baire}}$ where $H'$ is a total computable functional such that $\pi_2 \circ \Phi_p=\Phi_{H'(p)}$ for every $p$. Then for every pair of partial multi-valued functions $f$ and $g$ and every $(p,q) \in \dom(f \star g)$, we have $H(p,q) \in \dom(f \wstar g)$ and $f \wstar g (H(p,q))=\pi_2(f \star g(p,q))$.
    
    Moreover, $H$ admits a computable right inverse $K$. The pair $(K, \pi_2)$ witnesses $f \wstar g \leq_{\sW} f \star g$ for every pair $(f,g)$.
    \end{lemma}
    \begin{proof}
    
        For any $f$ and $g$ and any $(p,q) \in \dom(f \star g)$, for every $y \in \real{g}(q)$, we have $\pi_2(\Phi_p(y))\in \dom(\real{f})$, so that $\Phi_{H'(p)}(y) \in \dom(\real{f})$ and hence $H(p,q) \in \dom(f \wstar g)$. Moreover, $f \wstar g (H(p,q))=\pi_2(f \star g(p,q))$.
        
        Now let $K \colon (\baire)^2 \rightarrow (\baire)^2$ be a computable function which maps any pair $(x,y)$ to a pair $(x',y)$ where, for any $w \in \baire$, $\Phi_{x'}(w)= (\Phi_x(w), \Phi_x(w) )$. It is immediate to see that $H \circ K=\id {{(\baire)^2}}$, so $H$ is surjective and $K$ is one of its right inverses. Now let $(p,q) \in \dom(f \wstar g)$, by definition we have $K(p,q)=(p',q)$ for some $p'$ such that $\Phi_{p'}(w)= (\Phi_p(w), \Phi_p(w))$. So
        \[
        f \star g (K(p,q))=(\pi_1(\Phi_{p'}(\real{g}(q))),\real{f}(\pi_2(\Phi_{p'}(\real{g}(q)))))= (\Phi_p(\real{g}(q))),\real{f}((\Phi_p(\real{g}(q)))) 
        \]
        hence $\pi_2((f\star g)(K(p,q)))=(f\wstar g)(p,q)$ and the pair $(K, \pi_2)$ witnesses $f \wstar g \leq_{\sW} f \star g$.
    \end{proof}

    We now list some basic properties of $\wstar$.
    
    \begin{lemma}\label{Lem:propertiesofwstar}
        Let $f_0$, $f_1$, $g_0$, $g_1$ and $h$ be partial multi-valued functions. The following hold:
        \begin{enumerate}
        	\item $(f_1 \wstar g_1) \wstar h \equiv_{\W} f_1 \wstar (g_1 \wstar h)$,
        	\item $(\id {{\baire}} \times f_1) \wstar g_1 \equiv_{\W} f_1 \star g_1$,
        	\item if $g_0 \leq_{\W} g_1$, then $f_1 \wstar g_0 \leq_{\W} f_1 \wstar g_1$.
        	\item if $f_0 \leq_{\sW} f_1$ and $g_0 \leq_{\sW} g_1$, then $f_0 \wstar g_0 \leq_{\sW} f_1 \wstar g_1$,
        \end{enumerate}
    \end{lemma}
    \begin{proof}
    	Items (1) and (2) are immediate.
    	
    	For (3) let $H$ and $K$ witness $g_0 \leq_{\W} g_1$ so that for every $x$, $H(x, \real{g_1}(K(x))) \subseteq \real{g_0}(x)$, and let $H'$ be a computable function such that, for every $p, q, x \in \baire$, $\Phi_{H'(p,q)}(x)=\Phi_p(H(q, x))$. For every $(p,q) \in \dom(f_1 \wstar g_0)$ we have:
    	\[
    	\real{f_1} \circ \Phi_{H'(p,q)} \circ \real{g_1}(K(q)) = \real{f_1} \circ \Phi_{p} (H (q,\real{g_1}(K(q))) \subseteq \real{f_1} \circ \Phi_p \circ \real{g_0}(q)= f_1 \wstar g_0 (p,q).
    	\]
    	Therefore $\id{{\baire}}$ and $(p,q) \mapsto (H'(p,q),K(q))$ witness $f_1 \wstar g_0 \leq_{\W} f_1 \wstar g_1$.
    	
    	For (4), we first show that $f_0 \wstar g_0 \leq_{\sW} f_1 \wstar g_0$. Let $H$, $K$ witness $f_0 \leq_{\sW} f_1$, so that for every $x$, $H(\real{f_1}(K(x))) \subseteq \real{f_0}(x)$.
    	Let $K'$ be a computable function such that $K(\Phi_{p}(x))=\Phi_{K'(p)}(x)$ for all $p,x \in \baire$.
    	For any $(p,q) \in \dom(f_0 \wstar g_0)$, we have:
        \[
    	H \circ \real{f_1}\circ \Phi_{K'(p)} \circ \real{g_0}(q) = H \circ \real{f_1} \circ K \circ \Phi_p \circ \real{g_0} (q) \subseteq \real{f_0} \circ \Phi_p \circ \real{g_0}(q) =f_0 \wstar g_0 (p,q).
    	\]
    	Therefore $H$ and $K' \times \id {{\baire}}$ witness $f_0 \wstar g_0 \leq_{\sW} f_1 \wstar g_0$.
    	
    	Now we show that $f_1 \wstar g_0 \leq_{\sW} f_1 \wstar g_1$. Let $F$, $G$ witness $g_0 \leq_{\sW} g_1$, so that for every $x$, $F(\real{g_1}(G(x))) \subseteq \real{g_0}(x)$. Let $F'$ be a computable function such that for every $p,x \in \baire$, $\Phi_{F'(p)}(x)=\Phi_p(F(x))$. For every $(p,q) \in \dom(f_1 \wstar g_0)$ we have:
    	\[
    	\real{f_1} \circ \Phi_{F'(p)} \circ \real{g_1}(G(q)) = \real{f_1} \circ \Phi_{p} (F (\real{g_1}(G(q))) \subseteq \real{f_1} \circ \Phi_p \circ \real{g_0}(q)= f_1 \wstar g_0 (p,q).
    	\]
    	Therefore $\id{{\baire}}$ and $F' \times G$ witness $f_1 \wstar g_0 \leq_{\sW} f_1 \wstar g_1$, so that by transitivity we obtain $f_0 \wstar g_0 \leq_{\sW} f_1 \wstar g_1$.   	
	\end{proof}

        \begin{remark}\label{rem:wstarequalstar}
            We point out a consequence of Lemma \ref{Lem:propertiesofwstar}: if $f$ is a cylinder, i.e.\ $\id {{\baire}} \times f \equiv_{\sW} f$ and $g$ is any partial multi-valued function, then, using items (4) and (2) we have $f \wstar g \equiv_{\sW} (\id {{\baire}} \times f) \wstar g \equiv_{\W} f \star g$. 
        \end{remark}
 
		\begin{remark}
		It is not true that, for any $f_0,f_1,g$ such that $f_0 \leq_{\W} f_1$, we have $f_0 \wstar g \leq_{\W} f_1 \wstar g$. As an example, we can set $f_0=\id{{\baire}}$, $f_1 \colon \baire \rightarrow \baire$ given by $x \mapsto 0^{\omega}$ and $g$ to be the jump operator. By definition $f_1 \wstar g \colon (\baire)^2 \rightarrow \baire$ is the constant function with value $0^{\omega}$, while by Lemma \ref{Lem:propertiesofwstar}, $f_0 \wstar g = \id{{\baire}} \wstar g \equiv_{\W} \id{{\baire}} \star g \equiv_{\W} g$.
	\end{remark}
	
	\begin{corollary}\label{Cor:strongcompositionalproduct}
		Let $f,g$ be partial multi-valued functions with $g$ a cylinder. Then $\max_{\leq_{\sW}}\{f' \circ g' : f' \leq_{\sW} f \land g' \leq_{\sW} g\}$ exists and it is represented by $f \wstar g$.
	\end{corollary}
	\begin{proof}
		First, notice that, since $g$ is a cylinder, we have $\{f' \circ g' : f' \leq_{\sW} f \land g' \leq_{\sW} g\}=\{f' \circ g' : f' \leq_{\sW} f \land g' \leq_{\W} g\}$.
		Hence $f \wstar g \in \{f' \circ g' : f' \leq_{\sW} f \land g' \leq_{\sW} g\}$ as, if we let $g' \subc (\baire)^2 \rightarrow \baire$ be given by $(p,q) \mapsto \Phi_p(\real{g}(q))$ we have $g' \leq_{\W} g$ and $f \wstar g= f \circ g'$.
		
 Consider now $f_0 \leq_{\sW} f$, $g_0 \leq_{\sW} g$ such that $f_0 \circ g_0$ exists. Let $F_1, F_2$ and $G_1, G_2$ witness the two strong reductions, so that $F_1 \circ \real{f} \circ F_2(q) \subseteq \real{(f_0)}(q)$ for all $q \in \dom(\real{(f_0)})$ and $G_1 \circ \real{g} \circ G_2(q) \subseteq \real{(g_0)}(q)$ for all $q \in \dom(\real{(g_0)})$. For every $q \in \dom(\real{(f_0 \circ g_0)})$, we have:
		\[
		F_1 \circ \real{f} \circ F_2 \circ G_1 \circ \real{g} \circ G_2 (q) \subseteq \real{(f_0 \circ g_0)}(q).
		\]
		Now let $p \in \baire$ be computable such that $\Phi_p=F_2 \circ G_1$ and let $H$ be the computable map defined by $q \mapsto (p,G_2(q))$. We have for every $q \in \dom(\real{(f_0 \circ g_0)})$
\begin{multline*}
F_1 \circ f \wstar g \circ H(q) = F_1 \circ f \wstar g (p,G_2(q)) =\\
= F_1 \circ \real{f} \circ \Phi_p \circ \real{g} \circ G_2(q) = F_1 \circ \real{f} \circ F_2 \circ G_1 \circ \real{g} \circ G_2(q) \subseteq \real{(f_0 \circ g_0)} (q),   
\end{multline*}
		so that $f_0 \circ g_0 \leq_{\sW} f \wstar g$.
	\end{proof}
	
	Note that the definition $f \star_\mathrm{s} g = \max_{\leq_{\sW}}\{f' \circ g' : f' \leq_{\sW} f \land g' \leq_{\sW} g\}$ appears in \cite{UCCCT}, where the authors show that this maximum exists whenever $f$ and $g$ are cylinders (see also \cite[Corollary 11.5.5]{Brattka2021}). Corollary \ref{Cor:strongcompositionalproduct} shows that it suffices to assume that $g$ is a cylinder, and in this case provides an explicit representative of the strong Weihrauch degree.
	
	\begin{remark}
		Note that the characterization of $\max_{\leq_{\sW}}\{f' \circ g' : f' \leq_{\sW} f \land g' \leq_{\sW} g\}$ given in Corollary \ref{Cor:strongcompositionalproduct} may fail when $g$ is not a cylinder. 
  Indeed, consider the two following examples.
  
  Let $f= \id{{\baire}}$ and let $g \colon \baire \rightarrow \baire$ be the constant function with value $0^{\omega}$.
  If $f' \leq_{\sW} f$ and $g' \leq_{\sW} g$ are composable, then $f' \circ g'$ admits a realizer which consists of a constant function with computable value. 
  On the other hand, $f \wstar g$ is the function $(p,q) \mapsto \Phi_p(0^{\omega})$ and hence does not belong to $\{f' \circ g' : f' \leq_{\sW} f \land g' \leq_{\sW} g\}$.
  Note that in this case the set $\{f' \circ g' : f' \leq_{\sW} f \land g' \leq_{\sW} g\}$ does have a maximum with respect to $\leq_{\sW}$, namely the strong Weihrauch degree containing a total constant function with a computable value. However, this degree is not represented by $f \wstar g$.
  
  On the other hand, let $g \colon \baire \rightarrow \baire$ be the constant function with value $0^\omega$, $p_0, p_1 \in \baire$ be incomparable with respect to Turing reducibility, and $f \colon \baire \rightarrow \baire$ map $x$ to $p_0$ when $x(0)$ is even and to $p_1$ when $x(0)$ is odd.
  Now let $f' \leq_{\sW} f$ and $g' \leq_{\sW} g$ be composable.
  We show that $f' \circ g'$ is not the maximum of the set $\{f' \circ g' : f' \leq_{\sW} f \land g' \leq_{\sW} g\}$.
  By definition, $f' \circ g'$ admits a realizer which consists of a constant function with value some $c \leq_{\T} p_i$ for some $i \in \{0,1\}$; without loss of generality, assume $c \leq_{\T} p_0$.
  Now consider the function $g''= \Phi \circ g$ where $\Phi$ is a computable functional with $0^\omega \in \dom(\Phi)$ and such that $\Phi(0^\omega)=1^\omega$. The function $f \circ g'' \colon \baire \rightarrow \baire$ is constant with value $p_1$. Since $p_1 \nleq_{\T} p_0$ it is clear that $f \circ g'' \nleq_{\sW} f' \circ g'$.
  Therefore, for this choice of $f$ and $g$, the set $\{f' \circ g' : f' \leq_{\sW} f \land g' \leq_{\sW} g\}$ does not have a maximum with respect to $\leq_{\sW}$.
	\end{remark}

	We now recall known representations for different spaces of subsets of a given Polish space $X$. 
 If $(X,\tau)$ is a Polish space and $(B_n)_{n \in \omega}$ is an enumeration of one of its bases, we can define a representation $\delta_{\bSigma^0_1(X)}$ of $\bSigma^0_1(X)$ as $\delta_{\bSigma^0_1(X)}(p)=\bigcup_{n \in \omega}B_{p(n)}$.
	We can also define a representation for closed subsets of $X$ by stipulating that $\delta_{\bPi^0_1(X)}(p)=X \setminus \delta_{\bSigma^0_1(X)}(p)$,\footnote{It is well-known that, in the case $X=\baire$, the representation $\delta_{\bPi^0_1}$ defined here is computably equivalent to the representation $\delta_{\tr}$ which associates (the characteristic function of) a tree $T \subseteq \omega^{<\omega}$ to the closed set $[T] \subseteq \baire$. The same is true for the space $\ramsey$, the only difference being that we consider trees $T \subseteq [\omega]^{<\omega}$.} and exploit these definitions to cover all finite Borel levels as follows: for $n\geq1$, we say $\delta_{\bSigma^0_{n+1}(X)}(p)= \bigcup_{i \in \omega}\delta_{\bPi^0_n(X)}(p_i)$ and $\delta_{\bPi^0_{n+1}(X)}(p)=X \setminus \delta_{\bSigma^0_{n+1}(X)}(p)$. For the representations of the ambiguous classes $\bDelta^0_n$ for $n \geq 1$, we let $\delta_{\bDelta^0_n}(p)=A$ if and only if $p=\langle p_1, p_2 \rangle$ and $\delta_{\bSigma^0_n}(p_1)=\delta_{\bPi^0_n}(p_2)=A$.
	
	For Borel subsets of $X$ of transfinite rank, we use Borel codes.
	This is not new in the context of Weihrauch degrees: we essentially follow \cite{comparison2017}, which in turn is based on \cite[3H]{moschovakis}.
	
	\begin{definition}\label{Def:Borelcodes}
		By recursion define $\mathrm{BC}_0=\{p \in \baire: p(0)=0\}$ and, for every $\alpha <\omega_1$, $\mathrm{BC}_{\alpha}=\{1^{\smallfrown}\langle p_0, \dots, p_n, \dots \rangle :  \forall i \,\, p_i \in \mathrm{BC}_{\beta} \text{ for some } \beta < \alpha \}$. Finally set $\mathrm{BC}=\bigcup_{\alpha \in \omega_1}\mathrm{BC}_{\alpha}$ the set of Borel codes.
	\end{definition}
        We recursively define the notion of a \emph{sub-Borel code} as follows: if $p \in \mathrm{BC}_0$, then $\mathrm{sub}(p)=\{p\}$ and if $p=1^{\smallfrown}\langle p_0, \dots, p_n, \dots \rangle$, then $\mathrm{sub}(p)=\{p\} \cup \bigcup_{i \in \omega}\mathrm{sub}(p_i)$. From the definition it is immediate that the relation $\prec$ defined as $a \prec b$ if and only if ``$a \neq b$ and $a$ is a sub-Borel code of $b$'' is well-founded.
	
	\begin{definition}\label{Def:repBorel}
		Let $(X,\tau)$ be a Polish space, $\delta_{\bSigma^0_1(X)}$ a representation of the space of its open subsets obtained above. 
		We define a representation for the set of all Borel subsets of $X$ $\pi_{X} \colon \mathrm{BC} \rightarrow \mathcal{B}(X)$ by recursion:		
		\begin{itemize}
			\item when $p \in \mathrm{BC}_{0}$ we set $\pi_{X}(p)=\delta_{\bSigma^0_1}(p')$, where $p'$ is such that $p=0^{\smallfrown}p'$;
			\item if $p=1^{\smallfrown}\langle p_0, \dots, p_n, \dots \rangle \in \mathrm{BC}$ we set $\pi_{X}(p)= \bigcup_{i \in \omega} \left( X \setminus \pi_{X}(p_i)\right)$.
		\end{itemize}
	\end{definition}
         For every $\alpha <\omega_1$ we denote by $\pi_{\alpha,X}$ the restriction of $\pi_X$ to $\mathrm{BC}_{\alpha}$. Notice that the range of $\pi_{\alpha,X}$ is contained in $\bSigma^0_{1+\alpha}(X)$.
	The fact that each $\pi_{\alpha, X}$ is a representation of $\bSigma^0_{1+\alpha}(X)$ is proven by induction on $\alpha$.
	Similarly to what we did for Borel sets of finite rank, a $\bPi^0_{\alpha}$-name for a set $B$ is given by a $\bSigma^0_{\alpha}$-name for $X \setminus B$, and a $\bDelta^0_{\alpha}$-name for a set $C$ is a pair given by a $\bSigma^0_{\alpha}$-name for $C$ and a $\bPi^0_{\alpha}$-name for $C$. 
	It is immediate that, for all $\alpha < \omega$, the representations $\pi_{\alpha, X}$ and $\delta_{\bSigma^0_{\alpha+1}(X)}$ are computably equivalent.

 If $X$ is computable Polish, i.e.\ presented as a \emph{computable metric space}\footnote{Recall that a computable metric space $(X,d,\alpha)$ is a separable metric space $(X,d)$ with a dense sequence $\alpha \colon \omega \rightarrow X$ such that $d \circ (\alpha \times \alpha) : \omega^2 \rightarrow \mathbb{R}$ is a computable double sequence of real numbers.} with complete metric, then we can choose the basis $(B_n)_{n \in \omega}$ to consist of open balls with center in the dense set and rational radius. 
 In this case the formula $x \in \pi_{\alpha, X}(p)$ is $\Sigma^0_{1+\alpha}$. 
 This in particular holds for Baire space $\baire$, Cantor space $2^\omega$ and Ramsey space $\ramsey$. For these spaces, the natural presentation as computable metric spaces is equivalent to the one given by (restrictions of) the identity on Baire space.\medskip
	
	Now we give the definitions of the partial multi-valued functions we associate to fragments of the Galvin-Prikry theorem.
	The notation is consistent with that of \cite{marconevalenti}.
	
	\begin{definition}\label{ramseyfunctions}
		Let $\Gamma$ be a pointclass among $\{\bSigma^0_{\alpha}, \bPi^0_{\alpha}, \bDelta^0_{\alpha} : \alpha < \omega_1 \}$, we define:
		\begin{itemize}
			\item $\Gamma$-$\rt \colon \Gamma(\ramsey) \multif \ramsey$ as $\Gamma$-$\rt(A)=\mathrm{HS}(A)$,
			\item $\find_{\Gamma} \subc \Gamma(\ramsey) \multif \ramsey$ as $\find_{\Gamma}(A)=\mathrm{HS}(A) \cap A$, with domain $\dom(\find_{\Gamma})=\{A \in \Gamma(\ramsey) : \mathrm{HS}(A) \cap A \neq \emptyset\}$,
			\item $\wfind_{\Gamma} \subc \Gamma(\ramsey) \multif \ramsey$ as $\wfind_{\Gamma}(A)=\mathrm{HS}(A) \cap A$, with domain $\dom(\wfind_{\Gamma})=\{A \in \Gamma(\ramsey) : \mathrm{HS}(A) \subseteq A\}$.
		\end{itemize}
	\end{definition}
	
	We note that for every pointclass $\Gamma$, the function $\wfind_{\Gamma}$ is a restriction of both $\Gamma$-$\rt$ and $\find_{\Gamma}$ to a smaller domain. Therefore the reductions $\wfind_{\Gamma} \leq_{\sW} \Gamma$-$\rt$ and $\wfind_{\Gamma} \leq_{\sW} \find_{\Gamma}$ always hold.
	On the other hand, nothing can be said a priori about the relation between the ``strong'' (one-sided) function $\find_{\Gamma}$ and its two sided version $\Gamma$-$\rt$.
	
	We recall some useful benchmarks in the context of Weihrauch degrees.
	
	\begin{definition}
		Define the following functions:
		\begin{itemize}
			\item $\lim \subc \baire \rightarrow \baire$ defined as $\lim(\langle p_i \rangle_{i \in \omega})=p$ if $p$ is the pointwise limit of the sequence $(p_i)_{i \in \omega}$,
			\item for every $a \in \mathcal{O}$, $\J^{(a)} \colon \baire \rightarrow \baire$, given by $\J^{(a)}(x)=x^{(a)}$,
			\item $\CBaire \subc \bPi^0_1(\baire) \multif \baire$, given by $C \mapsto C$ with domain $\{C \in \bPi^0_1(\baire) : C \neq \emptyset\}$,
			\item  $\UCBaire \subc \bPi^0_1(\baire) \rightarrow \baire$, given by $\{x\} \mapsto x$ with domain $\{C \in \bPi^0_1(\baire) : \exists x \in \baire \,\, C=\{x\}\}$.
		\end{itemize}
	\end{definition}
	
	Note that, if $a,b \in \mathcal{O}$, and $|a|=|b|$, then $x^{(a)}$ is Turing equivalent to $x^{(b)}$ uniformly in $x$ (by relativization of \cite[Chapter II, Lemma 1.2]{sacks2017}), hence $\J^{(a)} \equiv_{\sW} \J^{(b)}$. For this reason we sometimes abuse the notation and write $\J^{(\alpha)}$ in place of $\J^a$ for some $|a|=\alpha$.
	Accordingly, we sometimes denote the function computing the $n$-th jump as $\J^n$.
	We summarize some known properties of the functions just introduced.
	\begin{proposition}\label{Prop:summary}
		The following hold:
		\begin{enumerate}
			\item for every $n \in \omega$, $\lim^{[n]} \equiv_{\sW} \J^n$,
			\item for all $a \in \mathcal{O}$, $\J^{(a)} \leq_{\W} \UCBaire$,
			\item $\UCBaire \star \UCBaire \equiv_{\W} \UCBaire$,
			\item $\CBaire \star \CBaire \equiv_{\W} \CBaire$,
			\item $\CBaire$ is parallelizable.
		\end{enumerate}
	\end{proposition}
	\begin{proof}
		For item (1), see \cite[Theorem 11.6.7]{Brattka2021}. Item (2) follows immediately from \cite[Theorem 3.13]{kmp}. Item (4) is \cite[Corollary 7.6]{closedchoice} and item (3) follows (as noted by the authors of \cite{closedchoice}) by the same proof as item (4). Lastly, item (5) is \cite[Proposition 11.7.49]{Brattka2021}.
	\end{proof}

        \begin{corollary}\label{Cor:choiceabsorbjumps}
         For any $a \in \mathcal{O}$, $\UCBaire \star \J^{(a)} \leq_{\W} \UCBaire$ and $\CBaire \star \J^{(a)} \leq_{\W} \CBaire$. 
	\end{corollary}
	
	We define the represented space $\tr$ as the set of trees on $\omega$, represented via their characteristic functions, and define $\wf \colon \tr \rightarrow 2$ as $\wf(T)=1$ if $[T]=\emptyset$ and $\wf(T)=0$ if $[T] \neq \emptyset$.
	Since every $\Pi^1_1$ subset of $\omega$  can be written as $\{n \in \omega : [T_n] = \emptyset\}$ for some sequence of trees $(T_n)_{n \in \omega}$ (see e.g. \cite[Lemma VI.1.1]{simpson}), the function $\widehat{\wf}$ is the natural analogue of $\pioo$ in the context of Weihrauch degrees (see e.g. \cite{CBWeihrauch}).
 
	We use $\widehat{\wf}$ and its iterates as benchmarks for Weihrauch classification of the functions introduced in Definition \ref{ramseyfunctions}. We start by recalling what is already known.
	
	\begin{proposition}\label{Prop:MV}
		$\wfind_{\bSigma^0_{1}} \equiv_{\W} \UCBaire$, $\wfind_{\bPi^0_{1}} \leq_{\W} \CBaire$ and $\CBaire \equiv^{\ari}_{\W} \wfind_{\bPi^0_{1}}$.
	\end{proposition}
	\begin{proof}
		 See \cite[Theorem 4.5, Proposition 4.9 and Corollary 5.5]{marconevalenti}.
	\end{proof}
	
    \begin{proposition}
        $\find_{\bPi^0_1}, \wfind_{\bPi^0_1}, \wfind_{\bSigma^0_1} <_{\sW} \widehat{\wf}$.        
    \end{proposition} 
    \begin{proof}
        To see that $\find_{\Pi^0_1} <_{\sW} \widehat{\wf}$, note that $\find_{\Pi^0_1} \equiv_{\sW} \CBaire <_{\sW} \widehat{\wf}$ (see \cite[Theorem 4.7]{marconevalenti} and \cite[Proposition 4.3]{CBWeihrauch}). 
        Similarly, by Proposition \ref{Prop:MV} $\wfind_{\bPi^0_1}, \wfind_{\bSigma^0_1} <_{\sW} \widehat{\wf}$.      
    \end{proof}
    
    \begin{proposition}\label{Prop:findabovesrt}
    	$\SRT 1 \leq_{\W} \CBaire \times \SRT 1 \leq_{\W} \find_{\bSigma^0_{1}} $ and $\CBaire \times \SRT 1 \nleq^{\ari}_{\W} \SRT 1$.
    \end{proposition}
    \begin{proof}
    	See the proof of \cite[Theorem 5.14]{marconevalenti}.
    \end{proof}
    
\section{The function $\bmd$}\label{bmd}
    
    In this section we show that $\widehat{\wf}$ is Weihrauch equivalent to the function computing the hyperjump of a set and also to a function, called $\bmd$, which takes $X \in 2^{\omega}$ into a $\beta$-model containing $X$.
    Moreover, we show how to use $\bmd$ to obtain choice functions on some collections of $\Sigma^1_1$ and $\Pi^1_1$ subsets of Baire space.
    These functions play a role in Section \ref{upperbounds} in our computation of upper bounds for the Weihrauch degrees of the functions introduced in Definition \ref{ramseyfunctions}.
    
    \subsection{Trees, $\Pi^0_1$ classes, and the hyperjump}\label{subs:trees}
    
    We begin by introducing tools and notations relative to uniform enumerations of trees and we relate these to uniform enumerations of effectively closed and $\Sigma^1_1$ sets.
    Recall that a formula $\varphi(f)$ in the language of second order arithmetic is $\Pi^0_1$ if there is a computable predicate $\theta(f,n)$ such that $\varphi(f)=\forall n\, \theta(f,n)$, while it is $\Sigma^1_1$ if it has the form $\exists g\, \psi(f,g)$ for some arithmetical $\psi(f,g)$.
    A $\Pi^0_1$ index for $\varphi$ is any code for a Turing machine computing $\theta$.
    
    \begin{theorem}[{Kleene's Normal Form Theorem, \cite[Lemma V.1.4]{simpson}}]\label{knf}
    	For any $\Sigma^1_1$ formula $\varphi(f)$ there exists a computable predicate $\theta(\sigma, \tau)$ such that $\aca$ proves \[\forall f (\varphi(f) \leftrightarrow \exists g\, \forall m\, \theta (g[m], f[m])).\]	
    \end{theorem}
    
    This implies that an enumeration of $\Sigma^1_1$ predicates can be directly extracted from one for $\Pi^0_1$ predicates.
    We focus on the latter to settle the notation for the rest of the paper.
    
    \begin{definition}\label{defuniformtrees}
    	For every $e \in \omega$, $X \in 2^{\omega}$ and $x,y \in \omega$ say that $\{e\}^X_y$ \emph{rejects} $x$ if $\{e\}^X(x)[y]\downarrow 0$.
    	Otherwise $\{e\}^X_y$ \emph{accepts} $x$.
    	For any $e \in \omega$ and $X \in 2^{\omega}$, define the tree $S^X_{e}=\{ \tau \in \omega^{<\omega} : \forall \tau' \sqsubseteq \tau \, ( \{e\}^X_{|\tau|} \text{ accepts } \tau' )\}$.
    \end{definition}
    \begin{lemma}\label{uniformtrees}
    	There is a computable function $s \colon \omega \rightarrow \omega$ such that $\{s(e)\}^X= \chi_{S^X_{e}}$ for all $X \subseteq \omega$.
    \end{lemma}
    \begin{lemma}\label{predicatestotrees}
    	Let $\varphi(f,X)$ be a $\Pi^0_1$ predicate with $\Pi^0_1$ index $e$, $f \in \baire$ and $X \in 2^{\omega}$.
    	Then $\varphi(f,X)$ holds if and only if $f \in [S^X_{e}]$.
    \end{lemma}
    Therefore any $\Pi^0_1$ predicate $\varphi(f)$ with parameter $X \in 2^{\omega}$ can be expressed in the form $f \in [S]$ for some $X$-computable tree $S$: a \emph{tree code} for $\varphi$ is any index $e$ such that $\{e\}^X=\chi_S$ for one such $S$.
    In this terminology Lemmas \ref{uniformtrees} and \ref{predicatestotrees} state that tree codes are computable uniformly from $\Pi^0_1$ indices. In the other direction, from a tree code for $\varphi$, it is easy to obtain a $\Pi^0_1$ index for a formula equivalent to $\varphi$.
    For this reason we use tree codes and $\Pi^0_1$ indices interchangeably.
    Accordingly, the $n$-th $\Pi^0_1(X)$ predicate (for some $n \in \omega$ and $X \in 2^{\omega}$), is the predicate ``being a branch through $S^X_{n}$''. Exploiting Kleene's Normal Form theorem, we say that $n$ is a code for a $\Sigma^1_1$ formula if it is a tree code for a formula $\theta$ as in the statement of Theorem \ref{knf}.
    
    We now consider the hyperjump function. It is well-known that the possible definitions of the hyperjump of a set $X$ (in terms of indices for $X$-computable well-orders, well-founded $X$-computable trees or in terms of ordinal notations) are all many-one equivalent uniformly in $X$. For definiteness we pick one of these definitions.
    
    \begin{definition}
    	Let $\hj \colon 2^{\omega} \rightarrow 2^{\omega}$ be defined by $\hj(X)=\{e \in \omega : [S^X_e]= \emptyset\}$.
    \end{definition} 
    
    The following Lemma is the strong Weihrauch version of \cite[Exercise VII.1.16]{simpson}.
    
    \begin{lemma}\label{Lem:wfquivhj}
    	$\widehat{\wf} \equiv_{\sW} \hj$.
    \end{lemma}
    \begin{proof}
    	To see that $\hj \leq_{\sW} \widehat{\wf}$, let $\Phi$ be the computable function defined as $\Phi(X)=(\chi_{S^X_e})_{e \in \omega}$ for every $X \in 2^{\omega}$.
    	Then $\hj(X)=\widehat{\wf}(\Phi(X))$.
    	For the converse, given $T=(T_n)_{n \in \omega}$, we compute $\hj(T)$.
    	It is clear that there is a computable function $\Psi$ which maps such $T$ to a function $f_T \in \baire$ such that $T_n=S^T_{f_T(n)}$, so that $f_T(n) \in \hj(T)$ if and only if $[T_n] = \emptyset$.
    	Since $X \leq_{\T} \hj(X)$ uniformly in $X$, it follows that $\Psi$ can be used to define $\Psi'$ with $\Psi'(\hj(T))=\widehat{\wf}(T)$ for all $T \in \dom(\widehat{\wf})$.
    \end{proof}
    
    \subsection{Definition and Weihrauch degree of $\bmd$}
    
    \begin{definition}\label{defofbmd}
    	Let $\bmd : 2^\omega \rightarrow 2^\omega \times \omega^{(\omega^{<\omega})}$ be the function mapping $X$ to a pair $(W, h)$ such that:
    	\begin{enumerate}
    		\item $W \subseteq \omega$ is the code for a countable $\beta$-model $\mathcal{M}$ with $X=W_0$,
    		\item $h \colon \omega^{<\omega} \rightarrow \omega$ is such that, letting $\psi(f)$ stand for $n_0$-th $\Pi^0_1(W_{n_1}, \dots, W_{n_k})$ predicate, $h(n_0, n_1, \dots, n_k)=0$ when $\psi(f)$ has no solution, otherwise $h(n_0, n_1, \dots, n_k)=m+1$ for some $m$ such that $\psi(W_m)$ holds.
    	\end{enumerate}
    \end{definition}
    
    A function $h$ as above plays the role of a Skolem function for $\Pi^0_1$ predicates with parameters in $W$. One such function is arithmetically definable, but not computable, using $W$ as a parameter. In the setting of reverse mathematics this means that, assuming the existence of $W$, $\aca$ proves the existence of $h$. In our setting it is natural to ask for $h$ as part of the output of $\bmd$ to make the $\beta$-models more effective.
    
    Notice that $\bmd$ is a \emph{cylinder}, i.e.\ $\bmd \times \id {{\baire}} \leq_{\sW} \bmd$ because we required $X=W_0$. This implies that any function which is Weihrauch reducible to $\bmd$ is also strongly Weihrauch reducible to it.
    
    Our first goal is to show that $\bmd$ is strongly Weihrauch reducible to $\widehat{\wf}$.
    
    \begin{definition}
    	For any $\tau \in \omega^{<\omega}$, $\sigma \in 2^{<\omega}$ and $X \in 2^{\omega}$, let
    	\[
    	T_{X,\sigma,\tau} = \{\nu \in \omega^{<\omega}:\nu \sqsubseteq \tau \lor [\tau \sqsubseteq \nu 
    	\land \forall e<|\sigma| (\sigma(e)=1 \rightarrow \nu_e \in S^{X\oplus \bigoplus_{i < e }\nu_i}_e ) ]\}\]
    	where $\nu_e$ stands for the sequence $i \mapsto \nu (\langle i,e \rangle)$.
    \end{definition}
    
    The idea is that $T_{X,\sigma,\tau}$ contains all sequences $\nu$ which are compatible with $\tau$ and such that, when $\sigma(e)=1$, the $e$-th computable function with oracle $X \oplus \bigoplus_{i < e} \nu_i$ does not reject $\nu_e$.
    If $g \in [T_{X,\sigma,\tau}]$ then, for all $e$ such that $\sigma(e)=1$, we have that $g_e \in [S^{X \oplus \bigoplus_{i <e}g_i}_{e}]$, or equivalently the $e$-th $\Pi^0_1(X \oplus \bigoplus_{i < e}g_i)$-predicate $\varphi$ holds of $g_e$.
    
    The following is an analogue of Lemma \ref{uniformtrees}.
    
    \begin{lemma}\label{preprocessing}
    	There is a total computable function $\Phi \colon \omega^{<\omega}\times 2^{<\omega} \rightarrow \omega$ such that, for every $\sigma,\tau$ and $X$, $\{\Phi(\sigma,\tau)\}^X$ is the characteristic function of $T_{X,\sigma,\tau}$.
    \end{lemma}
    
    We are now ready to show that $\bmd$ is strongly Weihrauch reducible to $\widehat{\wf}$. This mainly consists in a rewriting of \cite[Lemma VII.2.9]{simpson} in the language of Weihrauch degrees, and follows the same proof.
    
    \begin{lemma}\label{wfcomputesbmod}
    	$\bmd \leq_{\sW} \widehat{\wf}$.
    \end{lemma}
    \begin{proof}
    	By Lemma \ref{preprocessing} there exists a computable function $\Phi$ mapping $X \in 2^{\omega}$ to the sequence $(T_{X,\sigma,\tau})_{\sigma \in 2^{<\omega}, \tau \in \omega^{<\omega}}$.
    	Let $G(X)$ be the set such that $\chi_{G(X)}=1-\widehat{\wf}((T_{X,\sigma,\tau})_{\sigma, \tau})$.
    	Note that $(\sigma,\tau) \in G(X)$ if and only if $T_{X, \sigma, \tau}$ is ill-founded. If this holds, then there is some $j \in \omega$ such that $T_{X, \sigma, \tau^{\smallfrown}j}$ is also ill-founded (so $(\sigma, \tau^{\smallfrown}j) \in G(X)$).
    	Moreover $T_{X,\sigma^{\smallfrown}0,\tau}$ is always ill-founded when $(\sigma, \tau)\in G(X)$.
    	This implies that every pair of finite sequences in $G(X)$ can be extended to another pair in $G(X)$ where both strings are strictly longer.
    	
    	Given $G(X)$ we compute $f \in \baire$ and $s \in 2^{\omega}$ so that for every $n$ we have $(s \upto n, f \upto n) \in G(X)$ and $s$ is lexicographically maximal among the functions $s'$ which satisfy the above condition paired with $f$. Given $f \upto n$ and $s \upto n$, define $s(n)=1$ if $((s \upto n )^{\smallfrown}1, f\upto n) \in G(X)$, $s(n)=0$ otherwise, then define $f(n)=\min\{j : (s \upto (n+1), (f\upto n) ^{\smallfrown}j) \in G(X)\}$.
    	Let $W=\{n  \in \omega : f(n)=1\}$: the proof that $W$ is a code for a countable $\beta$-model containing $X$ is the same as in \cite{simpson}. Notice that there is a fixed index $i$ such that $X=W_i$ regardless of the starting set $X$ and that the recursive construction of $s$ and $f$ (and so of $W$) is uniform in $G(X)$.
    	
    	In the remainder of this proof we specify how to compute the function $h$ in the definition of $\bmd$ and we show that the construction does what we claim.
    	This is implicit in Simpson's proof.
    	
    	First notice that if $s(n)=1$, then, since by construction $T_{X,\sigma, \tau}$ is ill-founded for every $\sigma \sqsubset s$ and every $\tau \sqsubset f$, it follows that $f \in [T_{X, s \upto (n+1), \langle \rangle}]$ which in turn implies that $f_n \in \left[T^{X \oplus \bigoplus_{i < n}f_i}_{n}\right]$.
    	In particular this means that the $n$-th $\Pi^0_1(X \oplus \bigoplus_{i < n}f_i)$ predicate holds of $f_n$.
    	Conversely, if $s(n)=0$, then $T_{X, s\upto n^{\smallfrown} 1, f \upto n}$ is well-founded while $T_{s\upto n, f \upto n, X}$ is ill-founded.
    	This implies that there is no extension $f'$ of $f \upto n$ such that $T^{X \oplus \bigoplus_{i < n}f'_i}_{n}$ has a branch.
    	In particular, this holds for $f$ itself, hence $T^{X \oplus \bigoplus_{i < n}f_n}_{n}$ has no branches (so the $n$-th $\Pi^0_1(X \oplus \bigoplus_{i < n}f_i)$-predicate has no solution).
    	
    	We now show how to compute $h$.
    	It is clear that there is a computable function $h' \colon \omega^{<\omega} \rightarrow \omega$ such that for every $(n_0, n_1, \dots, n_k) \in \omega^{<\omega}$, 
    	$S^{X \oplus \bigoplus_{1 \leq i \leq k} W_{n_i}}_{n_0} = S^{X \oplus \bigoplus_{j < h'(n_0, n_1, \dots, n_k)} f_j}_{h'(n_0, n_1, \dots, n_k)}$.
     Intuitively $h'$ computes an index for the program which simulates the program of index $n_0$ using only a specifically selected portion of its oracle (this portion is determined by the indices $n_1, \dots, n_k$).
    	Therefore, we obtain that the $n_0$-th $\Pi^0_1(X \oplus \bigoplus_{1 \leq i \leq k} W_{n_i})$ predicate has a solution if and only if $s(h'(n_0, \dots, n_k))=1$, in which case $f_{h'(n_0, \dots, n_k)}$ is itself a solution.
    	 
    	Now for every $e$ consider the predicate $\psi_e(g)$ given by $g=\chi_{\mathrm{graph}(f_{e})}$: this predicate is $\Pi^0_1(f_e)$, moreover there is a computable function $c$ such that $\psi_e(g)$ is the $c(e)$-th $\Pi^0_1(X \oplus \bigoplus_{i<c(e)}f_i)$ predicate.
    	For every $e$, $\psi_e(g)$ always has a solution in $\mathcal{P}(\omega)$, hence, by construction of $f$, it follows that $f_{c(e)}$ is a solution to $\psi_e(g)$. In other words $f_{c(e)}=\chi_{\mathrm{graph}(f_e)}$ so, by definition of $W$, $\chi_{W_{c(e)}}=\chi_{\mathrm{graph}(f_{e})}$.
    	Hence we can set $h(n_0, \dots, n_k)=0$ if $s(h'(n_0, \dots, n_k))=0$ and $h(n_0, \dots, n_k)=c(h'(n_0, \dots, n_k))+1$ if $s(h'(n_0, \dots, n_k))=1$.
    	As the entire construction is uniform in the input $X$, we can define $\Psi$ to be the computable procedure which, on input $\widehat{\wf}((T_{\sigma,\tau,X})_{\sigma, \tau})$, outputs the pair $(W, h)$.\footnote{Technically we should also stipulate that $\Psi$ performs a computable rearrangement so that $X=W_0$, rather than $X=W_i$, and rearranges the function $h$ accordingly.}
    	Then $\Phi$ and $\Psi$ witness a strong Weihrauch reduction from $\bmd$ to $\widehat{\wf}$.
    \end{proof}
    
    We conclude the section by proving the converse of Proposition \ref{wfcomputesbmod}, thus characterizing the strong Weihrauch degree of $\bmd$ (for the reverse mathematics version of this equivalence, see Theorem VII.2.10 in \cite{simpson}).
    
    \begin{lemma}\label{bmdcomputeswf}
    	$\widehat{\wf} \leq_{\sW} \bmd$.
    \end{lemma}
    \begin{proof}
    	Given $T=(T_n)_{n \in \omega}$ a sequence of trees, let $(W, h) \in \bmd(T)$. By definition of $\bmd$, $T=W_0$, and it is clear that we have a computable function $c \colon \omega \rightarrow \omega$ such that the $c(n)$-th $\Pi^0_1(T)$ predicate is $f \in [T_n]$. Hence we can compute $\widehat{\wf}((T_n)_{n \in \omega})$ because $\widehat{\wf}((T_n)_{n \in \omega})(n)=0$ if $h(c(n))=0$ and $\widehat{\wf}((T_n)_{n \in \omega})(n)=1$ if $h(c(n)) \neq 0$.
    \end{proof}
    Thus we obtain: 
    \begin{proposition}\label{degreeofbmd}
    	$\bmd \equiv_{\sW} \widehat{\wf}$.
    \end{proposition}
    
    \subsection{Computing choice functions from $\bmd$}
    
    We now show how, from $(W, h) \in \bmd(X)$, we can compute (on the level of the indices) a choice function defined on all $\bPi^0_1$ subsets of Baire space which have a code in $W$.
    Then we exploit this to obtain choice functions for the $\bSigma^1_1$ and $\bPi^1_1$ sets which have a code in $W$.
    
    Note that, having access to $W$, the enumeration of $\Pi^0_1$ predicates with parameters used in Definition \ref{defofbmd} yields an enumeration of the closed subsets of Baire space coded in $W$. 
    
    \begin{definition}
    	For any $X \in 2^{\omega}$, $(W,h) \in \bmd(X)$ and $s= (n_0, n_1, \dots, n_k) \in \omega^{<\omega}$, let $C_{s,W}=\{x \in \baire: \psi_s(x)\}$ where $\psi_s(x)$ is the $n_0$-th $\Pi^0_1(W_{n_1}, \dots, W_{n_k})$ predicate.
    \end{definition}
    
    \begin{remark}
    	The sequence $(C_{s,W})_{s \in \omega^{<\omega}}$ is the sequence of all closed subsets of $\baire$ which have a code in $W$. 
    	Moreover, since the statement $\exists x \psi_s(x)$ is $\Sigma^1_1$ with parameters in $W$ we have that if $C_s \neq \emptyset$, then there exists some $Y$ in $W$ which belongs to $C_s$.
    	By definition of $h$ we have that for every $s=(n_0, n_1 \dots, n_k)$, $C_{s,W}=\emptyset$ iff $h(n_0,n_1, \dots, n_k)=0$. If $h(n_0,n_1, \dots, n_k)=m+1$, then $W_m \in C_{s,W}$.
    \end{remark}
    
    Now we define an enumeration for the analytic subsets of $\baire$ which have a code in the countable coded $\beta$-models we are considering and then show how we can obtain choice functions for them.
    
    \begin{definition}\label{Def:analytic}
    	Let $X \in 2^{\omega}$, $(W,h) \in \bmd(X)$.
    	For any $s \in \omega^{<\omega}$ define $A_{s,W}=\{f \in \baire: \exists g \in \baire \, \, f \oplus g \in C_{s,W}\}$.
    \end{definition}
    
    As above, the sequence $(A_{s,W})_{s \in \omega^{<\omega}}$ contains all analytic subsets of Baire space which have a code in $W$.
    By $\Sigma^1_1$-correctness of $W$ we have that, for every $s \in \omega^{<\omega}$, if $A_{s,W} \neq \emptyset$ then there exists $f \in \baire$ in $W$ such that $f \in A_{s,W}$.
    
    \begin{lemma}\label{computetruthofsigma11}
    	There is a computable function $c$ such that, for any $X \in 2^\omega$, $(W, h) \in \bmd(X)$ and $s \in \omega^{<\omega}$, $W_j \in A_{s,W}$ if and only if $h(c(j,s)) \neq 0$.
    \end{lemma}
    \begin{proof}
    	By $\Sigma^1_1$-correctness of $W$ we have that $W_j \in A_{s,W}$ if and only if there is some $k \in \omega$ such that $W_j \oplus W_k \in C_{s,W}$.
    	Define a computable function $c$ such that for every $j \in \omega$ and $s \in \omega^{<\omega}$, $A=\{g \in \baire : W_j \oplus g \in C_{s,W}\}=C_{c(j,s),W}$.
    	We have $W_j \in A_{s,W}$ if and only if $C_{c(j,s),W} \neq \emptyset$ if and only if $h(c(j,s)) \neq 0$.
    \end{proof}
    
    \begin{corollary}\label{computewitnessesofsigma11}
    	Let $X \in 2^\omega$ and let $(W, h) \in \bmd(X)$.
    	There exists an $h$-computable partial function $\achoice_h \subc \omega^{<\omega} \rightarrow \omega$ such that, for every $s \in \omega^{<\omega}$, $\achoice_h(s) \downarrow$ if and only if $A_{s,W} \neq \emptyset$, and if this is the case and $\achoice_h(s)=m$, then $W_m \in A_{s,W}$.
    	Moreover, there is a computable functional $\Phi$ such that, for every $X \in 2^\omega$ and every $(W,h) \in \bmd(X)$, $\Phi(h)=\achoice_h$.
    \end{corollary}
    \begin{proof}
    	Let $\achoice_h(s)= \min\{j \in \omega : h((c(j,s))) \neq 0\}$.
    	By Lemma \ref{computetruthofsigma11} (together with $\Sigma^1_1$ correctness of $W$), $\achoice_h(s) \downarrow$ if and only if $A_{s,W} \neq \emptyset$ and in this case, if $\achoice_h(s)=m$, then $W_m \in A_{s,W}$.
    \end{proof}
    
    We obtain a similar result for coanalytic sets which have a code in $W$. 
    However in this case nonempty sets might have no elements in $W$ (as being a nonempty $\bPi^1_1$ set is a $\bSigma^1_2$ property).
    
    \begin{corollary}\label{computewitnessesofpi11}
    	Let $X \in 2^\omega$ and let $(W, h) \in \bmd(X)$.
    	There exists an $h$-computable partial function $\cachoice_h \subc \omega^{<\omega} \rightarrow \omega$ such that, for every $s \in \omega^{<\omega}$, $\cachoice_h(s) \downarrow$ if and only if there exists $f$ in $W$ such that $f \notin A_{s,W}$ and, if this is the case and $\cachoice_h(s)=m$, then $W_m \notin A_{s,W}$. Moreover, there is a computable functional $\Phi$ such that, for every $X \in 2^\omega$ and every $(W,h) \in \bmd(X)$, $\Phi(h)=\cachoice_h$.
    \end{corollary}
    \begin{proof}
    	It suffices to define $\cachoice_h(s)=\min\{j \in \omega : h(c(j,s))=0\}$.
    \end{proof}
    
    \begin{remark}\label{turingreduction}
    	Corollary \ref{computewitnessesofsigma11} implies that if $e$ is an index for a Turing machine and $i$ an index such that $\{e\}^{W_i}=Y$ for some $Y \subseteq \omega$, then we can compute an index $j$ with $Y=W_j$; in other words $h$ allows us to perform Turing reductions in $W$ at the level of the indices.
    \end{remark}

 \section{Degree theoretic analysis of the Galvin-Prikry theorem}\label{secseparation}
	
	\noindent In this section we separate fragments of the Galvin-Prikry theorem on the basis of the Turing degrees of the homogeneous sets they assert the existence of.
	The key to obtaining most of our results is the formula $\varphi$ of the next subsection.
	Using $\varphi$ we show that if $\mathcal{I}$ is a jump ideal such that for every (lightface) $A \in \Delta^0_{n+1}(\ramsey)$ there exists $f \in \ramsey \cap \mathcal{I} \cap \mathrm{HS}(A)$, then $\hj^n(\emptyset) \in \mathcal{I}$.
	We also show that, not surprisingly, the converse does not hold: for every $n \in \omega$ there exists a (lightface) $\Delta^0_{n+1}$ set $B_n \subseteq \ramsey$ and a jump ideal $\mathcal{I}_n$ such that $\hj^n(\emptyset) \in \mathcal{I}_n$ yet no set $f \in \ramsey \cap \mathcal{I}_n$ is homogeneous for $B_n$.
	Further, we use an idea from \cite{everyhigherdegree} to extend our results to the transfinite levels of the Borel hierarchy: for any $\omega \leq \alpha <\omega^{\ck}_1$, if $\mathcal{I}$ is a $\mathsf{HYP}$-ideal such that for every (lightface) $A \in \Delta^0_{\alpha}(\ramsey)$ there exists $f \in \mathcal{I} \cap \mathrm{HS}(A)$, then $\hj^{\alpha}(\emptyset) \in \mathcal{I}$. Lastly, we obtain level by level upper bounds to the Turing degrees of homogeneous sets of partitions of transfinite Borel rank.
	
	\subsection{The formula $\varphi$ and its properties}
	
	We employ variations of a formula $\varphi$ found in \cite[Theorem 2.3]{tanaka}, where it is attributed to Simpson (also see \cite[Lemma VI.6.1]{simpson}).
	We say that a function $f \in \ramsey$ \emph{majorizes} a tree $T$ if for all $m \in \omega$ there exists $\tau \in T$ such that $|\tau|=m$ and $\forall i<m \, (\tau(i)\leq f(i))$ (we denote this as $\tau \leq f[m]$). 
    Note that by K\H{o}nig's Lemma a tree $T$ has a branch if and only if it is majorized by some function $f$.
	Given any $f \in \ramsey$ and $n \in \omega$ we denote by $f^n$ the function $i \mapsto f(n+i)$.
	Let 
 \[
 \varphi(f,g)= \forall h<f(0)(f^{1} \text{ majorizes } S^g_h \leftrightarrow f^{2} \text{ majorizes } S^g_h)
 \]
	where $(S^g_h)_{h \in \omega}$ is the $g$-computable sequence of all $g$-computable trees\footnote{There is a single computable function $s$ such that, for every $n \in \omega$ and every $g \subseteq \omega$, $\{s(n)\}^g=\chi_{S^g_n}$ (cf.\ Lemma \ref{uniformtrees}).}.
    Notice that $\varphi(f,g)$ is $\Delta^0_2$.
    
    We think of $g$ as a parameter which determines the sequence of trees of interest. For every $g \in \ramsey$, let $A_g=\{f \in \ramsey : \varphi(f,g)\}$.
	The following two results are implicit in \cite[Lemma VI.6.1]{simpson}. 
	
	\begin{lemma}\label{nothingavoidsphi}

		For every $g \in \ramsey$, there is no set $h \in \ramsey$ which avoids $A_g$, and this is provable in $\aca$.
	\end{lemma}

    \begin{definition}
        We say that $e \in \omega$ is an \emph{homogeneous-to-hyperjump index} (abbreviated as HtH-index) if $\{e\}^X$ is total for every $X \in 2^{\omega}$ and, for every $g \in \ramsey$ and every $f \in \mathrm{HS}(A_g)$, $\{e\}^{(f \oplus g)'}=\hj(g)$.
    \end{definition}
 
	\begin{lemma}\label{ebar}
		There exists an $\mathrm{HtH}$-index $\overline{e}$. This is provable in $\aca$.
	\end{lemma}
	\begin{proof}
		For completeness we explain how this follows from the proof of \cite[Lemma VI.6.1]{simpson}.
		There it is shown that, given any sequence of trees $T=(T_n)_{n \in \omega}$, if $\varphi_T(f)= \forall n < f(0) (f^1 \text{ majorizes } T_n \leftrightarrow f^2 \text{ majorizes } T_n)$, $A=\{f \in \ramsey: \varphi_T(f)\}$ and $f \in \mathrm{HS}(A)$, then for all $n \in \omega$ we have $[T_n] \neq \emptyset$ if and only if $f^{n+2}$ majorizes $T_n$.
		Therefore, for any $g \in \ramsey$, if $f \in \mathrm{HS}(A_g)$, then, for every $n \in \omega$, $[S^g_n] \neq \emptyset$ (i.e., $n \notin \hj(g)$) if and only if $f^{n+2}$ majorizes $S^g_n$.
		The latter condition is $\Pi^0_1(f \oplus g)$, and in particular there is some computable function $c(n)$ (independent of $f$ and $g$) such that $f^{n+2}$ majorizes $S^g_n$ if and only if $c(n) \notin (f \oplus g)'$.
		Therefore we can say that $\overline{e}$ is the index for a program which, on input $n$, checks whether $c(n)$ belongs to its oracle and outputs $1$ if this is the case, and $0$ otherwise.
	\end{proof}

	\begin{remark}
		The fact that HtH-indices are total implies that, if $\overline{e}$ is an HtH-index, then the predicate of $\sigma$, $n$ and $X$ given by $\sigma \in S_n^{\{\overline{e}\}^X}$ is computable.
	\end{remark}
	Lemma \ref{ebar} immediately implies that:
	
	\begin{corollary}\label{firstlevel}
		Let $g \in \ramsey$ and let $\mathcal{I}$ be a jump ideal such that for every $A \in \Delta^0_2(g)$ there exists some $f \in \mathcal{I} \cap \ramsey$ such that $f \in \mathrm{HS}(A)$.
		Then $\hj(g) \in \mathcal{I}$.
	\end{corollary}

	\subsection{The finite Borel levels}\label{sgeneralcase}
	
	We now obtain a sequence of formulas $(\chi_n)_{n \geq 1}$ such that, if $B_n=\{f \in \ramsey : \chi_n(f)\}$, then $B_n \in \Delta^0_{n+1}(\ramsey)$, $\mathrm{HS}(B_n) \subseteq B_n$, and for every $g \in \mathrm{HS}(B_n)$ we have $\hj^n(\emptyset) \leq_{\ari}  g$.
	
	\begin{definition}\label{def:Dn}
		Fix an $\mathrm{HtH}$-index $\overline{e}$. For any $f \in \ramsey$ we define the sequence of sets $(D_n(f))_{n \in \omega}$ recursively as follows: $D_0(f)=\emptyset$ and $D_{n+1}(f)=\{\overline{e}\}^{(f \oplus D_n(f))'}$.
	\end{definition}
	This sequence is always well defined as $\{\overline{e}\}^X$ is a total $\{0,1\}$-valued function for every $X \subseteq \omega$.
	\begin{lemma}\label{complexityofdn}
		The set $D_n(f)$ is uniformly $f^{(n)}$-computable for every $n \in \omega$ and $f \in \ramsey$, i.e.\ there exists a computable $d \colon \omega \rightarrow \omega$ such that, for every $f \in \ramsey$, $\{d(n)\}^{f^{(n)}}=D_n(f)$.
	\end{lemma}
 
	\begin{definition}\label{chin}
	For every $n$ we define a formula $\chi_n$ recursively as follows: $\chi_0(f)$ is $f=f$ and $\chi_{n+1}(f)= \chi_n(f) \land \varphi(f, D_n(f))$.
    Let moreover $B_n=\{f \in \ramsey : \chi_n(f)\}$.
	\end{definition}
	
	\begin{proposition}\label{generalcase}
		For every $n$ the formula $\chi_n$ is $\Delta^0_{n+1}$, $\mathrm{HS}(B_n) \subseteq B_n$ and, for every $g \in \mathrm{HS}(B_n)$, $\hj^n(\emptyset)=D_{n}(g)$.
	\end{proposition}
	\begin{proof}
		By induction.
        The case $n=0$ is trivial, so we assume that the statement holds for $n$ and show that it holds for $n+1$.
		First, to show that $\chi_{n+1}$ is $\Delta^0_{n+2}$, notice that it suffices to show that $\varphi(f, D_{n}(f))$ is $\Delta^0_{n+2}(f)$.
		To do so, we notice that $\varphi(f, D_{n}(f))$ is
		\[\forall h < f(0) \left(f^1 \text{ majorizes } S_h^{D_n(f)} \leftrightarrow f^2 \text{ majorizes } S_h^{D_n(f)}  \right)\] and ``$f^i \text{ majorizes } S_h^{D_n(f)}$'' is of the form $\forall m\, P(m, D_{n}(f) , f)$ where $P$ is a computable predicate of $m$, $f$ and $D_n(f)$.
		By Lemma \ref{complexityofdn}, $D_{n}(f)$ is $f^{(n)}$-computable uniformly in $f$, so $P(m, D_n(f), f)$ is $\Delta^0_{n+1}(f)$ and ``$f^i \text{ majorizes } S_h^{D_n(f)}$'' is $\Pi^0_{n+1}$.
		So, $\varphi(f, D_n(f))$ is $\Delta^0_{n+2}(f)$.
		
		We now show that $\mathrm{HS}(B_{n+1}) \subseteq B_{n+1}$ and if $g \in \mathrm{HS}(B_{n+1})$, then $\hj^{n+1}(\emptyset)=D_{n+1}(g)$.
		Let $g \in \mathrm{HS}(B_{n+1})$: by the inductive hypothesis no set can avoid $B_n$, so by Corollary \ref{relativization} there exists $h \in [g]^{\omega}$ be such that $h$ lands in $B_n$.
		Again by the inductive hypothesis $D_n(k)=\hj^{n}(\emptyset)$ for all $k \in [h]^{\omega}$, so $A_{\hj^n(\emptyset)} \cap [h]^{\omega}=\{k \in [h]^{\omega} : \varphi(k, D_n(k))\}$ (recall that $A_{\hj^n(\emptyset)} =\{k \in \ramsey : \varphi(k, \hj^{n}(\emptyset))\}$).
		By Lemma \ref{nothingavoidsphi} and Corollary \ref{relativization}, we obtain that there is some $h_0 \in [h]^{\omega}$ which lands in $A_{\hj^n(\emptyset)} \cap [h]^{\omega}$ (so, $h_0$ lands in $\{k \in [h]^{\omega} : \varphi(k, D_n(k))\}$).
		Therefore $h_0$ actually lands in $B_{n+1}$ and, since $h_0 \in [g]^{\omega}$ and $g \in \mathrm{HS}(B_{n+1})$, it follows that $g$ also lands in $B_{n+1}$.
		By the reasoning above we also obtain that $g$ lands in $A_{\hj^n(\emptyset)}$, so by Lemma \ref{ebar} we obtain that $\hj^{n+1}(\emptyset) = \{\overline{e}\}^{(g \oplus \hj^{n}(\emptyset))'}=\{\overline{e}\}^{(g \oplus D_n(g))'}=D_{n+1}(g)$.
	\end{proof}
	
	An immediate consequence of Proposition \ref{generalcase} is: 
	
	\begin{theorem}\label{reversal}
		Let $n \in \omega$ and let $\mathcal{I}$ be a Turing jump ideal such that for every $A \in \Delta^0_{n+1}(\ramsey)$ there exists some $g \in \mathrm{HS}(A) \cap \mathcal{I}$, then $\hj^n(\emptyset) \in \mathcal{I}$.
	\end{theorem}
	
	Exploiting the uniformity of Lemma \ref{complexityofdn} we obtain:
	
	\begin{corollary}\label{uniformreduction}
		There is a computable function $k$ such that $\hj^n(\emptyset) =\{k(n)\}^{g^{(n)}}$ for every $n \geq 1$ and $g \in \mathrm{HS}(B_n)$.
	\end{corollary}
	
	Note that it is possible to relativize all definitions from subsection \ref{sgeneralcase} as follows: given any $X$, we define $D_{0, X}(f)=X$ and $D_{n+1,X}=\{\overline{e}\}^{(f \oplus D_{n,X})'}$.
	These can be used to define the sequence of formulas $(\chi_{n,X})_{n \geq 1}$ and related sets $A_{n,X}$ as before.	
	From these, one can obtain relativized versions of all the results above in a straightforward manner.
	
	Now we show that the converse of Theorem \ref{reversal} does not hold.
	The construction is an easy application of a result of Simpson.
	
	\begin{lemma}[{\cite[Lemma 1]{everyhigherdegree}}] 
		For every $a \in \mathcal{O}$ there exists a $\Delta^0_1$ set $K_a \subseteq \ramsey$ such that $\mathrm{HS}(K_a) \subseteq K_a$ and if $f \in \mathrm{HS}(K_a)$ then $\emptyset^{(a)} \leq_{\T} f$.
	\end{lemma}
	
	\begin{remark}\label{anyjumps}
            Inspection of Simpson's construction shows that for every $a \in \mathcal{O}$ we can build a $\Delta^0_1$ formula $\psi_a(f,X)$ such that for any $X \in 2^\omega$, if $K_{a,X}=\{f \in \ramsey : \psi_a(f,X)\}$, then $\mathrm{HS}(K_{a,X}) \subseteq K_{a,X}$ and if $f \in \mathrm{HS}(K_{a,X})$, then $X^{(a)} \leq_{\T} f \oplus X$. This reduction is uniform in $X$, $f$, and $a$ in the following sense: there is a computable function $h \colon \omega \rightarrow \omega$ such that, if $a \in \mathcal{O}$, then $\{h(a)\}^{f \oplus X}=X^{(a)}$.
        \end{remark}

	If $\hat{a}$ is a notation for $\omega$, the formula $\psi_a(f,X)$ can be exploited in conjunction with the formulas $\chi_n$ to build lightface definable sets whose homogeneous sets compute $(\hj^n(\emptyset))^{(\omega)}$.
	
	\begin{lemma}\label{nonconverse}
		Let $\hat{a} \in \mathcal{O}$ be a notation for $\omega$. For all $n \geq 1$, $C_n=\{f \in \ramsey : \chi_n(f) \land \psi_{\hat{a}}(f, D_{n}(f))\}$.
		Then $C_n$ is $\Delta^0_{n+1}$, $\mathrm{HS}(C_n) \subseteq C_n$ and if $g \in \mathrm{HS}(C_n)$, then $(\hj^{n}(\emptyset))^{(\omega)} \leq_{\T} g^{(n)}$.
	\end{lemma}
	\begin{proof}
		The set $C_n$ is $\Delta^0_{n+1}$ as $\psi_{\hat{a}}$ is $\Delta^0_1$ and $D_n(f)$ is uniformly $f^{(n)}$-computable.
		Reasoning as in the proof of Proposition \ref{generalcase} we obtain that no set $g$ avoids $C_n$, so $\mathrm{HS}(C_n) \subseteq C_n$.
		If $g$ lands in $C_n$, then it lands in $B_n$, so by Proposition \ref{generalcase} $D_n(g)=\hj^n(\emptyset)$.
		Hence $g$ also lands in $H_{\hat{a},\hj^n(\emptyset)}=\{f \in \ramsey : \psi_{\hat{a}}(f,\hj^n(\emptyset))\}$, so $(\hj^n(\emptyset))^{(\omega)} \leq_{\T} g \oplus\hj^n(\emptyset) = g \oplus D_n(g) \leq_{\T} g^{(n)}$.
	\end{proof}

 	For every $n \in \omega$, let $\mathsf{ARITH}(\hj^n(\emptyset))=\{X \in 2^\omega : \exists k \, X \leq_{\T} (\hj^n(\emptyset))^{(k)}\}$ (where we set $\hj^0(\emptyset)=\emptyset$).
    For each $n$, $\mathsf{ARITH}(\hj^n(\emptyset))$ is the smallest Turing jump ideal which contains $\hj^n(\emptyset)$.
	\begin{corollary}\label{arith}
		For every $n \in \omega$ there exists a $\Delta^0_{n+1}$ set $C_n \subseteq \ramsey$ such that for every $g \in \ramsey \cap \mathsf{ARITH}(\hj^n(\emptyset))$, $g \notin \mathrm{HS}(C_n)$.
	\end{corollary}
	\begin{proof}
		Follows immediately from Lemma \ref{anyjumps} (case $n=0$) and Lemma \ref{nonconverse}.
	\end{proof}

   Compare our results with the following classical result of Solovay (\cite[Theorem 3.1]{hypencsets}):

   \begin{theorem}
   There exists an effectively closed set $W$ such that $\mathrm{HS}(W) \subseteq W$ but $\mathrm{HS}(W)$ contains neither $\Sigma^1_1$ nor $\Pi^1_1$ sets.    
   \end{theorem}

   Solovay's result readily relativizes, yielding a $\Pi^0_1$ formula $\eta(f,X)$ such that, for every $X \in 2^{\omega}$, letting $W_X=\{f \in \baire : \eta(f,X)\}$, we have that $\mathrm{HS}(W_X) \subseteq W_X$ but $\mathrm{HS}(W_X)$ contains neither $\Sigma^1_1(X)$ nor $\Pi^1_1(X)$ sets.
   This allows to lift Solovay's theorem to higher levels of the effective Borel hierarchy.
   
        \begin{proposition}\label{prop:relsolovay}
            For every $n \in \omega$, there exists a $\Pi^0_{n+1}$ set $E_n$ such that $\mathrm{HS}(E_n) \subseteq E_n$ and, if $g \in \mathrm{HS}(E_n)$, then $g$ is neither $\Sigma^1_1(\hj^n(\emptyset))$ nor $\Pi^1_1(\hj^n(\emptyset))$.          
        \end{proposition}
        \begin{proof}
            It suffices to set $E_n=\{f \in \ramsey : \chi_n(f) \land \eta(f, D_n(f))\}$. Again reasoning as in the proof of Proposition \ref{generalcase} we obtain that $\mathrm{HS}(E_n) \subseteq E_n$, and therefore if $g \in \mathrm{HS}(E_n)$, then, for all $h \in [g]^{\omega}$, $D_n(h)=\hj^n(\emptyset)$. Consequently any $g$ which lands in $E_n$ also lands in $W_{\hj^n(\emptyset)}$, so by the relativization of \cite[Theorem 3.1]{hypencsets} sketched above, $g$ is neither $\Pi^1_1(\hj^n(\emptyset))$ nor $\Sigma^1_1(\hj^n(\emptyset))$.
        \end{proof}

		Notice that Proposition \ref{prop:relsolovay} is an analogue of Corollary \ref{arith}, showing that for every $n \in \omega$, if we want to have homogeneous sets for all $\Pi^0_{n+1}$ partitions (which is more than asking for homogeneous sets for $\Delta^0_{n+1}$ partitions, as in Corollary \ref{arith}), we must look beyond $\Sigma^1_1(\hj^n(\emptyset)) \cup \Pi^1_1(\hj^n(\emptyset))$ (so, in particular, beyond $\mathsf{HYP}(\hj^n(\emptyset))$). 
        
	\subsection{The transfinite Borel levels}
	
	We extend the results from Section \ref{sgeneralcase} to the transfinite levels of the effective Borel hierarchy, obtaining, for all $\alpha < \omega^{\ck}_1$, lower bounds on the complexity of sets contained in $\mathsf{HYP}$-ideals which contain homogeneous sets for $\Delta^0_{\alpha}$ subsets of $\ramsey$.
	This closely follows the method used in \cite{everyhigherdegree}. We also obtain level by level upper bounds on the minimal complexity of homogeneous sets for Borel partitions given by a set of transfinite rank. This is achieved via an effectivization of the proof given by Galvin and Prikry which we obtain exploiting the tools of Section \ref{bmd}.
	
	First, we show that Proposition \ref{generalcase} can be extended to any computable ordinal.
    This requires extending the families $(D_n(f))_{n \in \omega} $and $(B_n)_{n \in \omega}$ (see Definitions \ref{def:Dn} and \ref{chin}).
	
	\begin{definition}
		Fix an $\mathrm{HtH}$-index $\overline{e}$. For any $a \in \mathcal{O}$ and any $f \in \ramsey$, define $E_a(f)$ as follows: if $a=1$, $E_1(f)=\emptyset$, if $a=2^b$, then $E_a(f)=\{\overline{e}\}^{(f \oplus E_b(f))'}$ and if $a= 3 \cdot 5^e$, then $E_a(f)=\oplus_{n \in \omega}E_{e(n)}(f^1)$.
	\end{definition}
	
	Lemma \ref{complexityofdn} extends straightforwardly:
	
	\begin{lemma}\label{complexityofda}
		For every $a \in \mathcal{O}$, then $E_a(f) \leq_{\T}f^{(a)}$, uniformly in $f$.
	\end{lemma}
	\begin{proof}
		By induction on $a \in \mathcal{O}$: the successor case is as in the proof of Lemma \ref{complexityofdn}.
		The limit case is immediate by definition: the $n$-th column in $f^{(a)}$ uniformly computes the $n$-th column of $E_a(f)$.
	\end{proof}
	
	\begin{definition}\label{def:transfinitelevels}
		For every $a \in \mathcal{O}$, define an effective Borel set $C_{a}$ as follows: if $a=1$, $C_1=B_0$, if $a=2^b$, then $C_a=\{f \in \ramsey : f \in C_b \land \varphi(f, E_b(f))\}$, and if $a=3 \cdot 5^e$, then $C_a=\{ f \in \ramsey : f^1 \in C_{e(f(0))}\}$. 
	\end{definition}
	
	Note that if $|a|=n$, then $C_a=B_n$ and $E_a(f)=D_n(f)$, so the families $(C_a)_{a \in \mathcal{O}}$ and $(E_a(f))_{a \in \mathcal{O}}$ extend respectively $(B_n)_{n \in \omega}$ and $(D_n(f))_{n \in \omega}$. 
	
	\begin{proposition}\label{transfinitecase}
		Let $\alpha < \omega^{\ck}_1$ and let $a \in \mathcal{O}$ be a notation for $\alpha$.
		The set $C_a$ is $\Delta^0_{1+\alpha}$, $\mathrm{HS}(C_a) \subseteq C_a$ and, if $g \in \mathrm{HS}(C_a)$, then $D_a(g)=\hj^a(\emptyset)$.
	\end{proposition}
	\begin{proof}
        By Lemma \ref{complexityofda} and by definition of $C_a$, it follows that the predicate $f \in C_a$ is equivalent (uniformly in $f$) to an $(f^{(a)})$-computable predicate, in other words, there is a $\Delta^0_1$ set $A$ such that $(f^{(a)}) \in A$ if and only if $f \in C_a$ for every $f \in \ramsey$. By \cite[Lemma 7.10 and Proposition 7.12]{dghtt} it follows that $C_a$ is $\Delta^0_{1+\alpha}$.
        Showing that nothing avoids $C_a$ and that $E_a(g)=\hj^a(\emptyset)$ is done similarly to the proof of Proposition \ref{generalcase} in case $a$ is a notation for a successor ordinal.
		When $a$ is a notation for a limit ordinal, both statements are easily derived from the inductive assumption (see \cite[Lemma 1]{everyhigherdegree}).
	\end{proof}
	\begin{theorem}\label{thm:transreversal}
		Let $\alpha < \omega^{\ck}_1$ and let $\mathcal{I}$ be a $\mathsf{HYP}$-ideal which contains homogeneous solutions for all $\Delta^0_{1+\alpha}$ sets.
		Then $\hj^{\alpha}(\emptyset) \in \mathcal{I}$.
	\end{theorem}
	
	Clearly we can relativize the definitions of the families $(E_a)_{a \in \mathcal{O}}$ and $(C_a)_{a \in \mathcal{O}}$ to any set $X \subseteq \omega$, obtaining $(E_a)_{a \in \mathcal{O}_X}$ and $(C_a)_{a \in \mathcal{O}_X}$. With these definitions, Proposition \ref{transfinitecase} and Theorem \ref{thm:transreversal} can be extended to all countable ordinals. \medskip
	
	Finally we show how to extend Corollary \ref{simpsonreversal} to the transfinite case. 
    Indeed, iterating the hyperjump operator in the transfinite allows for computation of homogeneous sets for Borel sets of any rank.
	To this end we show that, for any $X \subseteq \omega$ and any ordinal $\alpha < \omega^X_1$, if $\mathcal{M}$ is an $\omega$-model of $\atr$ which contains $\hj^{\omega^{1+\alpha}}(X)$, then $\mathcal{M} \vDash \Sigma^0_{\omega+\alpha}(X)$-$\mathsf{RT}$, and then apply this to the case when $\mathcal{M}$ is a $\beta$-model.
 	
	Our argument exploits what Tanaka called the Galvin-Prikry Lemma (see \cite[p.\ 87]{tanaka} or \cite[Lemma 9]{galvinprikry}), made effective by the use of the $\cachoice$ functions (cf.\ Corollary \ref{computewitnessesofpi11}).
	
	We now state and prove the Galvin-Prikry Lemma, then we show how we can modify its proof in order to have some control on the Turing degrees of the sets involved.
	
	We use the language of \emph{Mathias conditions}. A Mathias condition is a pair $(\sigma, f)$ where $\sigma \in [\omega]^{<\omega}$, $f \in \ramsey$ and $\forall n \in \ran (\sigma) \forall m \in \ran(f) \, \, n<m$. The set of Mathias conditions is ordered by the relation $\preceq$ where $(\sigma, f) \preceq (\tau, g)$ if $\ran(\tau) \subseteq \ran(\sigma) \subseteq (\ran(\tau) \cup \ran (g))$ and $\ran(f) \subseteq \ran(g)$. If $(\sigma, f) \preceq (\tau, g)$ we say that $(\sigma, f)$ \emph{extends} $(\tau, g)$.
	We also say that a set $A \subseteq \ramsey$ is \emph{completely Ramsey} if for every Mathias condition $(\sigma,f)$ there exists some $g \in [f]^{\omega}$ such that for all $h \in \ramsey$, $\sigma^{\smallfrown}(g \circ h) \in A$. The Galvin-Prikry Theorem is proved in \cite{galvinprikry} by showing that open sets are completely Ramsey, and then exploiting the lemma below to conclude that all Borel sets are completely Ramsey.
	\begin{lemma}[Galvin-Prikry]\label{galvinprikry}
		Let $(A_n)_{n \in \omega} \subseteq \ramsey$ be a sequence of completely Ramsey sets.
		There is some $f \in \ramsey$ such that $(\bigcup_{n \in \omega}A_n) \cap [f]^{\omega}$ is open in the subspace topology of $[f]^{\omega}$.
	\end{lemma}
	\begin{proof}
		We build $f$ recursively as follows: we define a sequence of functions $(f_n)_{n \in \omega}$ in $\ramsey$ and a sequence of finite sets $(\sigma_n)_{n \in \omega}$ in $ [\omega]^{<\omega}$, starting with $f_{0}=\id {\omega}$, $\sigma_0=\langle \rangle$.
		Assume we defined $f_0, \dots , f_n$ and $\sigma_0, \dots ,\sigma_n$ such that for all $i \leq n$: $|\sigma_i|=i$, $(\sigma_i, f^1_i)$ is a Mathias condition with $(\sigma_i, f^1_i)$ extending $(\sigma_j, f^1_j)$ for all $j < i$, and for all $i < n$ $(\sigma_{i+1}, f^1_{i+1})$ is homogeneous for $A_i$, i.e. for all $\sigma \subseteq \sigma_{i+1}$ either $\forall g[ \sigma^{\smallfrown}(f^1_{i+1} \circ g) \in A_i]$ or $\forall g [\sigma^{\smallfrown}(f^1_{i+1} \circ g) \notin A_i]$.
		
		Since $A_{n}$ is completely Ramsey, we obtain a finite sequence of functions $(g_{i,n})_{i \leq 2^n}$, starting with $g_{0,n}=f_n$ as follows: we enumerate all subsets of $\sigma_n$ as $(\tau_i)_{i<2^n}$ and for all $i < 2^n$, we pick $g_{i+1,n} \in [g_{i,n}]^{\omega}$ such that, either for all $h \in \ramsey$ we have $\tau_i^{\smallfrown}(g_{i+1,n} \circ h) \in A_n$ or for all $h \in \ramsey$ we have $\tau_i^{\smallfrown}(g_{i+1,n} \circ h) \notin A_n$.
		Note that 
        \[g_{i+1,n} \in [g_{i,n}]^{\omega} \land \forall h \in \ramsey [(\tau_i^{\smallfrown}(g_{i+1,n} \circ h) \in A_n) \lor (\tau_i^{\smallfrown}(g_{i+1,n} \circ h) \notin A_n)] 
        \]
        can be expressed as $g_{i+1,n}=g_{i,n} \circ h$ for some $h \in \mathrm{HS}(A_{n,{g_i}, \tau_i})$ with $A_{n, g_{i,n}, \tau_i}=\{k \in \ramsey : \tau_i^{\smallfrown}(g_{i,n} \circ k) \in A_{n}\}$.
        
		We define $f_{n+1}=g_{2^n,n}$, $\sigma_{n+1}=\sigma_n^{\smallfrown}f_{n+1}(0)$ and we finally set $\tilde{f}=\bigcup_{n \in \omega} \sigma_n$.
		By construction we have that for every $n \in \omega$ and every finite set $\tau \subseteq \sigma_{n+1}$, either for every $g \in \ramsey$ we have $(\tau^{\smallfrown}(f^1_{n+1}\circ g)) \in A_n$ or for every $g \in \ramsey$ we have $(\tau ^{\smallfrown}(f^1_{n+1} \circ g)) \notin A_n$.
		Therefore $(\bigcup_{n \in \omega}A_n) \cap [\tilde{f}]^{\omega}$ is open in the topology of $[\tilde{f}]^{\omega}$.
	\end{proof}
	
	We remark that the Galvin-Prikry Lemma is used to translate the problem of finding homogeneous sets for arbitrary countable unions of completely Ramsey sets into the easier problem of finding homogeneous sets for open sets.
	This is because, if $f$ and $(A_n)_{n \in \omega}$ are as in Lemma \ref{galvinprikry}, and if $F \colon \ramsey \rightarrow [f]^{\omega}$ is the $f$-computable homeomorphism of Lemma \ref{behaviourofcompositions} (so $F(g)=f \circ g$), then $P = F^{-1}(\bigcup_{n \in \omega} A_n \cap [f]^\omega)$ is open and if $g \in \mathrm{HS}(P)$, then $F(g) \in \mathrm{HS}(\bigcup_{n \in \omega}A_n)$.
	We now set out to show that if $A=\bigcup_{n \in \omega}A_n$ is a Borel set which is $\Sigma^0_{\omega+\alpha}$, we can choose $f$ to be such that $f \leq_{\T} \hj^{\omega^{1+\alpha}}(\emptyset)$.
	
	Now we show that we can recursively pass from a $z$-computable Borel code for a set $A \subseteq X$ (where $X \subseteq \baire$ is a computable Polish space) to a $\Delta^1_1(z)$ code for it. 
 This is essentially (the relativization of) \cite[Exercise 7B.7]{moschovakis} stated in the language of number codes for $z$-effectively Borel and $\Delta^1_1(z)$ sets, and it is similarly proven via effective transfinite recursion. We write the proof explicitly for the reader's convenience. 
 Note that by a $\Delta^1_1(z)$ code for a set $A$ we mean a pair natural numbers $(a,b)$ such that $a$ and $b$ are tree codes for $\Sigma^1_1$ formulas $\varphi_1(f,z)$ and $\varphi_2(f,z)$ such that $f \in A \leftrightarrow \varphi_1(f,z) \leftrightarrow \neg\varphi_2(f,z)$ (cf.\ the tree codes for $\Sigma^1_1$ formulas introduced in Section \ref{subs:trees}). 
 We recall that we denote the set of Borel codes for $\bSigma^{0}_{1+\alpha}$ sets of Definition \ref{Def:Borelcodes} by $\mathrm{BC}_{\alpha}$ and we write $\pi_{\alpha, X} \subc \baire \rightarrow \mathcal{B}(X)$ to denote the representation of Borel sets of $X$ of rank $1+\alpha$ of Definition \ref{Def:repBorel}.
	
	\begin{proposition}\label{translationcodes}
		There is an index $e$ such that, for all $z$, if $d \in \omega$ is such that $\{d\}^z$ is a Borel code for a set $A \subseteq X$, then $\{e\}^z(d)$ is a $\Delta^1_1(z)$ code for $A$. 
	\end{proposition}
	\begin{proof}
		First we prove the version of the proposition without the oracle $z$.
		In the course of this proof, if $\{a\}$ computes a total function, we denote by $\hat{a}$ the sequence $\hat{a}(n)=\{a\}(n)$.
		Consider the set $C=\{a \in \omega : \{a\} \text{ is total } \land \hat{a} \in \mathrm{BC}\}$ (we say that $a \in C$ codes a Borel code). On $C$ consider the relation $\prec$ given by $a \prec b$ if and only if $\hat{a} \neq \hat{b}$ and $\hat{a}$ is a sub-Borel code of $\hat{b}$ (in the sense of Definition \ref{Def:Borelcodes}).
		It is immediate from the definition of Borel codes that $\prec$ is well-founded.
		
		We want to obtain the index $e$ by use effective transfinite recursion on $\prec$ (\cite[Theorem I.3.2]{sacks2017}). To that end, we prove two intermediate results.
		
		\begin{claimn} There is an index $r_1$ such that, for every $z \in \baire$ and every $d \in \omega$, if $\{d\}^z \in \dom(\mathrm{BC}_0)$, then $\{r_1\}^z(d)$ is a $\Delta^1_1(z)$ code for $\pi_{0, X}(\{d\}^z)$.
		\end{claimn}
		
		\begin{proof}[Proof of Claim] This is essentially an effective version of Kleene's Normal Form Theorem (Theorem \ref{knf}) restricted to $\Sigma^0_1$ (resp.\ $\Pi^0_1$) formulas of the form $f \in \pi_{0, X}(\{d\}^z)$ (resp.\ $f \notin \pi_{0, X}(\{d\}^z)$).
		\end{proof}
		
		\begin{claimn} There is an index $r_U$ such that, for all $z$, if $\{d\}^z$ is total and lists a sequence of $\Delta^1_1(z)$ codes for sets $A_n \subseteq X$, then $\{f_U\}^z(d)$ is a $\Delta^1_1(z)$ code for $A = \bigcup_{n \in \omega} (X \setminus A_n)$.
		\end{claimn}
		
		\begin{proof}[Proof of Claim] Intuitively $\{f_U\}^z(d)$ outputs a $\Delta^1_1(z)$ code for the set $A$ by transforming $d$ into a $\Delta^1_1(z)$ code for $B=\{n^{\smallfrown}f \in \baire: f \notin A_n\}$, as then $f \in A \leftrightarrow \exists n\, n^{\smallfrown}f \in B$, and obtaining a $\Delta^1_1(z)$ code for $A$ from one for $B$ is a straightforward application of the algorithm used to prove Kleene's Normal Form theorem.
		\end{proof}
		
		Now we set $\{e\}(d)=\{f_1\}(d)$ if $\hat{d}(0)=0$ while, if $\hat{d}=1^{\smallfrown}\langle p_0, \dots, p_n, \dots \rangle$ we set $\{e\}(d)=\{f_U\}(d')$ where $d'$ is an index for the function $i \mapsto \{e\}(d_i)$ and each $d_i$ is an index such that $\hat{d_i}=p_i$, so that $\hat{d_i}(n)=\hat{d}(1+\langle i,n \rangle)$.
		
		For the relativized version, note that, for any fixed $z$, we would define $C_z$ and $\prec_z$ analogously to $C$ and $\prec$.
        Then the definition of the index working for oracle $z$ uses $\{f_1\}^z$ and $\{f_U\}^z$ but is formally the same as above, if we view the unrelativized version as having oracle $\emptyset$.
        Thus the same index works for all oracles.
	\end{proof}
	
	\begin{proposition}\label{Prop:transfiniteub}
There is a computable functional $\Phi$ such that, for all $z \in \baire$, $\alpha < \omega^z_1$, $d \in \omega$ such that $\{d\}^z$ is a Borel code for a $\Sigma^0_{\omega+\alpha}(z)$ set $B \subseteq \ramsey$, and $a \in \mathcal{O}_z$ such that $|a|_z = \omega^{1+\alpha}$, we have that $f=\Phi(d, z, a, \hj^{a}(z)) \in \ramsey$ and $B \cap [f]^{\omega}$ is open in the subspace topology of $[f]^\omega$.
	\end{proposition}
	\begin{proof}
        Throughout the proof if $a \in \mathcal{O}_z$ for some $z \in \baire$ is a notation for a successor ordinal, we denote by $a^-$ the unique predecessor of $a$ in $\mathcal{O}_z$ (so in particular $|a|_z=|a^-|_z+1$).
        
		It is easy to see that there is a computable function $c$ such that, for every $z$ and every $a \in \mathcal{O}_z$ such that $|a|_z$ is a limit, $c(a) \in \mathcal{O}_z$, $|c(a)|_z=|a|_z$, and if $c(a)=3 \cdot 5^e$, then $|\{e\}^z(n)|_z$ is a successor ordinal for every $n$.
		By the relativized analogue of \cite[Lemma II.1.2]{sacks2017} for the hyperjump, it follows that there is a computable functional $F$ such that, for all $z \in \baire$ and $a \in \mathcal{O}_z$, $F(a,z,\hj^a(z))=(c(a), z, \hj^{c(a)}(z))$.
		
		We can now describe what $\Phi$ does on its input and later check that it satisfies our requirements. The idea is that it tries to emulate the proof of Lemma \ref{galvinprikry}. The inputs to $\Phi$ we are interested in are of the form $(d, z, a, \hj^a(z))$ where $\{d\}^z$ is a Borel code for a set $B \in \Sigma^0_{\omega+\alpha}(z)$, $a \in \mathcal{O}_z$ and $|a|_z=\omega^{1+\alpha}$. Note that by definition $B= \bigcup_{n \in \omega} B_n$ for a sequence $(B_n)_{n \in \omega}$ such that for every $n \in \omega$, $\ramsey \setminus B_n$ is coded by a sub-Borel code of $\{d\}^z$. It is immediate to see that there is a computable function which, on input $d$ as above and $n \in \omega$, outputs an index $d^n$ such that $\{d^n\}^z$ is the sub-Borel code of $\{d\}^z$ which codes $\ramsey \setminus B_n$.
		
		On input $(d, z, a, \hj^a(z))$, we use the functional $F$ to obtain $F(a,z,\hj^a(z))=(c(a), z, \hj^{c(a)}(z))$ with $c(a)=3 \cdot 5^e$ for some $e$. 
		
		Since $\hj^{c(a)}(z)= \bigoplus_{i \in \omega}\hj^{e(i)}(z)$, we refer to $\hj^{e(i)}(z)$ as the $i$-th column of $\hj^{c(a)}(z)$. By definition of $c(a)$, the $i$-th column of $\hj^{c(a)}(z)$ is $\hj(\hj^{\{e\}(i)^-}(z))$. By Lemma \ref{wfcomputesbmod}, there is a computable functional $G$ which, on input $\hj(\hj^{\{e\}(i)^-}(z))$, computes $(M,h) \in \bmd(\hj^{\{e\}(i)^-}(z))$. So we can think that, using $F$ and $G$, we have access to a sequence $((M^i,h_i))_{i \in \omega}$ with $(M^i,h_i) \leq_{\T} \hj^{\{e\}(i)}(z)$ and $(M^i,h_i) \in \bmd(\hj^{(\{e\}(i))^-}(z))$.
		
		To output the first bit of $\Phi(d, z, a, \hj^{a}(z))$, we take $d^0$ and $z$ and, using the computable function of Proposition \ref{translationcodes}, compute a $\Delta^1_1(z)$ index $d$ for $B_0$. Then, exploiting the function $\{e\}$ and the index $d$, in parallel for all $i$, we compute finite strings $t_i$ such that (using the notation of Definition \ref{Def:analytic}), $A_{t_i, M^i}=\{f \in \ramsey : f \notin \mathrm{HS}(B_0)\}$, and start the computations of $\cachoice_{h_i}(t_i)$. 
        This continues indefinitely unless one of the computations terminates. If this ever happens, say $\cachoice_{h_{j_0}}(t_{j_0})\downarrow=m_0$, then $M^{j_0}_{m_0}=f_1 \in \ramsey$ we output $f_1(0)$ and store $j_0$ and $m_0$.
		
		Now to output the $n+1$-th bit of $\Phi(d, z, a, \hj^{a}(z))$ we proceed as follows: following the proof of Lemma \ref{galvinprikry} we need to extract homogeneous subsets $2^n$ times. For the first time, we take $j_n$ and $m_n$ and try to obtain some $h \in \mathrm{HS}(B_{n,g_{0,n},\tau_0})$ (the set $B_{n,g_{0,n},\tau_0}$ is defined from $B_n$, $g_{0,n}$, and $\tau_0$ as in our remark in the proof of Lemma \ref{galvinprikry}) where $g_{0,n}=f_n=M^{j_n}_{m_n}$ and  $\tau_0$ is the first element (with respect to some predetermined computable wellordering of $[\omega]^{<\omega}$) in the enumeration of the subsets of $\sigma_n$, the current content of the output tape.
		
		To that end, we compute a Borel code for $B_{n,g_{0,n}, \tau_0}$, convert it to a $\Delta^1_1(z \oplus f_n)$ code using Proposition \ref{translationcodes}, and compute a sequence of finite strings $t_i$ such that, for all $i \geq j_n$, $A_{t_i,M^i}=\{f \in \ramsey : f \notin \mathrm{HS}(B_{n,g_{0,n},\tau_0})\}$. 
        In parallel, we start the computations of $\cachoice_{h_i}(t_i)$ for every $i \geq j_n$, and continue computing until one of these converges. If the computation with index $k \geq j_n$ converges with value $m'$, then we set $g_{1,n}=M^k_{m'}$, and we proceed to define the following $g_{i,n}$ similarly until we obtain $g_{2^n,n}=f_{n+1}$. Then, we output $f_{n+1}(0)$.
		
		We now show that $\Phi$ computes indefinitely on all valid inputs, and that it produces a correct output.
		
		One thing worth stressing is that, when computing a sequence of finite strings $t_i$ such that $A_{t_i,M^i}=\{f \in \ramsey : f \notin \mathrm{HS}(B_{n,g_{k,n},\tau_k})\}$, it is essential to restrict the process to $i \geq j_n$ (where $j_n \in \omega$ is such that $M^{j_n}_m=g_{k,n}$, and we have the indices $j_n$, $m$ stored from the previous step) for the following reason: if $i < j_n$, we are not even guaranteed that $g_{k,n}$ is in $M^i$,\footnote{Recall that a sequence $t_i$ determines an analytic set $A_{t_i, M^i}$ by specifying a program code $\ell$ together with the indices $a_1, \dots, a_s$ such that $M^i_{a_1}, \dots, M^i_{a_s}$ are the oracles relevant in the definition of the set. Since $g_{k,n}$ may not belong to $M^{i}$ for $i < j_n$, we can be sure that these are correct only when $i \geq j_n$.} so the search for $t_i$ might give an incorrect string. 
        It follows that $\cachoice_{h_i}(t_i)$ could lead to a wrong answer and spoil the entire computation.
		
		On the other hand, if we restrict to $i \geq j_n$, using the function $\{e\}$, $\hj^{c(a)}(z)$ and the indices $j_n$ and $m$ it is possible to compute a sequence of indices $k_i$ such that $M^i_{k_i}=M^{j_n}_m$ for every $i \geq j_n$, and use these indices to correctly compute the sequence of the $t_i$.
		Once we have a sequence of $t_i$ such that $A_{t_i,M^i}=\{f \in \ramsey : f \notin \mathrm{HS}(B_{n,g_{k,n},\tau_k})\}$, Corollary \ref{computewitnessesofpi11} guarantees that, if one of the $\cachoice_{h_i}(t_i)$ converges, then its output is the index for a set which is homogeneous for $B_{n,g_{k,n},\tau_k}$.

		The discussion above shows that $\Phi$ does not make mistakes at the local level.
		We now prove that $\Phi$ computes indefinitely and produces a correct output by induction on $\alpha$.
		
		In the base case $\alpha=0$, we have that $|a|_z=\omega$, $\{d\}^z \in \mathrm{BC}_{\omega}$ codes a set $B \in \Sigma^0_{\omega}(z)$ given by $B=\bigcup_{n \in \omega}B_n$, and the sequence of sub-Borel codes of $\{d\}^z$ codes the sequence sets $\ramsey \setminus B_n$ where each $B_n$ is $\Sigma^0_{s(n)}(z)$ for some $s(n) \in \omega$.
		Since $|c(a)|_z=\omega$, letting $c(a)=3 \cdot 5^e$ and defining $\ell \colon \omega \rightarrow \omega$ as $\ell(n)=|\{e\}(n)|_z$, the sequence $((M^i,h_i))_{i \in \omega}$ of $\beta$-models introduced above is such that $(M^i,h_i) \leq_{\T} \hj^{\ell(i)}(z)$ and $\hj^{\ell(i)-1}(z) \in (M^i,h_i)$.
		
		By Corollary \ref{Cor:degupperbounds} if $q \leq_{\T} \hj^m(z)$ and $B' \in \Sigma^0_k(q)$, then there exists $g \leq_{\T} \hj^{m+k}(z)$ such that $g \in \mathrm{HS}(B')$.
		Notice that, by induction on the step of computation, using the fact that $\omega$ is additively closed, at every step of the computation of $\Phi$ the algorithm searches (inside each of the $\beta$-models $M^i$ for $i$ sufficiently large) for homogeneous sets for subsets of Ramsey space which are $\Sigma^0_k(p)$ for some $k \in \omega$ and some $p \leq_{\T} \hj^{m}(z)$ for some $m \in \omega$. Since the function $\ell$ is cofinal in $\omega$, by the previous observation, these searches must eventually succeed. The correctness of the output follows from the fact that $\Phi$ is emulating the construction in the proof of Lemma \ref{galvinprikry}, and does not make mistakes during the construction.
 
		Now in case $\alpha > 0$, we have that $|a|_z=\omega^{1+\alpha}$, $\{d\}^z \in \mathrm{BC}_{\omega+\alpha}$ codes a set $B \in \Sigma^0_{\omega+\alpha}(z)$ given by $B=\bigcup_{n \in \omega}B_n$, and the sequence of sub-Borel codes of $\{d\}^z$ codes the sequence of sets $\ramsey \setminus B_n$, where each $B_n$ is $\Sigma^0_{\omega+\gamma(n)}(z)$ for some $\gamma(n)< \alpha$. Since $|c(a)|_z=\omega^{1+\alpha}$, letting $c(a)=3 \cdot 5^e$, the sequence $((M^i,h_i))_{i \in \omega}$ of $\beta$-models introduced above is such that $(M^i,h_i) \leq_{\T} \hj^{\{e\}(i)}(z)$ and $\hj^{(\{e\}(i))^-}(z) \in (M^i,h_i)$.
		
		Notice that the inductive assumption implies in particular that for every $\gamma' < \alpha$, $z' \in \baire$ and every $\Sigma^0_{\omega+\gamma'}(z')$ set $B'$, there exists $f' \leq_{\T} \hj^{\omega^{1+\gamma'}}(z')$ such that $B' \cap [f']^{\omega}$ is open in the topology of $[f']^{\omega}$, and therefore if $M$ is a $\beta$-model containing $ \hj^{\omega^{1+\gamma'}}(z')$, then $M$ contains an homogeneous set for $B'$.
		
 		Again, by induction on the step of computation, at all times during the computation of $\Phi$ the algorithm searches (inside each of the $\beta$-models $M^i$ for sufficiently large $i$) for homogeneous sets for subsets of Ramsey space which are $\Sigma^0_{\omega+\gamma}(p)$ for some $p \leq_{\T} \hj^{\gamma'}(z)$, where $\gamma < \alpha$ and $\gamma' < \omega^{1+\alpha}$ (at each step the ordinal $\gamma'$ is bounded by a finite sum of ordinals strictly smaller than $\omega^{1+\alpha}$. Hence $ \gamma' <\omega^{1+\alpha}$ as the latter ordinal is additively closed). Since the function $\ell \colon \omega \rightarrow \omega^{1+\alpha}$ given by $\ell(n)=|\{e\}(n)|_z$ is cofinal in $\omega^{1+\alpha}$, again these searches must eventually succeed. The correctness of $\Phi$ is again a consequence of the fact that $\Phi$ follows (correctly) the construction of the Galvin-Prikry Lemma.
	\end{proof}
	
	Proposition \ref{Prop:transfiniteub} immediately yields: 
	
	\begin{theorem}\label{thm:transub1}
		Let $\alpha < \omega^{\ck}_1$ and let $\mathcal{M}$ be an $\omega$-model of $\atr$ containing $\hj^{\omega^{1+\alpha}}(\emptyset)$. Then $\mathcal{M} \vDash \sRT {\omega+\alpha}$.	
	\end{theorem}
	
	\begin{theorem}\label{thm:transub2}
		Let $\alpha < \omega^{\ck}_1$ and let $A \subseteq \ramsey$ be a $\Sigma^0_{\omega+\alpha}$ set. There exists $g \in \mathrm{HS}(A)$ with $g \leq_{\T} \hj^{\omega^{1+\alpha}+1}(\emptyset)$.
	\end{theorem}
	
	These two results can be seen as level-by-level degree theoretic upper bounds for solutions to instances of the Galvin-Prikry theorem, the first one being less explicit but in some sense sharper than the second.
    Note that these theorems relativize immediately to Borel sets which are not effectively Borel, as is evident from Proposition \ref{Prop:transfiniteub}.
	
\section{Weihrauch classification}\label{upperbounds}
	
	We classify the (arithmetical) Weihrauch degrees of functions related to restrictions of the Galvin-Prikry theorem.
	Important benchmarks in this classifications are $\widehat{\wf}$ and its iterates.
	In light of Lemma \ref{Lem:wfquivhj} and Proposition \ref{wfcomputesbmod}, we use the function $\bmd$, $\hj$ and $\widehat{\wf}$ interchangeably.
	
	\begin{lemma}
		For any $n \in \omega$, $\hj^n \equiv_{\sW} \hj^{[n]}$.
	\end{lemma}
	\begin{proof}
		First we show $\hj \star \hj \leq_{\sW} \hj^2$. This implies $\hj \star \hj \equiv_{\sW} \hj^2$ as the reduction $\hj^2 \leq_{\sW} \hj \star \hj$ is immediate. Since $\hj$ is a cylinder, by \cite[Lemma 3.10]{algebraicweihrauchdegrees}, we know that there is a computable $K$ such that $\hj \star \hj \equiv_{\sW} \hj \circ K \circ \hj$. Now by basic facts about the hyperjump it is clear that there is a computable functional $F$ such that for all $x \in \baire$ $F(\hj(x))=\hj(K(x))$. In particular, setting $x=\hj(y)$, then $F(\hj(\hj(y))=\hj(K(\hj(y))$ for all $y$. In other words, $\id {{\baire}}$ and $F$ witness $\hj \star \hj \leq_{\sW} \hj^2$.
		
		To show the equivalence for all $n \in \omega$, we use induction. So assume $\hj^{[n]} \leq_{\sW} \hj^n$, then $\hj^{[n+1]} \leq_{\W} \hj \star \hj^n$ by monotonicity of $\star$ and hence $\hj^{[n+1]} \leq_{\sW} \hj \star \hj^n$ because $\hj \star \hj^n$ is a cylinder. Now the reduction $\hj \star \hj^n \leq_{\sW} \hj^{n+1}$ holds and its proof is formally identical to the reduction $\hj \star \hj \leq_{\sW} \hj^{2}$, so we can conclude $\hj^{[n+1]} \leq_{\sW} \hj^{n+1}$. Again the reduction in the other direction is trivial.
	\end{proof}
	For this reason, we often use $\hj^n$ instead of $\widehat{\wf}^{[n]}$.
	
	We start by showing that, in analogy with the results from reverse mathematics, the function $\bmd$ can be exploited to compute infinite homogeneous sets for Borel subsets of $\ramsey$ of finite Borel rank.
    In the following we apply $\bmd$ to codes for Borel sets which, according to our definitions, are elements of $\baire$.
    This is done by implicitly identifying elements of $\baire$ with their graphs.

	\begin{lemma}\label{typeonecomputable}
		Let $X \in 2^{\omega}$, $(W, h) \in \bmd(X)$ and denote by $A_x \subseteq \ramsey$ the open set with $\delta_{\bSigma^0_1}$-name $W_x$ (if $W_x$ is a total function). 
		There is a $h$-computable function $\homo^1_h \subc \omega \rightarrow \omega \times 2$ such that for every $x \in \omega$, if $W_x$ is a total function, then $\homo^1_h(x)\downarrow$. In this case if $\homo^1_h(x)=(y,i)$, then $W_y$ lands in $A_x$ if $i=0$ and $W_y$ avoids $A_x$ if $i=1$.
		Moreover, there is a computable functional $\Phi$ such that, for every $X \in 2^\omega$ and every $(W, h) \in \bmd(X)$, $\Phi(h)=\homo^1_h$.
	\end{lemma}
	\begin{proof}
		Consider the $\Pi^1_1$ formulas $\mathrm{land}_1(g,x)= \forall h \, (g \circ h \in A_x)$ and $\mathrm{avoid}_1(g,x)= \forall h \, (g \circ h \notin A_x)$.
		By Kleene's Normal Form Theorem the formulas $\neg\mathrm{land}_1(g,x)$ and $\neg\mathrm{avoid}_1(g,x)$ can be written as, respectively, $\exists f \, f \oplus g \in C_{s_1(x),W}$ and $\exists f \, f \oplus g \in C_{s_2(x),W}$ for some computable functions $s_1, s_2 \colon \omega \rightarrow \omega^{<\omega}$.
		Since $\SRT1$ is true in $W$ (by Theorem \ref{openramsey} and because $W$ is a model of $\atr$) and $W$ is correct on $\Pi^1_1$ statements about its sets, at least one of $\mathrm{land}_1(g, x)$ and $\mathrm{avoid}_1(g, x)$ has a solution in $W$.
		Therefore we can compute $\homo^1_h(x)$ by running $\cachoice_h(s_1(x))$ and $\cachoice_h(s_2(x))$ in parallel to find the least index $y$ such that $\mathrm{land}_1(W_y,x)$ or $\mathrm{avoid}_1(W_y,x)$.
		We then output the pair $(y,i)$ where $i=0$ or $i=1$, according to which of the two computations halted.
		The existence of a uniform $\Phi$ follows directly from the uniformity of $h \mapsto \cachoice_h$ (Corollary \ref{computewitnessesofpi11}).
	\end{proof}
 
	\begin{corollary}\label{opencase}
		$\SRT1 \leq_{\sW} \widehat{\wf}$. 
	\end{corollary}
	\begin{proof}
		By Proposition \ref{degreeofbmd}, it suffices to show that $\SRT1 \leq_{\sW} \bmd$.
        This is an immediate consequence of Lemma \ref{typeonecomputable}: given $X$ coding an open set, let $(W, h) \in \bmd(X)$, and compute $\homo^1_h(0)$ ($X=W_0$ by definition of $\bmd$) to obtain an index $y$ for an homogeneous set $W_y$.
		This establishes a strong Weihrauch reduction thanks to the uniformity of $h \mapsto \homo^1_h$.
	\end{proof}

    In Proposition we show \ref{Prop:hjstrictlyabove} that $\widehat{\wf} \nleq_{\W} \SRT1$.
    
    Notice that the proof of Lemma \ref{typeonecomputable} uses only that $\SRT1$ holds in $W$ as this guarantees the termination of the program looking for homogeneous sets.
	This observation suffices to obtain upper bounds for versions the Galvin-Prikry theorem for $\bSigma^0_k$ sets in terms of iterations of $\widehat{\wf}$.
	
	\begin{lemma}\label{firstupperbound}
		Let $k \geq 1$.
		We have that $\SRT k \leq_{\sW} \widehat{\wf}^{[k]}$ (so, also $\wfind_{\bSigma^0_k}$ and $\wfind_{\Pi^0_k}$ reduce to $\widehat{\wf}^{[k]}$). 
	\end{lemma}
	\begin{proof}
		Fix $k \geq 1$ and let $X$ be a code for a $\bSigma^0_k$ set $A \in \ramsey$, so $A$ is $\Sigma^0_k(X)$.
		Apply the function $\bmd$ $k$ times starting from $X$ to obtain a sequence of countable coded $\beta$-models $W^1 \in W^2 \in \dots \in W^k$ and a function $h$ such that $(W^k, h) \in \bmd(W^{k-1})$. Notice that $X=W^1_0$ and therefore, by Remark \ref{turingreduction}, we can compute an index for $X$ in $W^k$ (from now on we use Remark \ref{turingreduction} silently).
		
  By Lemma \ref{manybetamodels}, $W^k$ satisfies the Galvin-Prikry theorem for sets which are $\Sigma^0_k(X)$. We now repeat the construction of Lemma \ref{typeonecomputable}, introducing formulas $\mathrm{land}_k$, $\mathrm{avoid}_k$, and a function $\homo^k_h$ generalizing the formulas and the function featured there.
		Notice that again we have to appeal to Kleene's Normal Form Theorem to obtain a computable function $s$ such that the formula $\neg\mathrm{land}_k(g,x)=\neg [\forall h (g \circ h \in A_x)]$, with $A_x$ the $\bSigma^0_k$ set with name $W^k_x$, is equivalent to $\exists f \, f \oplus g \in C_{s(x),W}$.
		Similarly we obtain a computable $s'$ for the analogous formula $\mathrm{avoid}_k$.
		We then use $s$ and $s'$ in the definition of $\homo^k_h$ to obtain the desired Weihrauch reduction.
	\end{proof}
	
	We can obtain a slightly better upper bound for $\SRT k$ reasoning as in the proof of Lemma \ref{manybetamodels}.
	
	\begin{lemma}\label{inductivestep1}
		For every $k \in \omega$ we have $\SRT{k+1} \leq_{\W} \SRT k \star \widehat{\wf}$.
	\end{lemma}
	\begin{proof}
		Given $X$ coding a $\Sigma^0_{k+1}(X)$ set $A \subseteq \ramsey$, we apply $\bmd$ to $X$ to obtain a pair $(W, h) \in \bmd(X)$.
		Inspection of the proof of Lemma \ref{manybetamodels} (see \cite[Lemma VI.6.2]{simpson}) shows that via a repeated use of the function $\homo^1_h$, one can compute in $W$ an index $j$ and a sequence of indices $(m_i)_{i \in \omega}$ such that each $W_{m_i}$ is a function and, if we let $f(i)=W_{m_i}(0)$ for every $i$, then $f \oplus W_j$ computes a $\bSigma^0_k$ name for a set $B$ such that, for all $g \in \ramsey$, if $g$ lands in (resp.\ avoids) $B$, then $f \circ g$ lands in (resp.\ avoids) $A$.
		This computation of $B$ from $(W,h)$ is uniform. 
		Therefore, any $g \in \SRT k (B)$ uniformly computes a solution to $\SRT{k+1}(X)$.
	\end{proof}
    \begin{corollary}\label{inductivestep2}
		For every $k \in \omega$ we have $\wfind_{\bSigma^0_{k+1}} \leq_{\W} \wfind_{\bSigma^0_{k}} \star \widehat{\wf}$ and  $\wfind_{\bPi^0_{k+1}} \leq_{\W} \wfind_{\bPi^0_{k}} \star \widehat{\wf}$.
    \end{corollary}
    \begin{proof}
        Using the notation of the proof of Lemma \ref{inductivestep1}, if $A \in \dom(\wfind_{\bSigma^0_{k+1}})$ then $B \in \dom(\wfind_{\bSigma^0_{k}})$.
    \end{proof}
\begin{corollary}\label{slightlylower}
		For every $k \in \omega$ we have $\SRT {k+1-j} \star \widehat{\wf}^{[j]} \leq_{\W} \SRT{k-j} \star \widehat{\wf}^{[j+1]}$, so in particular $\SRT{k+1} \leq_{\W}\SRT1 \star \widehat{\wf}^{[k]}$. Similarly $\wfind_{\bSigma^0_{k+1}} \leq_{\W} \wfind_{\bSigma^0_{1}} \star \widehat{\wf}^{[k]}$ and  $\wfind_{\bPi^0_{k+1}} \leq_{\W} \wfind_{\bPi^0_{1}} \star \widehat{\wf}^{[k]}$.
	\end{corollary}
	\begin{proof}
		Follows from Lemma \ref{inductivestep1} and Corollary \ref{inductivestep2} by induction.
	\end{proof}

 	We obtain more familiar upper bounds for $\wfind_{\bPi^0_{k+1}}$ and $\wfind_{\bSigma^0_{k+1}}$ exploiting the known reductions $\wfind_{\bPi^0_{1}} \leq_{\W} \CBaire$ and $\wfind_{\bSigma^0_{1}} \equiv_{\W} \UCBaire$ (Proposition \ref{Prop:MV}).
	
	\begin{corollary}\label{morefamiliar}
		For every $k \in \omega$, we have $\wfind_{\bPi^0_{k+1}} \leq_{\W} \CBaire\star \widehat{\wf}^{[k]}$ and $\wfind_{\bSigma^0_{k+1}} \leq_{\W} \UCBaire \star \widehat{\wf}^{[k]}$.
	\end{corollary}
 
	We show in Proposition \ref{Prop:charwfind} that these reductions are actually equivalences in the context of arithmetical Weihrauch degrees.
	
	The upper bound of Corollary \ref{slightlylower} is actually an improvement over that of Lemma \ref{firstupperbound}:
	
	\begin{proposition}\label{Prop:hjstrictlyabove}
		For every $k \in \omega$ we have that $\widehat{\wf}^{[k+1]} \nleq^{\ari}_{\W} \SRT1 \star \widehat{\wf}^{[k]}$.
	\end{proposition}
	\begin{proof}
		By Lemma \ref{Lem:wfquivhj} the reduction in the statement holds if and only if the Weihrauch reduction $\hj^{k+1} \leq^{\ari}_{\W} \SRT1 \star \hj^{k}$ holds.
		Now by contradiction assume $\Phi$ and $\Psi$ witness the latter reduction.
		Let $\mathcal{M}$ be a $\beta$-model such that $\hj^k(\emptyset) \in \mathcal{M}$ but $\hj^{k+1}(\emptyset) \notin \mathcal{M}$ as in Lemma \ref{smallbetamodels}.
		Let $p=\Phi(\emptyset)$: $p$ is arithmetical, so $\hj^k(p) \equiv_{\mathrm{HYP}} \hj^k(\emptyset)$, $\hj^k(p) \in \mathcal{M}$, and also any set hyperarithmetical in $\hj^k(p)$ is in $\mathcal{M}$.
		Since $\mathcal{M} \vDash \SRT1$ and $\mathcal{M}$ is a $\beta$-model, there is some $q \in \mathcal{M} \cap (\SRT1 \star \hj^k)(p)$.
		So, by assumption, $\Psi(q \oplus \emptyset)=\hj^{k+1}(\emptyset)$, which implies $\hj^{k+1}(\emptyset) \in \mathcal{M}$.
		This is a contradiction. 
	\end{proof}
	
	To obtain further results we exploit the analysis of Section \ref{secseparation}.
	
	\begin{proposition}\label{computingmorehyperjumps}
		For any $k \geq 1$ we have $ \widehat{\wf}^{[k]} <^{\ari}_{\W} \wfind_{\bDelta^0_{k+1}}$.
	\end{proposition}
	\begin{proof}
		The reduction $ \hj^k \leq^{\ari}_{\W} \wfind_{\bDelta^0_{k+1}}$ is an immediate consequence of the results in subsection \ref{sgeneralcase}.
		Indeed, it is clear (exploiting the uniformity of the function $X \mapsto D_{k,X}$ in Lemma \ref{complexityofdn}) that there is a Turing functional $\Phi$, which, given as input $X \in 2^{\omega}$, outputs a $\bDelta^0_{k+1}$-name for $A_{k,X}$.
		For the postprocessing, the relativization of Corollary \ref{uniformreduction} gives a Turing functional $\Psi$ such that for any $X$ and any $g \in \mathrm{HS}(A_{k,X})$, $\Psi((X \oplus g)^{(k)})=\hj^k(X)$.
		
		To prove that the converse reduction does not hold, let $\Phi$ and $\Psi$ be arithmetical functions.
		We show that they do not witness $\wfind_{\bDelta^0_{k+1}} \leq^{\ari}_{\W} \hj^k$.
		Let $p$ be a computable $\bDelta^0_{k+1}$-name for the set $C_k$ of Lemma \ref{nonconverse}. Since $\Phi(p)$ is arithmetical, $\hj^k(\Phi(p)) \equiv_{\T} \hj^k(\emptyset)$.
		Consequently, $\Psi(p, \hj^k(\Phi(p))) \in \mathsf{ARITH}(\hj^k(\emptyset))$, so by Corollary \ref{arith} we have that $\Psi(p, \hj^k(\Phi(p))) \notin \mathrm{HS}(C_k)$.
	\end{proof}

	\begin{lemma}\label{Lem:absorbjumps}
		For every $a \in \mathcal{O}$, we have $\wfind_{\bPi^0_1} \star \J^{(a)} \equiv_{\W} \wfind_{\bPi^0_1}$ and $\wfind_{\bSigma^0_1} \star \J^{(a)} \equiv_{\W} \wfind_{\bSigma^0_1}$.
	\end{lemma}
	\begin{proof}
		One reduction of both equivalences is obvious, while for the other one we use the explicit definition of $\star$ given in Definition \ref{Def:operations}. In the rest of this proof we denote with $\pi_1$ and $\pi_2$ the projections $(\baire)^2 \rightarrow \baire$ (so $\pi_i(\langle p_1, p_2\rangle)=p_i$).
	
		We first prove $\wfind_{\bPi^0_1} \star \J^{(a)} \leq_{\W} \wfind_{\bPi^0_1}$. Let $(p,q) \in \dom(\wfind_{\bPi^0_1} \star \J^{(a)})$. By definition we have that $\pi_2(\Phi_p(\J^{(a)}(q))) \in \dom(\real{(\wfind_{\bPi^0_1})})$, and, we can rephrase this as $\{u\}^{p \oplus \J^{(a)}(q)} \in \dom(\real{(\wfind_{\bPi^0_1})})$ for some index $u$ independent of $p$ and $q$.
		This is equivalent to the fact that $\{u\}^{p \oplus \J^{(a)}(q)}$ is the characteristic function of a tree $T(p,q)$ with $[T(p,q)]=\delta_{\bPi^0_1}(\{u\}^{p \oplus \J^{(a)}(q)})$.
		
		Now let $n$ be an index such that, for every oracle $g \in 2^\omega$, if $g$ is the characteristic function of a tree $T$, then $[S^g_n]=[T]$. For any index $e$ let $d_e$ be an index such that $\{d_e\}^X(m) =\{n\}^{\{e\}^X}(m)$ for all $m$. We obtain that, for every $g$, if $\{e\}^g$ is the characteristic function of a tree $T$, $[S^g_{d_e}]=[S^{\{e\}^g}_n]=[T]$.
		
		By Remark \ref{anyjumps}, we can compute a $\bPi^0_1$-name for the set $K_{a,q}$ and we are guaranteed that, if $f \in \mathrm{HS}(K_{a,q})$, then $\{h(a)\}^{f \oplus q}=\J^{(a)}(q)$. So, if we let $\{e\}$ be an index such that $\{e\}^{f \oplus p \oplus q}=\{u\}^{p \oplus \{h(a)\}^{f\oplus q}}$ we obtain that, if $f \in \mathrm{HS}(K_{a,q})$, then $\{e\}^{f \oplus p \oplus q}=\{u\}^{p \oplus \J^{(a)}(q)}$ and therefore $[S^{f \oplus p \oplus q}_{d_{e}}]=[T(p,q)]$.
		
		So, we can compute a name for a $\bPi^0_1$ set $A_{p,q}=\{f \in [\omega]^\omega : f \in [S^{f \oplus p \oplus q}_{d_{e}}] \cap K_{a,q}\}$.
		Reasoning as in Proposition \ref{generalcase} we obtain that nothing avoids $A_{p,q}$, so if $f \in \mathrm{HS}(A_{p,q})$ then it lands in $A_{p,q}$. This implies that $\{h(a)\}^{f \oplus q}=\J^{(a)}(q)$, and for every $g \in [f]^{\omega}$, $g \in [S^{g \oplus p \oplus q}_{d_{e}}]=[T(p,q)]$. Therefore $f$ also lands in $\delta_{\bPi^0_1}(\{u\}^{p \oplus \J^{(a)}(q)})$.
		
		Putting everything together, we can obtain the desired Weihrauch reduction as follows: given $(p,q) \in \dom(\wfind_{\bPi^0_1} \star \J^{(a)})$, we compute a $\bPi^0_1$-name for $A_{p,q}$. Given $f \in \wfind_{\bPi^0_1}(A_{p,q})$, together with our input $(p,q)$, we output the pair $(\pi_1(\Phi_p(\{h(a)\}^{f \oplus q})), f)=(\pi_1(\Phi_p(\J^{(a)}(q))), f)$.\smallskip
		
		The reduction $\wfind_{\bSigma^0_1} \star \J^{(a)} \leq_{\W} \wfind_{\bSigma^0_1}$ is proved similarly, so we only sketch the proof. 
		In analogy with the $\bPi^0_1$ case, we know that for every $(p,q) \in \dom(\wfind_{\bSigma^0_1} \star \J^{(a)})$, $\{u\}^{p \oplus \J^{(a)}(q)}$ is a $\bSigma^0_1$-name for an open set $A$, i.e.\ $\{u\}^{p \oplus \J^{(a)}(q)}$ is an element of $\baire$ which lists codes for basic clopen subsets of $\ramsey$ whose union is $A$.
		
		Now let $i \in \omega$ be a code for the empty set.
		It is easy to see that there is a computable function $u'$ such that, $\{u'(e)\}^{p'} \in \baire$ for every oracle $p'$ and index $e$ and, if $\{e\}^{p'} \in \baire$ and $\{u\}^{\{e\}^{p'}} \in \baire$, then $\ran(\{u'(e)\}^{p'})=\ran(\{u\}^{\{e\}^{p'}}) \cup \{i\}$ (intuitively $\{u'(e)\}^p$ attempts to compute $\{u\}^{\{e\}^{p}}(m)$ in parallel for every $m$ and periodically outputs the code $i$, so it is guaranteed to always produce an element of $\baire$).
		Therefore we can stipulate that the preprocessing function computes a $\bSigma^0_1$-name for the set $\delta_{\bSigma^0_1}(\{u'\}^{p \oplus {h(a)}^{f \oplus q}}) \cap K_{a,q}$. The rest of the proof is the same as the proof of the $\bPi^0_1$ case.
	\end{proof}
	
	\begin{remark}\label{Rem:aristar}
 		Let $f,f'$ be such that $f \star \J^{(\omega)} \leq_{\W} f$, $f' \star \J^{(\omega)} \leq_{\W} f'$ and $f \equiv^{\ari}_{\W} f'$ (for two examples, by Corollary \ref{Cor:choiceabsorbjumps} and Lemma \ref{Lem:absorbjumps}, one can pick as $(f,f')$ the pair ($\wfind_{\bPi^0_1}$, $\CBaire$) or the pair ($\wfind_{\bSigma^0_1}$, $\UCBaire$)). Then for any $g, g'$ such that $g \leq^{\ari}_{\W} g'$, we have $f \star g \leq^{\ari}_{\W} f' \star g'$.
 		This is an immediate consequence of the definition of compositional product together with the characterization $f \leq^{\ari}_{\W} g \iff \exists n \,\, f \leq_{\W} \lim^{[n]} \star g \star \lim^{[n]}$ (see \cite[pg.\ 344]{marconevalenti}).
	\end{remark}
	
	\begin{proposition}\label{Prop:charwfind}
		For every $k \in \omega$, we have $\wfind_{\bPi^0_{k+1}} \equiv^{\ari}_{\W} \CBaire \star \widehat{\wf}^{[k]}$ and  $\wfind_{\bSigma^0_{k+1}} \equiv^{\ari}_{\W} \UCBaire \star \widehat{\wf}^{[k]}$.
	\end{proposition}
	\begin{proof}
		
		The case $k=0$ is Proposition \ref{Prop:MV}, so we can assume $k > 0$.
		
		In light of Corollary \ref{morefamiliar}, to show $\wfind_{\bPi^0_{k+1}} \equiv^{\ari}_{\W} \CBaire \star \widehat{\wf}^{[k]}$ we only need to show $\CBaire \star \hj^k \leq^{\ari}_{\W} \wfind_{\bPi^0_{k+1}}$.
		Since $\CBaire \leq^{\ari}_{\W} \wfind_{\bPi^0_1}$ by Remark \ref{Rem:aristar} it suffices to show $\wfind_{\bPi^0_1} \star \hj^k \leq^{\ari}_{\W} \wfind_{\bPi^0_{k+1}}$.
		
		Notice that if $(p,q) \in \dom(\wfind_{\bPi^0_1} \star \hj^k)$, then $\pi_2(\Phi_p(\hj^k(q)))$ is the characteristic function of a tree $T_{p,q}$ such that $\mathrm{HS}([T_{p,q}]) \subseteq [T_{p,q}]$.
        Let $e$ be an index such that $\pi_2(\Phi_x(y))=\{e\}^{x \oplus y}$ for all $x, y \in \baire$, let $n$ be an index such that for every $A \in 2^{\omega}$, if $A$ is the characteristic function of a tree $T_A$, $[S^A_n]=[T_A]$, and let $r$ be an index such that $\{r\}^X(m)=\{n\}^{\{e\}^X}(m)$ for all $m$. 
		Therefore we obtain that if $(p,q) \in \dom(\wfind_{\bPi^0_1} \star \hj^k)$, then $[S^{p \oplus (\hj^k(q))}_r]=[S^{\{e\}^{p \oplus \hj^k(q)}}_n]=[S^{\pi_2(\Phi_p(\hj^k(q)))}_n]=[T_{p,q}]$.
		
		Now consider the formula
		\[
		\mu(f,g,h)=\chi^h_{k}(f) \land  f \in [S^{g \oplus D_{k,h}(f)}_r],
		\]
        where $\chi^h_{k}$ is the formula obtained by relativizing the formula in Definition \ref{chin} to the oracle $h$ as we did in Proposition \ref{computingmorehyperjumps}.
		The formula $\mu$ is $\Pi^0_{k+1}$ because the first conjunct is $\Delta^0_{k+1}$ by (the relativization of) Proposition \ref{generalcase}, while the second conjunct is $\Pi^0_1$ in $f^{(k)}, g, h$ and hence $\Pi^0_{k+1}$ in $f,g,h$.
		By arguments similar to those used in the proof of Proposition \ref{generalcase}, for every $(p,q) \in \dom(\wfind_{\bPi^0_1} \star \hj^k)$ nothing avoids the set $C_{p,q}=\{f \in \ramsey : \mu(f,p,q)\}$. For any $f \in \mathrm{HS}(C_{p,q})$, we have $D_{k,q}(f)=\hj^k(q)$, so $[S_r^{p \oplus D_{k,q}(f)}]=[S_r^{p \oplus \hj^k(q)}]=[T_{p,q}]$ and therefore $f$ lands in the closed set $[T_{p,q}]$.
		
		So we get our arithmetical Weihrauch reduction as follows: the preprocessing takes an input $(p,q) \in \dom(\wfind_{\bPi^0_1} \star \hj^k)$ and produces (computably) a $\bPi^0_{k+1}$-name for $C_{p, q}$. The postprocessing takes $f \in \mathrm{HS}(C_{p,q})$ and $(p, q)$ and outputs the pair $(\pi_1(\Phi_p(D_{k,q}(f))), f)=(\pi_1(\Phi_p(\hj^k(q)), f)$.\smallskip
		
		The proof of $\wfind_{\bSigma^0_{k+1}} \equiv^{\ari}_{\W} \UCBaire \star \widehat{\wf}^{[k]}$ is obtained by modifying the proof above exactly as in Lemma \ref{Lem:absorbjumps}.
	\end{proof}
	
	\begin{corollary}\label{arisepwfind}
		For all $k \in \omega$, $\wfind_{\bSigma^0_{k+1}} <^{\ari}_{\W} \wfind_{\bPi^0_{k+1}}$.
	\end{corollary}
	\begin{proof}
		By Proposition \ref{Prop:charwfind} we only need to show that, for every $k \in \omega$, $\UCBaire \star \hj^k <^{\ari}_{\W} \CBaire \star \hj^k$.
		The reduction $\UCBaire \star \hj^k \leq^{\ari}_{\W} \CBaire \star \hj^k$ holds trivially as the left side is a restriction of the right side to a smaller domain.
		For the separation, note that (by the relativization of a known fact for the functions $\UCBaire$ and $\CBaire$) for all arithmetical $(p,q) \in \dom(\UCBaire \star \hj^k)$, any $y \in (\UCBaire \star \hj^k)(p,q)$ is hyperarithmetical in $\hj^k(\emptyset)$, while there are computable instances of $(\CBaire \star \hj^k)$ which do not admit solutions hyperarithmetical in $\hj^k(\emptyset)$.
		This establishes  $\CBaire \star \hj^k \nleq^{\ari}_{\W} \UCBaire \star \hj^k$.
	\end{proof}
	
	The proofs of Lemma \ref{Lem:absorbjumps} and Proposition \ref{Prop:charwfind} essentially exploit the fact that the functions $\wfind_{\Gamma}$ are one-sided. This feature can be used to obtain analogous lower bounds for functions of the form $\find_{\Gamma}$.
	
	\begin{proposition}\label{Prop:lowerboundfind}
		For every $k \in \omega$, we have $\find_{\bPi^0_1} \star \widehat{\wf}^{[k]}\leq^{\ari}_{\W} \find_{\bPi^0_{k+1}}$ and  $\find_{\bSigma^0_1} \star \widehat{\wf}^{[k]}\leq^{\ari}_{\W} \find_{\bSigma^0_{k+1}}$.	
	\end{proposition}	
	\begin{proof}
	    Note that the proof of Proposition \ref{Prop:charwfind} actually shows $\wfind_{\bPi^0_{1}} \star \hj^k \leq^{\ari}_{\W} \wfind_{\bPi^0_{k+1}}$ and $\wfind_{\bSigma^0_{1}} \star \hj^k \leq^{\ari}_{\W} \wfind_{\bSigma^0_{k+1}}$.
	    We claim that the same reductions also witness $\find_{\bPi^0_{1}} \star \hj^k \leq^{\ari}_{\W} \find_{\bPi^0_{k+1}}$ and $\find_{\bSigma^0_{1}} \star \hj^k \leq^{\ari}_{\W} \find_{\bSigma^0_{k+1}}$.
		Indeed, let $G$ be the preprocessing function witnessing the reduction $\wfind_{\bPi^0_{1}} \star \hj^k \leq^{\ari}_{\W} \wfind_{\bPi^0_{k+1}}$ described in Proposition \ref{Prop:charwfind} and let $(p,q) \in \dom(\find_{\bPi^0_1} \star \hj^k)$. Then $G(p,q)$ is a $\bPi^0_{k+1}$ name for the set
		\[
		C_{p,q}=\{f \in \ramsey : \chi^q_k(f) \land f \in [S^{g \oplus D_{k,q}(f)}_r]\},
		\]
		and, letting $B_{p,q}=\delta_{\bPi^0_1}(\pi_2(\Phi_p(\hj^k(q))))$, by the proof of Proposition \ref{Prop:charwfind} we have that if $x \in \ramsey$ lands in the set $A_{k,q}=\{f \in \ramsey : \chi^q_k(f)\}$, then $x \in C_{p,q}$ if and only if $x \in B_{p,q}$.
		
		Now by assumption we have that $\mathrm{HS}(B_{p,q}) \cap B_{p,q} \neq \emptyset$, so we can pick $g \in \mathrm{HS}(B_{p,q}) \cap B_{p,q}$ and, applying Corollary \ref{relativization}, we can pick some $h \in [g]^{\omega}$ which is homogeneous for $A_{k,q}=\{f \in \ramsey : \chi^q_k(f)\}$.
		By Proposition \ref{generalcase}, we have that $h$ lands in $A_{k,q}$, hence, for all $x \in [h]^{\omega}$, $x \in C_{p,q}$ if and only if $x \in [S^{g \oplus D_{k,q}(f)}_r]$ if and only if $x \in B_{p,q}$. 
		Since $h \in [g]^{\omega}$ and $g$ lands in $B_{p,q}$, it follows that for every $x \in [h]^{\omega}$, $x \in B_{p,q}$.
		Putting everything together, we obtain that $x$ lands in $C_{p,q}$, so $\delta_{\bPi^0_{k+1}}(G(p,q))=C_{p,q} \in \dom(\find_{\bPi^0_{k+1}})$.
		
		The proof of Proposition \ref{Prop:charwfind} shows that, conversely, any $h \in \mathrm{HS}(C_{p,q}) \cap C_{p,q}$ lands in $A_{k,q}$ as well as in $B_{p,q}$. This shows that the functions witnessing  $\wfind_{\bPi^0_{1}} \star \hj^k \leq^{\ari}_{\W} \wfind_{\bPi^0_{k+1}}$ in the above proof also witness $\find_{\bPi^0_{1}} \star \hj^k \leq^{\ari}_{\W} \find_{\bPi^0_{k+1}}$. The same reasoning applies to the reduction $\find_{\bSigma^0_{1}} \star \hj^k \leq^{\ari}_{\W} \find_{\bSigma^0_{k+1}}$.
		\end{proof} 
	
	\begin{corollary}\label{Cor:findaboveall}
		For any $k \in \omega$, we have $\SRT {k+1} \leq^{\ari}_{\W}\find_{\bSigma^0_{k+1}}$.
	\end{corollary}
	\begin{proof}
		We have 
  $$\SRT {k+1} \leq_{\W} \SRT 1 \star \widehat{\wf}^{[k]} \leq_{\W} \find_{\bSigma^0_{1}} \star \widehat{\wf}^{[k]} \leq^{\ari}_{\W} \find_{\bSigma^0_{k+1}}$$
  where we use Lemma \ref{slightlylower}, Proposition \ref{Prop:findabovesrt} and Proposition \ref{Prop:lowerboundfind}.
	\end{proof}
	
	\begin{lemma}\label{Lem:MV415}
		For any $k \in \omega$, $\SRT 1 \star \widehat{\wf}^{[k]} \nleq^{\ari}_{\W} \CBaire \star \widehat{\wf}^{[k]}$.
	\end{lemma}
	\begin{proof}
		It follows from \cite[Proposition 4.14]{marconevalenti} that $\wf \leq_{\W} \J^1 \star \SRT 1$, so by monotonicity of $\star$ we have $\wf \star \hj^k \leq_{\W} \J^1 \star \SRT 1 \star \hj^k$ and hence $\wf \star \hj^k \leq^{\ari}_{\W} \SRT 1 \star \hj^k$.
		Therefore, as in the proof of \cite[Corollary 4.15]{marconevalenti}, we have that if $\SRT 1 \star \hj^k \leq^{\ari}_{\W} \CBaire \star \hj^k$, then $\wf \star \hj^k \leq^{\ari}_{\W} \CBaire \star \hj^k$.
		This in turn would imply that $\widehat{(\wf \star \hj^k)} \leq^{\ari}_{\W} \widehat{(\CBaire \star \hj^k)}$.
		
		Now, it is immediate to see that $\hj^{k+1} \equiv_{\W} \widehat{\wf} \star \hj^k \leq_{\W}\widehat{(\wf \star \hj^k)}$.
		On the other hand, by Lemma \ref{Lem:distributivity} and since $\CBaire$ and $\hj^k$ are both parallelizable, we have $\widehat{(\CBaire \star \hj^k)} \leq_{\W} \widehat{\CBaire} \star \widehat{\hj^k} \equiv_{\W} \CBaire \star \hj^k$. Putting everything together we arrive at $\hj^{k+1} \leq^{\ari}_{\W} \CBaire \star \hj^k$, but this reduction cannot hold as, by (the relativization to $\hj^k(\emptyset)$ of) Gandy's Basis Theorem \cite[Theorem III.1.4]{sacks2017} for any computable $x \in \dom(\CBaire \star \hj^k)$ there exists some $y \in (\CBaire \star \hj^k)(x)$ such that $\hj^{k+1}(\emptyset)$ is not hyperarithmetical in $y$.
	\end{proof}
	
	\begin{corollary}\label{Cor:findstrictabove}
		For any $k \in \omega$, $\find_{\bSigma^0_{k+1}} \nleq^{\ari}_{\W} \SRT 1 \star \widehat{\wf}^{[k]}$.
	\end{corollary}
	\begin{proof}
		Using Proposition \ref{Prop:findabovesrt} and reasoning as in  Corollary \ref{Cor:findaboveall} we obtain $(\CBaire \times \SRT 1) \star \widehat{\wf}^{[k]} \leq^{\ari}_{\W} \find_{\bSigma^0_{k+1}}$, so it suffices to show $(\CBaire \times \SRT 1) \star \hj^k \nleq^{\ari}_{\W} \SRT 1 \star \hj^k$. 
  We claim that $(\CBaire \times \SRT 1) \star \hj^k \leq^{\ari}_{\W} \SRT 1 \star \hj^k$ implies $\SRT 1 \star \hj^k \leq^{\ari}_{\W} \CBaire \star \hj^k$, contradicting Lemma \ref{Lem:MV415}.
  In the proof of \cite[Theorem 5.14]{marconevalenti} the claim is established for $k=0$, and the same argument works also when $k>0$.
	\end{proof}
	
	We remark that Proposition \ref{computingmorehyperjumps} and Lemma \ref{firstupperbound} together imply that the principles $\SRT k$ for $k \in \omega$ form a strictly increasing chain in the arithmetical Weihrauch degrees. Concentrating on the different functions related to a single Borel rank $k+1$, we obtain the picture of Figure \ref{figure1}.
 
 \medskip
	
	One may wonder whether the arithmetical Weihrauch reduction of Proposition \ref{computingmorehyperjumps} (and other reductions of the form $\widehat{\wf}^{[k]} \leq^{\ari}_{\W} f$ with $f \leq_{\W} \SRT {\alpha}$) can be improved to Weihrauch reductions.
	We show that this is impossible for very general reasons.
	\begin{definition}
		A set $A \subseteq \omega$ is called \emph{recursively encodable} if for every $B \in \ramsey$ there exists $C \in [B]^{\omega}$ such that $A \leq_{\T}C$.
	\end{definition}
	
	\begin{lemma}[{Jockush-Solovay, \cite{Jockush}, \cite{hypencsets}}]\label{recursivelyencodable}
		A set $A \in \ramsey$ is recursively encodable iff it is hyperarithmetic.
	\end{lemma}
	
	\begin{definition}\label{def:mes}
		We say that a function $f \subc \baire \rightarrow \baire$ \emph{has non-hyperarithmetic outputs} if there is some computable $x \in \dom(f)$ such that $f(x)$ is not hyperarithmetic.
		We say that a partial multi-valued $f \subc \baire \multif \baire$ \emph{meets every infinite set} if, for every $p \in \dom(f)$ and every $g \in \ramsey$, $[g]^{\omega} \cap f(p) \neq \emptyset$.
	\end{definition}
	
	Using this terminology, the function $\widehat{\wf}$ has non-hyperarithmetic outputs and, for every $\alpha \in \omega_1$, $\real{(\SRT {\alpha})}$ meets every infinite set (this is basically Corollary \ref{relativization}).
	
	\begin{proposition}\label{prop:nononhyp}
		If $g \subc {\omega}^{\omega} \multif {\omega}^{\omega}$ meets every infinite set and $f \subc {\omega}^{\omega} \rightarrow {\omega}^{\omega}$ has non-hyperarithmetic outputs, then $f \nleq_{\W} g$.
	\end{proposition}
	\begin{proof}
		Suppose $\Phi$ and $\Psi$ witness $f \leq_{\W} g$ and let $x$ be a computable input for $f$ such that $f(x)$ is not hyperarithmetic.
		Then $\Phi(x) \in \dom(g)$ and for all $y \in g(\Phi(x))$, $\Psi(x,y)=f(x)$.
		Since $x$ is computable, this in particular implies that $f(x) \leq_{\T} y$ for all $y \in g(\Phi(x))$.
		Now $f(x)$ is not hyperarithmetic, hence by Lemma \ref{recursivelyencodable} there exists some $h \in \ramsey$ such that $f(x) \not \leq_{\T} k$ for all $k \in [h]^\omega$.
		Since $g$ meets every infinite subset, there is some $z \in g(\Phi(x)) \cap [h]^{\omega}$.
		Now $z \in \Phi(g(x))$ implies $f(x) \leq_{\T} z$, but by definition of $h$ we have $f(x) \not\leq_{\T}z$.
	\end{proof}
	
	\begin{corollary}\label{nohyperjumps}
		$\widehat{\wf} \nleq_{\W} \SRT {\alpha}$ for all $\alpha \in \omega_1$.
	\end{corollary}
	
	\begin{corollary}\label{Cor:strict}
		For every $k \geq 1$ we have $\SRT {k+1} <_{\W} \SRT 1 \star \widehat{\wf}^{[k]}$.
	\end{corollary}
	\begin{proof}
		The reduction holds by Corollary \ref{slightlylower}, and it is strict by Corollary \ref{nohyperjumps}.
	\end{proof}
	
	By Corollary \ref{nohyperjumps} we cannot strengthen Proposition \ref{Prop:charwfind} substituting arithmetical Weihrauch reducibility with Weihrauch reducibility.
    Moreover, together with the separation of Proposition \ref{computingmorehyperjumps}, Corollary \ref{nohyperjumps} also implies that $\SRT k$ and $\widehat{\wf}^{[k]}$ are incomparable for every $k$.

    Despite these negative results, a closer look at the proof of Proposition \ref{Prop:charwfind} shows that we can prove more about the Weihrauch degrees of one-sided principles.
	
    \begin{proposition}\label{Prop:wstarWreductions}
    	For every $k \in \omega$, the following hold:
    	\begin{enumerate}
    		\item $\wfind_{\bPi^0_{1}} \wstar \widehat{\wf}^{[k]} \leq_{\W} \wfind_{\bPi^0_{k+1}}$, $\wfind_{\bSigma^0_{1}} \wstar \widehat{\wf}^{[k]} \leq_{\W} \wfind_{\bSigma^0_{k+1}} $,
    		\item $\find_{\bPi^0_{1}} \wstar \widehat{\wf}^{[k]} \leq_{\W} \find_{\bPi^0_{k+1}}$, $\find_{\bSigma^0_{1}} \wstar \widehat{\wf}^{[k]} \leq_{\W} \find_{\bSigma^0_{k+1}}$.
    	\end{enumerate}
    \end{proposition}
    \begin{proof}
    	For (1) the proof of Proposition \ref{Prop:charwfind} shows the reduction $\wfind_{\bPi^0_{1}} \star \widehat{\wf}^{[k]} \leq^{\ari}_{\W} \wfind_{\bPi^0_{k+1}}$ via a computable preprocessing, call it $F$, and a non-computable arithmetical postprocessing, call it $G$. Looking at this specific $G$ and letting $\pi_2 \colon (\baire)^2 \rightarrow \baire$ denote the projection on the second component, we see that $\pi_2  \circ G=\id {{\baire}}$. Therefore, if we let $K$ be the computable functional of Lemma \ref{lem:relwstarstar}, we obtain that $F \circ K$ and $\id{{\baire}}=\pi_2 \circ G$ witness $\wfind_{\bPi^0_{1}} \wstar \widehat{\wf}^{[k]} \leq_{\W} \wfind_{\bPi^0_{k+1}}$. The same applies to the case of $\wfind_{\bSigma^0_{k+1}}$.
    	
    	The proof of (2) is the same as the proof of (1).
    \end{proof}
    
    \begin{corollary}\label{Cor:findWLB}
    	We have $\find_{\bPi^0_{1}} \star \widehat{\wf}^{[k]} \leq_{\W} \find_{\bPi^0_{k+1}}$ and $\find_{\bSigma^0_{1}} \star \widehat{\wf}^{[k]} \leq_{\W} \find_{\bSigma^0_{k+1}}$.
    \end{corollary}
    \begin{proof}
    	Since $\find_{\bPi^0_{1}}$ and $\find_{\bSigma^0_{1}}$ are cylinders (\cite[Theorem 4.7, Corollary 4.25]{marconevalenti}), using Remark \ref{rem:wstarequalstar} together with Proposition \ref{Prop:wstarWreductions} we obtain 
        \[
    	\find_{\bPi^0_{1}} \star \widehat{\wf}^{[k]} \equiv_{\W} \find_{\bPi^0_{1}} \wstar \widehat{\wf}^{[k]} \leq_{\W} \find_{\bPi^0_{k+1}}
    	\] 
    	and 
     \[
    	\find_{\bSigma^0_{1}} \star \widehat{\wf}^{[k]} \equiv_{\W} \find_{\bSigma^0_{1}} \wstar \widehat{\wf}^{[k]} \leq_{\W} \find_{\bSigma^0_{k+1}}.\qedhere
    	\]
    	
    \end{proof}
	
	\begin{corollary}\label{Cor:upleftarrow}
		We have $\SRT {1} \star \widehat{\wf}^{[k]} <_{\W} \find_{\bSigma^0_{k+1}}$.
	\end{corollary}
	\begin{proof}
		By the proof of Corollary \ref{Cor:findaboveall} we have $\SRT {1} \star \widehat{\wf}^{[k]} \leq_{\W} \find_{\bSigma^0_{1}} \star \widehat{\wf}^{[k]}$, so $\SRT {1} \star \widehat{\wf}^{[k]} \leq_{\W} \find_{\bSigma^0_{k+1}}$ follows from Corollary \ref{Cor:findWLB}. The strictness of the reduction follows from Corollary \ref{Cor:findstrictabove}.
	\end{proof}
	\begin{corollary}\label{Cor:strongsep}
		We have $\wfind_{\bPi^0_{k+1}} <_{\W} \find_{\bPi^0_{k+1}}$.
	\end{corollary}
	\begin{proof}
		The reduction is immediate by definition. The strictness follows from the fact that, by Corollary \ref{Cor:findWLB}, $\widehat{\wf}^{[k]} \leq_{\W} \find_{\bPi^0_{k+1}}$, while $\widehat{\wf}^{[k]} \nleq_{\W} \wfind_{\bPi^0_{k+1}}$ by Corollary \ref{nohyperjumps}. 
	\end{proof}
    
    We conclude with a treatment of the uniform relationship between the Galvin-Prikry theorem for Borel sets of transfinite rank and transfinite iterations of the hyperjump. 
    
    \begin{definition}
    	For any $\alpha < \omega_1$, define $\hj^{\alpha} \subc 2^{\omega} \rightarrow 2^{\omega}$, with domain given by $\{X \in 2^{\omega}: \omega^X_1 > \alpha\}$ as $\hj^{\alpha}(X)=\{(a,\hj^a(X)) : a \in \mathcal{O}_X \land |a|_X=\alpha\}$.
    \end{definition}
   
	With this definition, translating Proposition \ref{Prop:transfiniteub} in the language of Weihrauch reductions we obtain:
	
	\begin{proposition}
		For any $\alpha < \omega_1$, we have $\SRT{\omega+\alpha}  \leq_{\W} \SRT 1 \star  \hj^{\omega^{1+\alpha}}$.
	\end{proposition}
	
	Going in the other direction, one could wonder whether Proposition \ref{transfinitecase} also corresponds to a Weihrauch reduction.
	Indeed one would expect that, with computable preprocessing and hyperarithmetical postprocessing, one should be able to obtain $\hj^{\alpha}$ from $\SRT{1+\alpha}$, for every $\alpha< \omega_1$.
	
	One possible way to approach the question is the following.
	We define the relativization of the sets $C_{a}$ in Definition \ref{def:transfinitelevels} to an oracle $X$ in the natural way to obtain, for any given $X$, a sequence of sets $(C_{a,X})_{a \in \mathcal{O}_X}$.
	It is easy to see that there is a computable functional $\Phi'$ which, given $a$ and $X$, outputs a Borel code for $C_{a,X}$.
	
	In order to obtain a reduction from $\hj^{\alpha}$ to $\SRT{1+\alpha}$ using Proposition \ref{transfinitecase} we need a computable functional $\Phi$ which, given any $X$ with $\alpha < \omega^X_1$, outputs a Borel code for a set of the form $C_{a,X}$ where $a \in \mathcal{O}_X$ and $|a|_X=\alpha$. 
	Clearly, if $\alpha <\omega_1^{\ck}$, then $\Phi$ exists and can be easily defined from $\Phi'$. This yields:
	\begin{proposition}
		Let $\alpha < \omega^{\ck}_1$. There are a computable functional $\Phi$ and an hyperarithmetical functional $\Psi$ such that, given any $f \vdash \SRT {1+\alpha}$, the function $x \mapsto \Psi(x, f(\Phi(x)))$ is a realizer of $\hj^\alpha$.
	\end{proposition}
	In case $\alpha \geq \omega^{\ck}_1$ we do not know of a way to obtain a functional $\Phi$ as above.
	
	\section{Reverse mathematics of the lightface Galvin-Prikry theorem}\label{revmath}
	
	We use the results of the previous section to obtain reverse mathematical corollaries relating fragments of the lightface Galvin-Prikry theorem to iterated hyperjumps.
	
	\begin{theorem}\label{equivalences}
		For any $k \in \omega$, the following are equivalent for $\beta$-models:
		\begin{enumerate}
			\item $\hj^k(\emptyset)$ exists,
			\item there exists a sequence of countable coded $\beta$-models $M^1 \in M^2 \in \dots \in M^k$,
			\item the lightface $\Sigma^0_{k+1}$ Galvin-Prikry theorem,
			\item the lightface $\Delta^0_{k+1}$ Galvin-Prikry theorem.
		\end{enumerate}
	\end{theorem}
	\begin{proof}
		$1 \rightarrow 2$ follows from Lemma \ref{smallbetamodels}, $2 \rightarrow 3$ is Corollary \ref{simpsonreversal}, $3 \rightarrow 4$ is trivial.
		The implication $4 \rightarrow 1$ is Theorem \ref{reversal}.
	\end{proof}
	
	\begin{corollary}\label{betamodelseparation}
		For every $n \in \omega$ there exists a $\beta$-model $M^n$ such that $M^n \vDash \sRT n$ but $M^n \nvDash \sRT{n+1}$. 
	\end{corollary}
	\begin{proof}
		By Theorem \ref{equivalences}, it suffices to show that for every $n \in \omega$ there are $\beta$-models $M^n$ which contain $\hj^n(\emptyset)$ but do not contain $\hj^{n+1}(\emptyset)$.
		These are the models of Lemma \ref{smallbetamodels}. 
	\end{proof}
	
	Notice that the implications $1 \rightarrow 2 \rightarrow 3 \rightarrow 4$ and $2 \rightarrow 1$ in Theorem \ref{equivalences} hold for any model of $\atr$.
	It is natural to ask whether the same can be said for $4 \rightarrow 1$.
	
	\begin{openquestion}\label{bigopenq}
		Does the lightface $\Delta^0_{k+1}$ Galvin-Prikry theorem imply the existence of $\hj^k(\emptyset)$ over $\atr$? 
	\end{openquestion}
	
	For a fixed $k \in \omega$, the proof of implication $4 \rightarrow 1$ in Theorem \ref{equivalences} goes through in the theory $\atr + \forall g \exists h \chi_k(g \circ h)$ (where $\chi_k$ is the formula of Definition \ref{chin}), as this theory suffices to prove Proposition \ref{generalcase}.
	Indeed, for all $n \in \omega$, $\aca$ proves that, if $A_n=\{f \in \ramsey : \chi_n(f)\}$ and $g \in \mathrm{HS}(A_n) \cap A_n$, then $D_n(g)=\hj^n(\emptyset)$ (this is essentially an iteration of Lemma \ref{ebar}).
	Hence:
 
 \begin{corollary}\label{Cor:T}
     Statements $1$ to $4$ in Theorem \ref{equivalences} are equivalent over the theory $T= \atr + \{\forall g \exists h \, \chi_k(g \circ h) : k \in \omega \}$.
 \end{corollary}
	Since $T$ holds in any $\beta$-model, it follows that $T$ is strictly weaker than $\pioo$.
	Thus we could ask if $T$ is also strictly stronger than $\atr$.
	A negative answer to this question would yield a positive answer to Question \ref{bigopenq}.
	
	A partial answer to Question \ref{bigopenq} was given by Suzuki and Yokoyama in \cite{suzukiyokoyama}, where it was shown that an analogue of Proposition \ref{generalcase} (their Lemma 6.13) can be proven in (less than) $\atr$, yielding:
	
	\begin{lemma}[{\cite[Lemma 6.14]{suzukiyokoyama}}]
		For every $k \in \omega$ there exists $m_k$ such that for all $X$, $\Delta^0_{m_k}(X)$-$\rt$ implies the existence of $\hj^k(X)$ over $\atr$.
	\end{lemma}
	
	Note that from the result in \cite{suzukiyokoyama} we obtain that $m_k>k+1$ for all $k \geq 2$, so it still does not provide a full answer to our question.
	
	We mention another possible strategy to tackle Question \ref{bigopenq} which the reader might deem more natural or interesting.
	
	Our proof of Proposition \ref{generalcase} (and thus, of implication $4 \rightarrow 1$ in Theorem \ref{equivalences}) relies on Corollary \ref{relativization}. 
	The proof of Corollary \ref{relativization} in turn is based on the fact that, for any $f \in \ramsey$ if $B \subseteq \ramsey$ is $\Delta^0_n$ for some $n \in \omega$, and $F \colon \ramsey \rightarrow [f]^\omega$ is given by $g \mapsto f \circ g$, then $B_f=F^{-1}[B \cap [f]^\omega]$ is $\Delta^0_n(f)$, and therefore $\DRT n$ guarantees that $\mathrm{HS}(B_f) \neq \emptyset$.
	
	Notice that, since the function $F$ is a computable homeomorphism between $\ramsey$ and $[f]^{\omega}$ which is very well behaved with respect to subsets, the set $B_f$ is a essentially a copy of $B \cap [f]^\omega$ in the ``renamed'' Ramsey space $[f]^{\omega}$. 
    Therefore, it is sensible to expect that a model of (say) $\aca$ which contains the set $f$ should perceive this similarity and consequently ``believe'' that, if $B$ admits homogeneous sets, then so does $B_f$.
	
	We could then formulate, for any class of formulas $\Gamma$, axiom schemata $(\ra {\Gamma})$ ($\mathsf{RA}$ for \emph{Ramsey absoluteness}) such that $M \vDash \ra {\Gamma}$ if and only if
	\[
	M \vDash \mathrm{HS}(B) \neq \emptyset \rightarrow \forall f \in \ramsey \,\, (\mathrm{HS}(B_f) \neq \emptyset)
	\]  
	where $B$ ranges over all $\Gamma$-definable sets of the model $M$\footnote{Formally $\ra {\Gamma}$ is the set of sentences of the form 
	\[
	(\exists h \forall g \varphi(h \circ g) \lor \exists h \forall g \neg \varphi(h \circ g)) \rightarrow \forall f (\exists h \forall g \psi_{\varphi}(h \circ g, f) \lor \exists h \forall g \neg \psi_{\varphi}(h \circ g, f))
	\]
	where $\varphi$ ranges over $\Gamma$ formulas of second order arithmetic and $\psi_{\varphi}(x,f)$ is the formula $\varphi(f \circ x)$. We point out that although $\Delta^0_n$ is not strictly speaking a class of formulas, this can be dealt with by modifying the sentence above in the usual way, so we can also define $\ra {\Delta^0_n}$. Notice that, e.g., $\sRT{n} + \ra {\Sigma^0_n}$ corresponds to the restriction of Corollary \ref{relativization} to $\Sigma^0_n$ sets. 
	}.
	
	It follows from the discussion above that for a given $n \in \omega$, the axiom schema $\sRT{n} + \ra {\Sigma^0_n}$ (resp.\ $\dRT n + \ra {\Delta^0_n}$) could be a more robust way to express that $\Sigma^0_n$ (resp.\ $\Delta^0_n$) sets are Ramsey.
	
	At this point it is natural to wonder about the proof theoretic strength of $\ra {\Sigma^0_n}$ and $\ra {\Delta^0_n}$:
	
	\begin{openquestion}\label{q:ra1}
 		Is it the case that $\aca$ proves $\bigcup_{n \in \omega} \ra {\Sigma^0_n}$? 
	\end{openquestion}
	
	\begin{openquestion}\label{q:ra2}
		Is it the case that $\aca + \sRT n$ (resp.\ $\aca + \dRT n$) proves $\ra {\Sigma^0_n}$ (resp.\ $\ra {\Delta^0_n}$)?
	\end{openquestion}
	
	Clearly, a positive answer to Question \ref{q:ra1} implies a positive answer to Question \ref{q:ra2} which in turn trivially implies that $\sRT {n}$ is equivalent to $\sRT {n} + \ra {\Sigma^0_n}$ over $\atr$. This would in some sense confirm that the schema $\sRT n + \ra {\Sigma^0_n}$ is the appropriate way to think about the statement ``$\Sigma^0_n$ sets are Ramsey''.
	Perhaps more importantly, it would also yield a positive answer to Question \ref{bigopenq}, and therefore it would allow us to weaken the theory $T$ in the hypothesis of Corollary \ref{Cor:T} by substituting it with the more natural $\atr$.
	
	On the other hand, a countermodel providing a negative answer to either Question \ref{q:ra1} or Question \ref{q:ra2} would provide an interesting and arguably pathological behaviour of the Ramsey property.
	 
	Lastly, recall that in Lemma \ref{nonconverse} and Corollary \ref{arith} we showed that, for every $n \in \omega$, there exists a $\Delta^0_{n+1}$ set $C_n$ such that $\mathrm{HS}(C_n) \cap \mathsf{ARITH}(\hj^n(\emptyset))= \emptyset$. 
	Inspection of the proof of Lemma \ref{nonconverse} reveals that it can be proven in $\aca + \ra {\Delta^0_{n+1}}$.
	Now, for every $n \in \omega$, $\mathsf{ARITH}(\hj^n(\emptyset)) \vDash \aca$.
	Hence the model $\mathsf{ARITH}(\hj^n(\emptyset))$ shows that $\aca$+``$\hj^n(\emptyset)$  exists'' $\nvdash \dRT{n+1} + \ra {\Delta^0_{n+1}}$.
	
	\begin{openquestion}\label{q:arith}
	 	Is it the case that $\aca$+``$\hj^n(\emptyset)$  exists'' $\nvdash \dRT{n+1}$? Can this be witnessed by an $\omega$-model?
	\end{openquestion}
	
	A positive answer to Question \ref{q:ra2} implies a positive answer to Question \ref{q:arith}.
	\bibliographystyle{alpha}
	\bibliography{REFERENCES}
	
\end{document}